\documentclass[12pt]{amsart}         
\usepackage{amscd}                   
\usepackage{amssymb}   
\usepackage{eucal}     

\newtheorem{thm}{Theorem}[section]
\newtheorem{lem}[thm]{Lemma}

\newtheorem{cor}[thm]{Corollary}
\newtheorem{claim}[thm]{Claim}
\renewcommand{\thestep}{}
\newtheorem{prop}[thm]{Proposition}

\theoremstyle{definition}

\renewcommand{\thecase}{}

\newtheorem{rmk}[thm]{Remark}

\theoremstyle{remark}

\makeatletter
\def\alphenumi{
  \def\theenumi{\alph{enumi}}
  \def\p@enumi{\theenumi}
  \def\labelenumi{(\@alph\c@enumi)}}
\makeatother

\makeatletter
\def\thecase{\@arabic\c@case}
\makeatother

\numberwithin{equation}{section}

\makeatletter
\def\thestep{\@arabic\c@step}
\makeatother

\newenvironment{pf}{\begin{proof}[\proofname]}{\end{proof}}
\newenvironment{pf*}[1]{\begin{proof}[#1]}{\end{proof}}

\newcommand\baralpha{{\bar\alpha}}

\newcommand\barbeta{{\bar\beta}}

\newcommand\barM{{\bar{M}}}

\newcommand\barrd{{\overline{\partial}}}

\newcommand\ubarP{{\underline{P}}}

\newcommand\fs{{\mathfrak{S}}}
\newcommand\fS{{\mathfrak{S}}}
\newcommand\ubarfS{{\underline{\mathfrak{S}}}}
\newcommand\ubarfV{{\underline{\mathfrak{V}}}}

\newcommand\CC{\mathbb{C}}

\newcommand\FF{\mathbb{F}}
\newcommand\GG{\mathbb{G}}
\newcommand\HH{\mathbb{H}}
\newcommand\II{\mathbb{I}}
\newcommand\JJ{\mathbb{J}}

\newcommand\NN{\mathbb{N}}

\newcommand\PP{\mathbb{P}}

\newcommand\RR{\mathbb{R}}

\newcommand\bE{{\mathbf{E}}}

\newcommand\bH{{\mathbf{H}}}

\newcommand\bV{{\mathbf{V}}}

\newcommand\bz{{\mathbf{z}}}

\newcommand\ssG{{{}^\circ\mathcal{G}}}

\newcommand{\cov}{\nabla}

\newcommand{\Dirac}{{\mathcal{D}}}

\newcommand\half{{\textstyle{\frac{1}{2}}}}

\newcommand\quarter{{\textstyle{\frac{1}{4}}}}
\newcommand\threehalf{{\textstyle{\frac{3}{2}}}}

\newcommand\ff{{\mathfrak{f}}}

\newcommand\fm{{\mathfrak{m}}}
\newcommand\fM{{\mathfrak{M}}}

\newcommand\fV{{\mathfrak{V}}}

\newcommand\fX{{\mathfrak{X}}}

\newcommand\de{\delta}

\newcommand\La{\Lambda}

\newcommand\om{\omega}
\newcommand\Om{\Omega}

\newcommand\hatubarP{{\widehat{\underline{P}}}}

\newcommand\gl{{\mathfrak{g}\mathfrak{l}}}

\newcommand\so{{\mathfrak{s}\mathfrak{o}}}

\newcommand\su{{\mathfrak{s}\mathfrak{u}}}
\newcommand\fu{{\mathfrak{u}}}

\newcommand\GL{\operatorname{GL}}
\newcommand\Or{\operatorname{O}}

\newcommand\PU{\operatorname{PU}}

\newcommand\SO{\operatorname{SO}}

\newcommand\Spin{\operatorname{Spin}}
\newcommand\SU{\operatorname{SU}}
\newcommand\U{\operatorname{U}}

\newcommand{\8}{\infty}

\newcommand\ad{{\operatorname{ad}}}

\newcommand\asd{{\operatorname{asd}}}

\newcommand\Cl{\operatorname{Cl}}

\newcommand\codim{\operatorname{codim}}

\newcommand\Coker{\operatorname{Coker}}

\newcommand\End{\operatorname{End}}
\newcommand\Ext{\operatorname{Ext}}
\newcommand\Flag{\operatorname{Flag}}

\newcommand\Fred{\operatorname{Fred}}

\newcommand\Hom{\operatorname{Hom}}

\newcommand\Ind{\operatorname{ind}}

\newcommand\Ker{\operatorname{Ker}}

\newcommand\Met{\operatorname{Met}}

\newcommand\Real{\operatorname{Re}}

\newcommand\Ran{\operatorname{Ran}}

\newcommand\Rank{\operatorname{Rank}}

\newcommand\Sym{\operatorname{Sym}}

\newcommand\tr{\operatorname{tr}}

\newcommand\can{{\mathrm{can}}}

\newcommand\fc{{\mathrm{fc}}}
\newcommand\id{{\mathrm{id}}}
\newcommand\LC{{\mathrm{LC}}}

\newcommand\sm{{\mathrm{sm}}}

\newcommand\spinc{\text{$\text{spin}^c$ }\allowbreak}
\newcommand\Spinc{\mathrm{Spin}^c}

\newcommand\sA{{\mathcal{A}}}

\newcommand\sC{{\mathcal{C}}}
\newcommand\sD{{\mathcal{D}}}

\newcommand\sF{{\mathcal{F}}}
\newcommand\sG{{\mathcal{G}}}

\newcommand\sJ{{\mathcal{J}}}

\newcommand\sM{{\mathcal{M}}}
\newcommand\sN{{\mathcal{N}}}

\newcommand\sP{{\mathcal{P}}}

\newcommand\sU{{\mathcal{U}}}

\newcommand\sZ{{\mathcal{Z}}}

\newcommand\tsC{{\tilde\sC}}

\newcommand\tII{{\tilde\mathbb{I}}}
\newcommand\tJJ{{\tilde\mathbb{J}}}

\newcommand\tfM{{\tilde{\mathfrak{M}}}}

\newcommand\tubarP{{\widetilde{\underline{P}}}}

\begin{document}

\title[Generic metrics and transversality]
{Generic Metrics, Irreducible Rank-One PU(2) Monopoles, and
Transversality}   
 
\author[Paul M. N. Feehan]{Paul M. N. Feehan}
\address{Department of Mathematics\\
Ohio State University\\
Columbus, OH 43210}
\email{feehan@math.ohio-state.edu}

\thanks{Supported in part by an NSF Mathematical 
Sciences Postdoctoral Fellowship under grant DMS 9306061
and by NSF grant DMS-9704174}

\date{August 3, 2000. First version: August 27, 1997.
Comm. Anal. Geom. {\bf 8} (2000), to appear; math.DG/9809001.} 

\maketitle

\section{Introduction}
Our main purpose in this article is to prove that the moduli space of
solutions to the $\PU(2)$ monopole equations is a smooth manifold of
the expected dimension for simple, generic parameters such as (and
including) the Riemannian metric on the given four-manifold: see Theorem
\ref{thm:Transversality}. In 
\cite{FL1} we proved transversality using an extension of the
holonomy-perturbation methods of Donaldson, Floer, and Taubes
\cite{DonOrient}, \cite{DonPoly}, \cite{Floer}, \cite{TauCasson}, 
together with the existence of
an Uhlenbeck compactification for the perturbed moduli space.  In
\cite{FLGeorgia}, \cite{FL2b} we discuss applications of these results to
the $\PU(2)$ monopole program for proving the equivalence between
Donaldson and Seiberg-Witten invariants conjectured in
\cite{Witten}, \cite{MooreWitten} 
(see, for example, \cite{FLGeorgia}, \cite{OTVortex}, 
\cite{OTQuaternion}, \cite{OTSurvey}, \cite{PTCambridge}, \cite{PTLocal}). 
However, it is an important and interesting question to see whether
there are simpler alternatives to the holonomy perturbations and this
is the issue we address here.

The idea that one should be able to use PU(2) monopoles to prove
Witten's conjecture concerning the relation between the two types of
four-manifold invariants was proposed by Pidstrigach and Tyurin in
1994; see \cite{FL1}, \cite{FL2b}, \cite{FLGeorgia}, \cite{OTVortex},
\cite{OTQuaternion}, \cite{OTSurvey}, \cite{PTCambridge}, \cite{PTLocal},
 \cite{TelemanMonopole} for work on this program due independently to the
author and Leness, Pidstrigach and Tyurin, and Okonek and Teleman.
The idea itself can be informally described quite quickly, the key
point being that the moduli space of PU(2) monopoles contains the
moduli space of anti-self-dual connections together with copies of the
various Seiberg-Witten moduli spaces, these forming singular loci in
the higher-dimensional moduli space of PU(2) monopoles. In principle,
then, one should be able to use intersection theory on this
higher-dimensional moduli space to relate the two kinds of
invariants. In practice, despite the simplicity of this basic idea,
the difficulties surrounding its implementation are daunting.  For
this program to succeed one needs to know that the moduli space
of $\PU(2)$ monopoles --- away from the anti-self-dual and
Seiberg-Witten points --- is a smooth manifold of the expected
dimension. This ensures that these exceptional points are the only
singularities and that the $\PU(2)$ monopole moduli space forms a
smooth (though non-compact, because of bubbling) cobordism between the
links of the singularities.

For $\PU(2)$ monopoles, we would ideally 
like an analogue of the Freed-Uhlenbeck generic metrics
theorem for the anti-self-dual equation \cite{DK}, \cite{FU} or the
generic parameter result introduced by Witten for the Seiberg-Witten
equations \cite{KMThom}, \cite{Witten}. However, results of this kind
for the $\PU(2)$ monopole equations appear to be much harder to
prove. One of the outcomes of our joint work with Leness \cite{FL1}
was a proof that one could nonetheless obtain a useful transversality
result via holonomy perturbations by extending related ideas of
Donaldson and Taubes \cite{DonOrient}, \cite{DonPoly}, 
\cite{TauCasson}. Such holonomy perturbations are 
important when considering three-manifold versions of the monopole
equations.  While all of the intersection theory
calculations on the moduli space of $\PU(2)$ monopoles described in
\cite{FL2b} could be carried out using holonomy perturbations, the
generic-parameter result (Theorem \ref{thm:Transversality})
established in the present article represents
a very significant simplification and has been the basis of our
continuing work on the project to mathematically verify Witten's
conjecture \cite{FL2b}, \cite{FL3}, \cite{FL4}. 
Our generic-parameter approach (see \S \ref{subsec:Outline}
for an outline) makes essential
use of certain unique continuation
properties for reducible solutions to the $\PU(2)$ monopole equations
which we developed in our earlier work \cite{FL1},
although not the holonomy perturbations themselves.

Issues of transversality also appear to play a significant role in the
ongoing work of Kronheimer and Mrowka to complete a three-dimensional
analogue of the Pidstrigach-Tyurin program and use $\PU(2)$ monopoles
to prove ``Property P'' for knots, via a comparison of Yang-Mills and
Seiberg-Witten Floer homologies \cite{KronheimSurvey}, \cite{KMPropertyP}.  
Moving outside the realm of
low-dimensional manifolds per se, a technical issue which plagued the
definition of Gromov-Witten invariants for general symplectic
manifolds concerned the absence of generic-parameter transversality
results for boundary components of the compactification. For
semi-positive symplectic manifolds such transversality difficulties do
not arise \cite{RuanTian}.  
The problem was eventually addressed --- via an important
differential-geometric extension of certain algebraic excess
intersection theory methods --- by Fukaya and Ono \cite{FukayaOno}, Li
and Tian \cite{LiTian}, Liu and Tian \cite{LiuTian}, Ruan
\cite{RuanGW}, and Siebert \cite{SiebertGW}, using a variety of different
approaches.  A solution was also announced by Hofer and Salamon.
In the case of $\PU(2)$
monopoles it is already a difficult problem to obtain transversality
away from the exceptional solutions, even when no bubbling has
occurred. In general, it is not possible to ensure that the loci of
exceptional solutions are unobstructed, so we must still use
differential-geometric excess intersection theory techniques
\cite{FL2b} broadly similar to those of \cite{LiTian}, \cite{RuanGW}
and going back to ideas of
\cite{AS4}, \cite{DonPoly}, \cite{DK}, \cite{MMR}, \cite{TauIndef}.

\subsection{PU(2) monopoles}
Throughout this article, $X$ denotes a closed, connected, oriented, 
smooth four-manifold. In order to state our results, we briefly
recall the description of the moduli space of $\PU(2)$ monopoles from
\cite{FL1}, \cite{FLGeorgia}. We give $X$ a Riemannian metric and 
consider Hermitian two-plane bundles $E$ over $X$ whose determinant line
bundles $\det E$ are isomorphic to a fixed Hermitian line bundle
endowed with a fixed $C^\8$, unitary connection $A_e$.  Let
$(\rho,W^+,W^-)$ be a
\spinc structure on $X$, where $\rho:T^*X\to\Hom_\CC(W^+,W^-)$ is the
Clifford map, and the Hermitian four-plane bundle $W := W^+\oplus W^-$ is
endowed with a $C^\8$, unitary connection. The unitary connection on $W$
uniquely determines a Riemannian connection on $T^*X$, via the Clifford map
$\rho$, and a unitary connection on $\det W^+$; conversely, a choice of
Riemannian connection on $T^*X$ and unitary connection on $\det W^+$ induce
a unitary connection on $W$. The connection on $W$ is called \spinc if it
induces the Levi-Civita connection on $T^*X$ for the given Riemannian metric.

Let $k\ge 2$ be an integer and let $\sA_E$ be the space of $L^2_k$
connections $A$ on the $\U(2)$ bundle $E$ all inducing the fixed
determinant connection $A_e$ on $\det E$.  Equivalently, following \cite[\S
2(i)]{KMStructure}, we may view $\sA_E$ as the space of $L^2_k$ connections
$A$ on the $\PU(2)=\SO(3)$ bundle $\su(E)$.  We shall pass back and forth
between these viewpoints, via the fixed connection on $\det E$, relying on
the context to make the distinction clear.  Given a 
unitary connection $A$ on $E$
with curvature $F_A\in L^2_{k-1}(\La^2\otimes\fu(E))$, then $(F_A^+)_0 \in
L^2_{k-1}(\La^+\otimes\su(E))$ denotes the traceless part of its self-dual
component. Equivalently, if $A$ is a connection on $\su(E)$ with curvature
$F_A\in L^2_{k-1}(\La^2\otimes\so(\su(E)))$, then $\ad^{-1}(F_A^+) \in
L^2_{k-1}(\La^+\otimes\su(E))$ is its self-dual component, viewed as a
section of $\La^+\otimes\su(E)$ via the isomorphism
$\ad:\su(E)\to\so(\su(E))$. When no confusion can arise, the isomorphism
$\ad:\su(E)\to\so(\su(E))$ will be implicit and so we regard $F_A$ as a
section of $\La^+\otimes\su(E)$ when $A$ is a connection on $\su(E)$.

For an $L^2_k$ section $\Phi$ of $W^+\otimes E$, let $\Phi^*$ be its
pointwise Hermitian dual and let $(\Phi\otimes\Phi^*)_{00}$ be the
component of the Hermitian endomorphism $\Phi\otimes\Phi^*$ of $W^+\otimes
E$ which lies in $\su(W^+)\otimes\su(E)$. The Clifford multiplication
$\rho$ defines an isomorphism $\rho:\La^+\to\su(W^+)$ and thus an
isomorphism $\rho=\rho\otimes\id_{\su(E)}$ of $\La^+\otimes\su(E)$ with
$\su(W^+)\otimes\su(E)$.

Let $\sG_E$ be the Hilbert Lie group of $L^2_{k+1}$ unitary gauge
transformations of $E$ with {\em determinant one\/}. It is often
convenient to take quotients by a slightly larger symmetry group than
$\sG_E$ when discussing pairs, so let $S_Z^1$ denote the center of $\U(2)$
and set
$\ssG_E := S_Z^1\times_{\{\pm\id_E\}}\sG_E$, 
which we may view as the group of $L^2_{k+1}$ unitary gauge transformations
of $E$ with constant determinant.  The stabilizer of a unitary
connection on $E$ in $\ssG_E$ always contains the center
$S^1_Z\subset\U(2)$. 

We call an $L^2_k$ pair $(A,\Phi)$ in the pre-configuration space,
$$
\tsC_{W,E} := \sA_E\times L^2_k(W^+\otimes E),
$$
a $\PU(2)$ monopole if $\fs(A,\Phi)=0$, where the 
$\ssG_E$-equivariant map $\fs:\tsC_{W,E}\to L^2_k(\La^+\otimes\su(E))\oplus
L^2_k(W^-\otimes E)$ is defined by
\begin{equation}
\fs(A,\Phi)
:=
\begin{pmatrix}
F_A^+ - \tau\rho^{-1}(\Phi\otimes\Phi^*)_{00} 
\\
\Dirac_A\Phi + \rho(\vartheta)\Phi
\end{pmatrix},
\label{eq:PT}
\end{equation}
where $\Dirac_A:L^2_k(W^+\otimes E)\to L^2_{k-1}(W^-\otimes E)$ is the
Dirac operator, while $\tau \in C^\8(\GL(\Lambda^+))$ and $\vartheta \in
\Omega^1(X,\CC)$ are perturbation parameters.
We let $M_{W,E} := \fs^{-1}(0)$ be the moduli space of
solutions cut out of the configuration space,
$$
\sC_{W,E} := \tsC_{W,E}/\ssG_E,
$$ 
by the section \eqref{eq:PT}, where $u\in\ssG_E$ acts by $u(A,\Phi) :=
(u_*A,u\Phi)$. 

As customary, we say that an $\SO(3)$ connection $A$ on $\su(E)$ is {\em
irreducible\/} if its stabilizer in $\sG_E$ is $\{\pm\id_E\}$,
corresponding to the center of $\SU(2)$, and {\em reducible\/} otherwise;
we say that a pair $(A,\Phi)$ on $(\su(E),W^+\otimes E)$ is irreducible
(respectively, reducible) if the connection $A$ is irreducible
(respectively, reducible). We let $\sC_{W,E}^*\subset \sC_{W,E}$ be the
open subspace of gauge-equivalence classes of irreducible pairs.  We say
that a section $\Phi$ of $W^+\otimes E$ has {\em rank $r$\/} if, when
considered as a section of $\Hom_\CC(W^{+,*},E)$, the section $\Phi(x)$ has
complex rank less than or equal to $r$ at every point $x\in X$, with equality
at some point; we say that a pair $(A,\Phi)$ on $(\su(E),W^+\otimes E)$ has
rank $r$ if the section $\Phi$ has rank $r$; if $(A,\Phi)$ has rank zero,
that is $\Phi\equiv 0$ on $X$, we call $(A,\Phi)$ a {\em zero-section\/}
pair.  We let $\sC_{W,E}^0\subset \sC_{W,E}$ be the open subspace of
gauge-equivalence classes of non-zero-section pairs and recall that
$\sC_{W,E}^{*,0} := \sC_{W,E}^*\cap\sC_{W,E}^0$ is a Hausdorff, Hilbert
manifold \cite[Proposition 2.8]{FL1} represented by pairs with stabilizer
$\{\id_E\}$ in $\ssG_E$. Let $M_{W,E}^{*,0} = M_{W,E}\cap\sC_{W,E}^{*,0}$
be the open subspace of the moduli space $M_{W,E}$ represented by
irreducible, non-zero-section $\PU(2)$ monopoles.

If $(A,\Phi)$ is a $\PU(2)$ monopole then 
(see \cite[Lemma 5.21]{FL1}) 
$$
\text{$(A,\Phi)$ reducible} \Rightarrow \text{$(A,\Phi)$ has rank less than or
equal to one}.
$$
However, it is an important observation, due to Teleman
\cite{OTSurvey}, \cite{TelemanMonopole}, that
$$
\text{$(A,\Phi)$ reducible} \nLeftarrow \text{$(A,\Phi)$ rank one}. 
$$
Indeed, counterexamples are easily constructed (at least for the
unperturbed $\PU(2)$ monopole equations) when $X$ is a K\"ahler manifold
with its canonical \spinc structure \cite{OTSurvey}, \cite{TelemanMonopole}
(see Appendix \ref{app:Teleman}). That is, a reducible
$\PU(2)$ monopole necessarily has rank less than or equal to one but in
general, {\em there exist irreducible, rank-one $\PU(2)$ monopoles\/}. If
$\Phi\equiv 0$, then the $\PU(2)$ monopole equations 
\eqref{eq:PT} just imply that $A$ is an
anti-self-dual connection on $\su(E)$; the locus of reducible $\PU(2)$
monopoles in $M_{W,E}$ can be identified with a union of Seiberg-Witten
moduli spaces. 

As we shall see in \S \ref{sec:Ranktwo} it is not too difficult to prove
that $M_{W,E}^{*,0}$ is a smooth manifold of the expected dimension away
from the locus of irreducible, rank-one solutions using the perturbation
parameters $(\tau,\vartheta)$ alone.  However, as irreducible, rank-one
solutions to \eqref{eq:PT} may be present in $M_{W,E}^{*,0}$, it seems
impossible to prove that the entire space $M_{W,E}^{*,0}$ is a smooth
manifold of the expected dimension using these parameters alone.  A similar
problem arises in the proof of transversality for the `\spinc-ASD'
equations given in
\cite[Proposition I.3.5]{PTDirac}; a version of these equations can be
obtained from the equations \eqref{eq:PT} by omitting the quadratic term
$\tau\rho^{-1}(\Phi\otimes\Phi^*)_{00}$. In the proof of \cite[Proposition
I.3.5]{PTDirac} it is claimed that if $D_A\Phi=0$ and $\Phi$ is rank one,
then $A$ is reducible [p. 277]: Teleman's counterexample shows that this
claim is incorrect and he points out an error in their argument
\cite{TelemanMonopole}.

On the other hand, the fact that the counterexamples of
\cite{TelemanMonopole} occur when $X$ is a complex, K\"ahler surface suggests that 
irreducible, rank-one solutions might not be present in the moduli
spaces $M_{W,E}^{*,0}$ for generic Riemannian metrics and compatible
Clifford maps, so this is a feature of the equations \eqref{eq:PT}
which we shall explore further in this article.

\subsection{Statement of results}
We can now state our Uhlenbeck compactness and transversality results
for the moduli space of $\PU(2)$ monopoles with the perturbations
discussed in the preceding section. As we remarked earlier,
applications of these results to the $\PU(2)$ monopole program for
proving the equivalence of Donaldson and Seiberg-Witten invariants are
described in \cite{FLGeorgia}, \cite{FL2b}.

\begin{thm}
\label{thm:Compactness}
Let $X$ be a closed, oriented, smooth four-manifold with $C^\8$ Riemannian
metric, \spinc structure $(\rho,W^+,W^-)$ with \spinc connection on $W =
W^+\oplus W^-$, and a Hermitian two-plane bundle $E$ with unitary
connection on $\det E$.  Then there is a positive integer $N_p$, depending
at most on the curvatures of the fixed connections on $W$ and $\det E$
together with $c_2(E)$, such that for all $N\ge N_p$ the topological space
$\barM_{W,E}$ is second countable, Hausdorff, compact, and given by the
closure of $M_{W,E}$ in
$\cup_{\ell=0}^{N}(M_{W,E_{-\ell}}\times\Sym^\ell(X)$ with respect to the
Uhlenbeck topology, where $E_{-\ell}$ is a Hermitian two-plane bundle over
$X$ with $\det E_{-\ell} = \det E$ and $c_2(E_{-\ell})=c_2(E)-\ell$
for each integer $\ell\ge 0$.
\end{thm}

Theorem \ref{thm:Compactness} is simply a special case of the
more general result proved in \cite{FL1} for the moduli space of solutions
to the perturbed $\PU(2)$ monopole equations in the presence of holonomy
perturbations, so no separate proof is required. 
We include the statement here since it is more accessible in
the absence of holonomy perturbations and because we appeal to it in
\cite{FL2b}, \cite{FL3}, \cite{FL4}.

\begin{rmk}
The existence of an Uhlenbeck compactification for the moduli space of
solutions to the unperturbed $\PU(2)$ monopole equations
\eqref{eq:PT} was announced by Pidstrigach
\cite{PTCambridge} and an argument was outlined in \cite{PTLocal}.
A similar argument for the equations 
\eqref{eq:PT} (without perturbations)
was outlined by Okonek and Teleman in
\cite{OTQuaternion}. An independent proof of Uhlenbeck compactness for
\eqref{eq:PT} and other perturbations of these equations is
also given in \cite{TelemanMonopole}.
\end{rmk}

The perturbation parameters $(\tau,\vartheta)$
occurring in the statement of Theorem \ref{thm:Transversality} below
appear explicitly in the $\PU(2)$ monopole equations
\eqref{eq:PT} and are described further in \S
\ref{subsec:TopStratumTransv}. The perturbation parameter 
$f \in C^\8(\GL(T^*X))$ is a
variation of the Clifford map $\rho:T^*X\to\Hom_\CC(W^+,W^-)$ and of the
Riemannian metric $g$ on $T^*X$ by an automorphism of $T^*X$ and is
described further in \S \ref{sec:VarDirac}; the perturbed $\PU(2)$
monopole equations with the three perturbation parameters
$(f,\tau,\vartheta)$ are given in 
equation \eqref{eq:DiracParamMonopole}. Let
$M^{\asd}_E$ denote the moduli space of anti-self-dual $\SO(3)$ connections on
$\su(E)$. 

\begin{thm}
\label{thm:Transversality}
Let $X$ be a closed, oriented, smooth four-manifold with $C^\8$ Riemannian
metric $g$, \spinc structure $(\rho,W^+,W^-)$ with unitary connection on 
$\det W^+$, and a Hermitian line bundle $\det E$ with unitary
connection. Then there is a first-category subset $\sP^\8_{\fc}$ of the
Fr\'echet space $\sP^\8$ of $C^\8$ perturbation parameters
$(f,\tau,\vartheta)$ such that for all $(f,\tau,\vartheta)$ in $\sP^\8 -
\sP^\8_{\fc}$ the following holds: For each parameter $(f,\tau,\vartheta)$
in $\sP^\8 - \sP^\8_{\fc}$ and Hermitian two-plane bundle $E$ over $X$, the
moduli space $M^{*,0}_{W,E}(f,\tau,\vartheta)$ of $\PU(2)$ monopoles is a
smooth manifold of the expected dimension,
\begin{align*} 
\dim M^{*,0}_{W,E}(f,\tau,\vartheta)
&=
\dim M^{*,\asd}_E + 2\Ind_\CC\sD_A -1
\\
&= 
-2p_1(\su(E))-\threehalf(e(X)+\sigma(X)) 
\\
&\quad + \half p_1(\su(E))+\half((c_1(W^+)+c_1(E))^2-\sigma(X))-1.
\end{align*}
\end{thm}

As explained further in \S \ref{subsec:CliffordDiracVar}, the Levi-Civita
connection on $T^*X$ for the metric $f^*g$, Clifford map $\rho\circ f$, and
unitary connection on $\det W^+$ define a \spinc connection on $W$.  The
rest of our article is devoted to proving Theorem
\ref{thm:Transversality}.

\subsection{Outline of the proof and of the remainder of the article}
\label{subsec:Outline}
We prove Theorem \ref{thm:Transversality} in two steps, both of which have
analogues in the case of the moduli space of anti-self-dual $\SO(3)$
connections over a four-manifold $X$ with $b^+(X) > 0$ (although that
constraint is not required here): 
\begin{enumerate}
\item For generic
parameters $(\tau,\vartheta)$ (and any parameter $f$), we show
that the moduli space $M_{W,E}^{*,\natural}(f,\tau,\vartheta)$ of
$\PU(2)$ 
monopoles $(A,\Phi)$ with $A$ irreducible and $\Phi$ rank two is a smooth
manifold of the expected dimension. 
\item For generic parameters
$(f,\vartheta)$ (and any parameter $\tau$), we show that the moduli
space $M_{W,E}^{*,0}(f,\tau,\vartheta)$ of $\PU(2)$ monopoles
$(A,\Phi)$ with $A$ irreducible and $\Phi\not\equiv 0$ contains no pairs
$(A,\Phi)$ with $\Phi$ having rank one, that is, 
$M_{W,E}^{*,0}(f,\tau,\vartheta)
= M_{W,E}^{*,\natural}(f,\tau,\vartheta)$.
\end{enumerate}
The first result, Theorem
\ref{thm:GenericProjection}, is proved in \S \ref{sec:Ranktwo} using
the Sard-Smale theorem \cite{SmaleSard} and Aronszajn's unique continuation
theorem \cite{Aronszajn}, together with some linear algebra calculations
reminiscent of those arising in the proof of the Freed-Uhlenbeck generic
metrics theorem \cite{FU}. The second result, Theorem
\ref{thm:NoRankOneIrreducibles}, is considerably
more difficult and is proved
in \S \ref{sec:Period} using an infinite-dimensional analogue, for the
family of perturbed Dirac operators $\Dirac_{A,f,\vartheta}$, of the
well-known `period map' argument for the family of operators $d^{+,g}$
\cite[\S VI]{DonConn}, \cite[Corollary 4.3.15]{DK}; see also
\cite[Corollary 3.21]{FU}, \cite[Appendix B]{TauIndefPrelim}.  Combining
these two results yields Theorem \ref{thm:Transversality}. In the case of
the moduli space of $\SO(3)$ anti-self-dual connections, the first step is
accomplished by the Freed-Uhlenbeck generic metrics theorem: that is, for
generic metrics, $M_E^{\asd,*}(g)$ is a smooth manifold of the expected
dimension \cite[Corollary 4.3.18]{DK}, \cite[Theorem 3.17]{FU}. The second
step is achieved by the (finite-dimensional) period-map argument: for
generic metrics and $b^+(X)>0$, the moduli space $M_E^{\asd}(g)$ contains
no reducible connections \cite[\S VI]{DonConn}, \cite[Corollary
4.3.15]{DK}.

The detailed proof of Theorem \ref{thm:Transversality} is, of course, more
complicated than the preceding synopsis can convey, so we indicate how the
remainder of this article is organized. Section \ref{sec:Ranktwo} contains
the proof that the moduli space of irreducible, rank-two $\PU(2)$ monopoles
is regular. In \S \ref{subsec:TopStratumTransv} we define a
parametrized moduli space of $\PU(2)$ monopoles and, assuming this is
regular away from the loci of reducible or lower-rank $\PU(2)$ monopoles,
we use the Sard-Smale theorem to show that the individual moduli spaces of
irreducible, rank-two $\PU(2)$ monopoles are regular for generic
parameters. The heart of the first step, then, is given in \S
\ref{subsec:SmoothParamModuli}, where we show that the parametrized moduli
space of irreducible, rank-two $\PU(2)$ monopoles is regular. Section
\ref{sec:VarDirac} describes the perturbations of the Dirac operator: in \S
\ref{subsec:CliffordDiracVar} we compute the differential of the resulting
family of Dirac operators and in \S \ref{subsec:ParamModuliDirac} we define
the $\PU(2)$ monopole equations and associated parametrized moduli space
with the full set of perturbations. Section \ref{sec:Period} contains the
proof that, for generic parameters, the moduli spaces of
irreducible $\PU(2)$ monopoles contain no rank-one pairs. Sections
\ref{subsec:Grassmann}, \ref{subsec:Fredholm}, and \ref{subsec:RankOneLoci}
describe the loci of irreducible, rank-one $\PU(2)$ monopoles in terms of
incidence correspondences --- with infinite codimension --- in certain
Banach flag and Grassman manifolds, by analogy with standard
finite-dimensional constructions of algebraic geometry \cite{Harris}. Section
\ref{subsec:Grassmann} describes general incidence correspondences using
Banach Grassman manifolds and \S
\ref{subsec:Fredholm} reinterprets these correspondences using spaces of
Fredholm operators, which are more suitable for our purposes. In \S
\ref{subsec:RankOneLoci} we define the loci of irreducible, rank-one
$\PU(2)$ monopoles in terms of these  incidence correspondences. These
rank-one loci are in turn detected with aid of the Dirac-operator period
map, which is defined and whose differential is shown to have the requisite
surjectivity properties in \S \ref{subsec:DiffPeriodMap}. Section
\ref{subsec:SardSmale} contains a proof of the parametric transversality
theorem: the result is well-known, but we include the proof as no single
reference contains the precise statement (and proof) we need. Finally, in \S
\ref{subsec:Kuranishi}, we complete the proof that the moduli spaces of
irreducible $\PU(2)$ monopoles contain no rank-one pairs for generic
parameters: we apply the parametric transversality theorem to show, in
essence, that the loci of irreducible, rank-one $\PU(2)$ monopoles must
have infinite codimension in their respective moduli spaces for generic
parameters and so these loci are empty. The
stratification of the space of Fredholm operators (by kernel dimension)
\cite{Koschorke}, employed in \S \ref{subsec:Kuranishi}, 
was also used by Maier in \cite{Maier}. At the beginning of \S
\ref{subsec:Kuranishi} we first give an outline of the main argument completing
the proof of Theorem \ref{thm:NoRankOneIrreducibles}, as the detailed
argument is rather long and technical.

Theorem \ref{thm:NoRankOneIrreducibles} serves the auxiliary purpose of
giving another method of completing the gap in the transversality
argument used by Pidstrigach-Tyurin in their definition of the \spinc
polynomial invariants \cite[\S I]{PTDirac}; the same end is achieved
via the holonomy perturbations of \cite{FL1}. We briefly describe this
application and Teleman's counterexample in the Appendix, together with
some facts we need from linear algebra.

\subsection{Other approaches to transversality}
As we remarked at the beginning of the Introduction, versions of Theorem 
\ref{thm:Compactness} and \ref{thm:Transversality} were proved by the
author and Leness in \cite{FL1} for the $\PU(2)$ monopole equations
\eqref{eq:PT} with additional holonomy perturbations; see
\cite{FLGeorgia} for a more concise account of this method.  The possible
presence or absence of irreducible, rank-one solutions to the $\PU(2)$
monopole equations considered in \cite{FL1} makes no difference to the
argument there, as the perturbations are strong enough to yield
transversality without a separate analysis of the locus of irreducible,
rank-one solutions. A preliminary version \cite{FL96} of the article
\cite{FL1} relied on the incorrect assertion of \cite[p. 277]{PTDirac}
described above and used only the perturbation parameters
$(\tau,\vartheta)$. Transversality results for the $\PU(2)$ monopole
equations, with perturbation parameters including $(\tau,\vartheta)$ and
the Riemannian metric $g$ on $T^*X$, were conjectured by Pidstrigach and
Tyurin in \cite{PTCambridge}, \cite{PTLocal}.

Teleman has explored another approach to the transversality
problem, quite different from those of \cite{FL1}, \cite{PTLocal}, using
certain {\em ad hoc\/} perturbations of the principal symbol of the
linearization of the $\PU(2)$ monopole equations
\cite{TelemanMonopole}. Like the 
$\PU(2)$ monopole equations \eqref{eq:PT}, the equations of
\cite{TelemanMonopole} employ the perturbation parameters 
$(\tau,\vartheta)$, together with a term of the form
$\sum_i\rho^{-1}\tau_i\cov_{A,Y_i}(\Phi\otimes\Phi^*)_{00}$ in the
curvature equation in \eqref{eq:PT}, where $\{Y_i\}$ is a finite set of
vector fields spanning $TX$ at every point of $X$ and the coefficients
$\tau_i$ are in $\Omega^0(\gl(\su(W^+)))$. The approach of
\cite{TelemanMonopole} illustrates the significance of principal-symbol
perturbations. However, as noted in the introduction to
\cite{TelemanGenericMetric}, it appears to be difficult or impossible to
show that solutions to the $\PU(2)$ monopole equations of
\cite{TelemanMonopole} which are reducible on an open subset of $X$ are
necessarily 
reducible over all of $X$, as required by the transversality argument in
\cite{TelemanMonopole}. Indeed, Proposition 3.1.2 in \cite{TelemanMonopole}
assumes 
incorrectly that that the Agmon-Nirenberg unique continuation theorem for
an ordinary differential equation on a Hilbert space \cite{AgmonNirenberg}
applies to 
the $\PU(2)$ monopole equations of
\cite{TelemanMonopole}, as an examination of the hypotheses of
\cite{AgmonNirenberg} reveals; 
the Agmon-Nirenberg theorem is used by Donaldson-Kronheimer in their proof
of the corresponding unique continuation property for solutions to the
anti-self-dual equation and by the author and Leness for our proof of a
restricted unique continuation property for the holonomy-perturbed $\PU(2)$
monopole equations of \cite{FL1}. Varying the Riemannian metric $g$ on
$T^*X$ and the Clifford map $\rho:T^*X\to\Hom_\CC(W^+,W^-)$ by
automorphisms $f$ of $T^*X$, as we do in \S \ref{sec:VarDirac}, perturbs
the principal symbol of the Dirac operator; the key application of this
perturbation occurs in the proof of Proposition
\ref{prop:CokerDifferentialPeriod}, where we show that a certain partial
differential of the Dirac-operator period map is surjective.
Variations of the Dirac operator with
respect to the Riemannian metric have been described by Bourguignon and
Gauduchon in \cite{BG}. Their technique, with some enhancements, was
recently used by Maier to prove certain generic metrics results for Dirac
operators on low-dimensional spin or \spinc manifolds \cite{Maier}. The
variational formulas of \cite{BG}, \cite{Maier} were helpful for the
development of our present article, though the variation of the Dirac
operator we construct in \S \ref{sec:VarDirac} is quite different from
that of \cite{BG}. The unique continuation property for reducible
solutions to the anti-self-dual equation is an important ingredient in
the proof of the Freed-Uhlenbeck generic metrics theorem given in
\cite{DK} (see their Lemma 4.3.21); it is proved for the $\PU(2)$
monopole equations \eqref{eq:PT} by the author and Leness in
\cite{FL1} (see Theorem \ref{thm:LocalToGlobalReducible} in
\ref{subsec:DiffPeriodMap} here for a precise statement) and plays an
important role in the proof of Theorem \ref{thm:Transversality} here
as well. 

The present article was completed by August 1997, at which time it
was widely circulated and submitted to a print journal. Shortly after
this, the author received a preprint \cite{TelemanGenericMetric} from
A. Teleman addressing, independently, the issue of transversality for
the $\PU(2)$ monopole equations with perturbations similar to those
we employ here.

\section*{Acknowledgements}
It is a pleasure to thank Thomas Leness for his generous support during the
preparation of this article. We are very grateful to Tom Mrowka and
Peter Ozsv\'ath for their encouragement and many helpful comments.  We
warmly thank Clifford Taubes for his help and encouragement over the years;
we also thank him and Dieter Kotschick for bringing the work of
Stefan Maier concerning generic metrics results for Dirac operators to our
attention and thank Stefan for sending us the preprint version of
\cite{Maier}. The preprint \cite{TelemanMonopole} of Andrei Teleman played a
crucial role in the development of our article by pinpointing the gap in
\cite{PTDirac}.  We would like to thank the Mathematics
Departments at Harvard and Ohio State University and the National Science
Foundation for their generous sponsorship during the preparation of this
article.


\section{Transversality on the complement of the loci of lower-rank
PU(2) monopoles} 
\label{sec:Ranktwo}
Let $\sC^{*,\natural}_{W,E}\subset \sC_{W,E}$ be the open subspace given by
gauge-equivalence classes of pairs $(A,\Phi)$ with $A$ irreducible and
$\Phi$ rank two, and set
$$
M^{*,\natural}_{W,E} := M_{W,E}\cap \sC^{*,\natural}_{W,E}.
$$
Our goal in this section is to show that
$M^{*,\natural}_{W,E}(\tau,\vartheta)$ is a regular manifold of the expected
dimension for generic $C^\8$ parameters $(\tau,\vartheta)$: this result,
Theorem \ref{thm:GenericProjection}, is proved in \S
\ref{subsec:TopStratumTransv} under the assumption
that the parametrized moduli space is regular while the latter result,
Theorem \ref{thm:SmoothParamModuliSpace}, is proved in \S
\ref{subsec:SmoothParamModuli}.

\subsection{The parametrized moduli space}
\label{subsec:TopStratumTransv}
Our argument relies, as with most standard applications in gauge theory
\cite{DK}, \cite{FU}, \cite{KMThom}, \cite{Witten}, on the Sard-Smale theorem
for a Fredholm map of Banach manifolds \cite{SmaleSard}.

We shall need to be precise here and throughout our article about the
sense in which a parameter is `generic'. A subset $S$ of a topological
space $\sP$ is called a set of the {\em first category}
\footnote{A {\em second-category} subset is not necessarily the complement of a
first-category subset --- although the term is sometimes inaccurately used
that way: it is simply a subset which is not of the first category
\cite{Rudin}. To eliminate any possible confusion, we avoid using the term
here.} if its complement $\sP-S$ is a countable intersection of dense 
open sets or, equivalently, if $S$ is a countable union of closed subsets
of $\sP$ with empty interior; if $\sP$ is a complete metric space, then
Baire's theorem implies that $\sP-S$ is dense in $\sP$ \cite{Rudin}. In our
applications, $\sP$ will either be a Banach or Fr\'echet manifold (with
a complete metric), so $\sP-S$ will always be dense if $S$ is a
first-category subset.

Throughout this section we use a 
Banach manifold of $C^r$ perturbation parameters (with $r$ large) given by
\begin{equation}
\sP := C^r(\GL(\La^+)) \oplus C^r(\La^1\otimes\CC).
\label{eq:RankTwoParamSpace}
\end{equation}
We could, of course, also include the space $C^r(\GL(T^*X))$ of perturbations
of the Riemannian metric on $X$ and of the Clifford map
in our parameter space \eqref{eq:RankTwoParamSpace}, as the
latter perturbations are required for our period-map argument in \S
\ref{sec:Period}. However, these additional parameters are not
used to prove the main result of this section (Theorem
\ref{thm:GenericProjection}), so we omit them for the sake of clarity. In
the same vein, the space $C^r(\GL(\La^+))$ of perturbation parameters
$\tau$ is not required in our proof of the main result of \S
\ref{sec:Period} (Theorem \ref{thm:NoRankOneIrreducibles}), so in that
section we omit the space $C^r(\GL(\La^+))$ from the definition of the
Banach manifold of perturbation parameters used there, again for reasons of
clarity. To conclude the proof of our main transversality result (Theorem
\ref{thm:Transversality}), though, we can clearly assume that the same
space of perturbations (containing the three types of parameters) has been
used throughout sections \ref{sec:Ranktwo} and \ref{sec:Period}.

We define a $\ssG_E$-equivariant map
$$
\ubarfS:= (\ubarfS_1,\ubarfS_2)
:\sP\times \tsC_{W,E}
\to L^2_{k-1}(\La^+\otimes(\su(E))\oplus L^2_{k-1}(W^-\otimes E)
$$
by setting
\begin{equation}
\label{eq:ParamMonopole}
\ubarfS(\tau,\vartheta,A,\Phi) 
=
\begin{pmatrix}
\ubarfS_1(\tau,\vartheta,A,\Phi) \\
\ubarfS_2(\tau,\vartheta,A,\Phi)
\end{pmatrix} 
:= 
\begin{pmatrix}
F_A^+ - \tau\rho^{-1}(\Phi\otimes\Phi^*)_{00} \\
\Dirac_A\Phi+\rho(\vartheta)\Phi
\end{pmatrix}, 
\end{equation}
where $(A,\Phi)$ is a pair on $(\su(E),W^+\otimes E)$ and the isomorphism
$\ad:\su(E)\simeq \so(\su(E))$ is implicit, $\ssG_E$ acts trivially on the
space of perturbations $\sP$, and so $\ubarfS^{-1}(0)/\ssG_E$ is a subset
of $\sP\times\sC_{W,E}$.  We let $\fM_{W,E}$ denote the {\em parametrized
moduli space\/} $\ubarfS^{-1}(0)/\ssG_E$ and let
$\fM_{W,E}^{*,\natural}:=\fM_{W,E}\cap (\sP\times \sC_{W,E}^{*,\natural})$.
We fix a $C^\8$ Riemannian metric $g$ on $T^*X$, noting that it is not
varied in this section.

The $\ssG_E$-equivariant map $\ubarfS$
defines a section of a Banach vector bundle 
$\ubarfV$ over
$\sP\times \sC_{W,E}^{*,\natural}$ with total space
$$
\ubarfV := \sP\times\tsC_{W,E}^{*,\natural}\times_{\ssG_E}
\left(L^2_{k-1}(\La^+\otimes\su(E))\oplus L^2_{k-1}(W^-\otimes E)\right), 
$$
so $\fS :=
\ubarfS(\tau,\vartheta,\cdot)$ is a section over $\sC_{W,E}^{*,\natural}$ of
the Banach vector bundle $\fV:=\ubarfV|_{(\tau,\vartheta)}$.  In
particular, the parametrized moduli space $\fM_{W,E}^{*,\natural}$ is the
zero set of the section $\ubarfS$ of the vector bundle $\ubarfV$ over
$\sP\times
\sC_{W,E}^{*,\natural}$.

\begin{thm}\label{thm:SmoothParamModuliSpace}
The zero set in $\sP\times \sC_{W,E}^{*,\natural}$ of the section $\ubarfS$
is regular and, in particular, the moduli space $\fM^{*,\natural}_{W,E}$ is
a smooth Banach submanifold of $\sP\times \sC_{W,E}^{*,\natural}$.
\end{thm}

To preserve continuity, we defer the proof of Theorem
\ref{thm:SmoothParamModuliSpace} to \S \ref{subsec:SmoothParamModuli}.
The differential $D\ubarfS := (D\ubarfS)_{[\tau,\vartheta,A,\Phi]}$ 
of the section $\ubarfS$ at a point
$[\tau,\vartheta,A,\Phi]$ in 
$\sP\times\sC^{*,\natural}_{W,E}$ is given by
\begin{equation}
\label{eq:LinearDfS}
D\ubarfS(\delta\tau,\delta\vartheta,a,\phi)
= \begin{pmatrix}
D\ubarfS_1(\delta\tau,\delta\vartheta,a,\phi) \\
D\ubarfS_2(\delta\tau,\delta\vartheta,a,\phi) 
\end{pmatrix}, 
\end{equation}
where $(a,\phi)\in \Ker d_{A,\Phi}^{0,*}
\subset L^2_k(\La^1\otimes\su(E))\oplus L^2_k(W^+\otimes E)$
represents a vector in the 
tangent space $(T\sC^{*,\natural}_{W,E})_{[A,\Phi]}
= T_{[A,\Phi]}\sC^{*,\natural}_{W,E}$
and $(\delta\tau,\delta\vartheta)\in\sP$. Here, $d_{A,\Phi}^{0,*}$ is the $L^2$
adjoint of the differential 
$$
d_{A,\Phi}^0:L^2_{k+1}(\su(E))\oplus i\RR_Z
\to 
L^2_k(\su(E))\oplus L^2_k(W^+\otimes E)
$$
of the map $\ssG_E\to \tsC_{W,E}$, $u\mapsto u(A,\Phi)$, at the identity.

The differentials in \eqref{eq:LinearDfS} are given explicitly by
\begin{align}
D\ubarfS_1\left(\delta\tau,\delta\vartheta,a,\phi\right)
&= 
d^+_Aa -(\delta\tau)\tau\rho^{-1}(\Phi\otimes\Phi^*)_{00} 
\label{eq:LinearOfS1}\\
&\quad-\tau\rho^{-1}(\phi\otimes \Phi^*+\Phi\otimes\phi^*)_{00}, 
\notag
\\
D\ubarfS_2\left(\delta\tau,\delta\vartheta,a,\phi\right) 
&=
\Dirac_A\phi + \rho(\vartheta)\phi + \rho(a)\Phi +\rho(\delta\vartheta)\Phi.
\label{eq:LinearOfS2}
\end{align}
We have $D\ubarfS(\cdot,d_{A,\Phi}^0\zeta)=0$ for all
$\zeta\in L^2_{k+1}(\su(E))\oplus i\RR_Z$ since $\ubarfS$ is
$\ssG_E$-equivariant. By the regularity results of \cite[\S 3]{FL1}
we may assume, without loss of generality, that
the pair $(A,\Phi)$ in $\tsC_{W,E}^{*,\natural}$ is a $C^r$ representative for
the point $[A,\Phi]$ in the zero set 
$\fS^{-1}(0)\subset\sC_{W,E}^{*,\natural}$.  
Since the tangent space  
$(T\sC^{*,\natural}_{W,E})_{[A,\Phi]}$ may be
identified with $\Ker d_{A,\Phi}^{0,*}$ by the slice result
\cite[Proposition 2.8]{FL1}, we have
$$
D\ubarfS(0,0,a,\phi) := d_{A,\Phi}^1(a,\phi) 
= (d_{A,\Phi}^{0,*} + d_{A,\Phi}^1)(a,\phi),
$$
for $(a,\phi)\in \Ker d_{A,\Phi}^{0,*}$, so the differential 
$D\ubarfS|_{\{0\}\times T\sC^{*,\natural}_{W,E}}$ is Fredholm, where 
$\{0\}\times T\sC^{0,*}_{W,E} = T(\{\tau,\vartheta\}\times
\sC^{*,\natural}_{W,E})$.  We recall that $d_{A,\Phi}^0$ and $d_{A,\Phi}^1$
are the two differentials in the elliptic deformation complex for the
$\PU2)$ monopole equations \eqref{eq:PT}.
Thus, $\ubarfS$ is a Fredholm section when restricted to the fixed-parameter
fibers $\{\tau,\vartheta\}\times\sC^{*,\natural}_{W,E}\subset
\sP\times\sC^{*,\natural}_{W,E}$ and so the Sard-Smale theorem (in the form of
Proposition 4.3.11 in
\cite{DK}) implies that there is a first-category subset $\sP_{\fc}\subset\sP$
such that the zero sets in $\sC^{*,\natural}_{W,E}$ of the sections 
$\fS := \ubarfS(\tau,\vartheta,\cdot)$
are regular for all perturbations $(\tau,\vartheta)\in\sP-\sP_{\fc}$. Now 
$$
M^{*,\natural}_{W,E}(\tau,\vartheta)
=  
\fS^{-1}(0)\cap\sC^{*,\natural}_{W,E}
$$
and so for generic $C^r$ parameters $(\tau,\vartheta)$, the moduli space
$M^{*,\natural}_{W,E}(\tau,\vartheta)$ is a smooth manifold. Finally, the
argument 
of \cite[\S 5.1.2]{FL1} implies that the parameters $(g,\tau,\vartheta)$ can be
assumed without loss of generality to be $C^\8$ if $X$ is a smooth manifold.
Let
$$
\sP^\8 := C^\8(\GL(\La^+)) \times C^\8(\La^1\otimes\CC),
$$
denote the Fr\'echet manifold of $C^\8$ perturbation parameters. The
expected dimension of $M^{*,\natural}_{W,E}(\tau,\vartheta)$ is given by
Proposition 2.28 in \cite{FL1} and is obtained by computing the index of
the elliptic deformation complex for the $\PU(2)$ monopole equations
\eqref{eq:PT}. In summary, we have:

\begin{thm}\label{thm:GenericProjection}
Let $X$ be a closed, oriented, smooth four-manifold with $C^\8$ Riemannian
metric. Then there is a first-category subset $\sP^\8_{\fc} \subset \sP^\8$
such that for all $(\tau,\vartheta)$ in $\sP^\8 - \sP^\8_{\fc}$ the
following holds: The zero set of the section
$\fS:=\ubarfS(\tau,\vartheta,\cdot)$ is regular and the moduli space
$M^{*,\natural}_{W,E}(\tau,\vartheta) =
\fS^{-1}(0)\cap\sC^{*,\natural}_{W,E}$ is a smooth submanifold of
of the expected dimension.
\end{thm}

\subsection{Smoothness of the parametrized moduli space}
\label{subsec:SmoothParamModuli}
We prove Theorem \ref{thm:SmoothParamModuliSpace} in this section by
showing that the $\ssG_E$-equivariant map
$\ubarfS:\sP
\times\tsC_{W,E}^{*,\natural}\to L^2_{k-1}(\La^+\otimes\su(E))\oplus
L^2_{k-1}(W^-\otimes E)$ vanishes transversely and so the parametrized
moduli space $\fM^{*,\natural}_{W,E}=\ubarfS^{-1}(0)/\ssG_E$ is a smooth Banach
manifold. We may suppose without loss of generality that
$(\tau,\vartheta,A,\Phi)$ is a $C^r$ representative for a point in
$\ubarfS^{-1}(0)$ and denote $D\ubarfS := (D\ubarfS)_{\tau,\vartheta,A,\Phi}$
for convenience.

Note that $\Ran D\ubarfS_2$ is a real subspace of $L^2_{k-1}(W^-\otimes E)$
according to \eqref{eq:LinearOfS2}: indeed, $D\ubarfS_2(\delta\tau,0,a,0)$
spans a real subspace as $a$ varies, while
$D\ubarfS_2(\delta\tau,\delta\vartheta,0,\phi)$ spans a complex subspace as
$(\delta,\vartheta,\phi)$ varies. Hence,
$$
\Ran(D\ubarfS)\subsetneqq
L^2_{k-1}(\La^+\otimes\su(E))\oplus L^2_{k-1}(W^-\otimes E)
$$ 
if and only if there exists a non-zero pair $(v,\psi)$ such that for all
$(\delta\tau,\delta\vartheta,a,\phi)$ we have
\begin{align}
&(D\ubarfS(\delta\tau,\delta\vartheta,a,\phi),(v,\psi))_{L^2(X)}  
\label{eq:CokernelCondition}\\
&= (D\ubarfS_1(\delta\tau,\delta\vartheta,a,\phi),v)_{L^2(X)}
+ \Real(D\ubarfS_2(\delta\tau,\delta\vartheta,a,\phi),\psi)_{L^2(X)} = 0.
\notag
\end{align}
Suppose that $(v,\psi)\in L^2_{k-1}(\La^+\otimes\su(E))\oplus
L^2_{k-1}(W^-\otimes E)$ is real $L^2$-orthogonal to 
$\Ran D\ubarfS$, so $(v,\psi)$ lies in $\Ker(D\ubarfS)^*$ by
\eqref{eq:CokernelCondition}. Then
elliptic regularity for the Laplacian $D\ubarfS(D\ubarfS)^*$, with
$C^{r-1}$ coefficients, implies that $(v,\psi)$ is in
$C^{r+1}$ \cite[\S 3]{FL1}. Moreover, Aronszajn's theorem \cite[Remark 3,
p. 248]{Aronszajn}, \cite[Theorem 1.8]{Kazdan} implies that pairs
$(v,\psi)$ in the kernel of $D\ubarfS(D\ubarfS)^*$ have the unique continuation
property \cite[Lemma 5.9]{FL1}. 

We begin with a simple observation from linear algebra; though elementary,
it has an important application in the Freed-Uhlenbeck proof of the
generic metrics theorem for the anti-self-dual equation (see, for example,
Lemma 3.7 in \cite{FU}). We shall make extensive use of it here.

\begin{lem}
\label{lem:OrthogImage}
Let $U$, $V$, $W$ be finite-dimensional Hilbert spaces (either all real or
all complex) with $\dim U\leq \dim V$. Suppose $M\in \Hom(U,W) = W\otimes
U^*$ and $N \in \Hom(V,W) = W\otimes V^*$.  If $\langle MP, N
\rangle_{W\otimes V^*} = 0$ for all $P\in\Hom(V,U)$, then $\Ran M \perp
\Ran N$ in $W$ and so $\Rank M + \Rank N \leq \dim W$.
\end{lem}
\vfil\eject

\begin{proof}
Let $\{v_i\}$ be an orthonormal basis for $V$, with dual orthonormal basis
$\{v_i^*\}$ for $V^*$ given by $v_i^* := \langle\cdot,v_i\rangle$, and let
$\{u_j\}$ be an orthonormal basis for $U$. We then 
have  
\begin{align*}
0 &= \langle MP, N \rangle_{W\otimes V^*} 
\\
&= \sum_{i,j}\langle MPv_i\otimes v_i^*, 
Nv_j\otimes v_j^* \rangle_{W\otimes V^*} 
\\
&= \sum_{i,j}\langle MPv_i, Nv_j\rangle_W\cdot
\langle v_i^*, v_j^* \rangle_{V^*} 
\\
&= \sum_i\langle MPv_i, Nv_i\rangle_W.
\end{align*}
Since $P\in\Hom(V,U)$ is arbitrary, for each $i'\in\{1,\dots,\dim V\}$
and $j\in\{1,\dots,\dim U\}$
we may choose $P = P_{i'j}$ such that $P_{i'j}v_i=0$ for all $i\neq i'$ and
$P_{i'j}v_{i'} = u_j$, so the preceding identity implies that
$$
\langle Mu_j, Nv_{i'}\rangle_W = 0,
\quad\text{for all }1\le i'\le\dim V \text{ and }1\le j\le\dim U. 
$$
Hence, $\langle Mu, Nv\rangle_W = 0$ for all $u\in U$, $v\in V$, and the
assertion follows. 
\end{proof}

Setting $(\de\tau,a,\phi)=0$ in 
\eqref{eq:LinearOfS1}, \eqref{eq:LinearOfS2}, and
\eqref{eq:CokernelCondition} yields 
$$
\Real(\rho(\de\vartheta)\Phi,\psi)_{L^2(X)} = 0,
\quad
\text{for all }\delta\vartheta\in C^r(\La^1\otimes\CC).
$$
Using the identity $z=\Real(z)+i\Real(-iz)$ for $z\in\CC$, we see that
$$
(\rho(\de\vartheta)\Phi,\psi)_{L^2(X)}
=
\Real(\rho(\de\vartheta)\Phi,\psi)_{L^2(X)}
+
i\Real(\rho(-i\de\vartheta)\Phi,\psi)_{L^2(X)},
$$
and so we obtain the full complex $L^2$-orthogonality condition
\begin{equation}
(\rho(\de\vartheta)\Phi,\psi)_{L^2(X)} = 0,
\quad
\text{for all }\delta\vartheta\in C^r(\La^1\otimes\CC).
\label{eq:PhiPsiOrthog}
\end{equation}
Recall that the ranks of sections of the complex bundle $W^+\otimes E$ and
real bundles $\su(W^+)\otimes\su(E)$ and $\La^+\otimes \su(E)$ are defined
by considering them as sections of $\Hom_\CC(W^{+,*},E)$,
$\Hom_\RR(\su(W^+)^*,\su(E))$, and $\Hom_\RR(\La^{+,*},\su(E))$,
respectively.

\begin{lem}
\label{lem:HermitianOrthog}
\cite{FL96}
Suppose $\Phi\in C^0(X,W^+\otimes E)$ and $\Psi\in C^0(X,W^-\otimes E)$
satisfy \eqref{eq:PhiPsiOrthog}. Then $\Phi$ and $\Psi$, considered as
elements of $\Hom_\CC(W^{+,*},E)$ and $\Hom_\CC(W^{-,*},E)$
respectively, have complex-orthogonal images in $E$ at every point of $X$
and thus $\Rank_\CC\Phi(x) + \Rank_\CC\Psi(x) \leq 2$ at each point $x\in X$.
\end{lem}

\begin{pf}
Since $\vartheta$ is an arbitrary, $C^r$ complex-valued one-form, we obtain
the pointwise identity $\langle\rho(\vartheta_x)\Phi_x,\Psi_x\rangle_x=0$
for all $\vartheta_x\in (T^*X)_x\otimes\CC$. Let $\{\phi_1,\phi_2\}$,
$\{\psi_1,\psi_2\}$ be orthonormal frames for $W^+|_x$, $W^-|_x$,
respectively. The Clifford map is complex linear and restricts to give a
complex-linear isomorphism \cite[p. 89]{MorganSWNotes},
$$
\rho:(T^*X)_x\otimes_\RR\CC\to\Hom_\CC(W^+,W^-)_x.
$$
(This is proved as Lemma \ref{lem:MorganLemma} in Appendix
\ref{subsec:LinearAlgebra}.)  The conclusion now follows from the
first assertion in Lemma \ref{lem:OrthogImage}.  
\end{pf}

Since $(A,\Phi)\in\tsC^{*,\natural}_{W,E}$ by hypothesis and $\Phi$ is $C^r$,
there is a non-empty open subset $U\subset X$ on which
$\Rank_\CC\Phi(x)=2$, for all $x\in U$. Equation
\eqref{eq:PhiPsiOrthog} and Lemma \ref{lem:HermitianOrthog}
imply that $\Rank_\CC\psi(x) = 0$ for all $x\in U$, that is, $\psi \equiv
0$ on $U$. Similarly, by varying $\tau$ we see that \eqref{eq:LinearOfS1},
\eqref{eq:LinearOfS2}, and \eqref{eq:CokernelCondition} yield
\begin{equation}
((\delta\tau)\tau\rho^{-1}(\Phi\otimes\Phi^*)_{00},v)_{L^2(X)}=0,
\quad
\text{for all } \de\tau\in C^r(\gl(\La^+)),
\label{eq:QuadPhiZero}
\end{equation} 
by setting $(\delta\vartheta,a,\phi)=0$.

\begin{rmk}
Even though $\Phi$ solves an elliptic equation, namely
$(\Dirac_A+\rho(\vartheta))\Phi = 0$, and $\Rank_\CC\Phi(x)=2$ for all $x$ in
the open set $U$, it does not necessarily follow that $\Rank_\CC\Phi(x)=2$
at all points $x$ in $X$; see \cite{JinKazdan} and \cite[p. 668]{Kazdan}
for counterexamples in the case of harmonic maps. If all the perturbation
parameters are analytic and $(X,g)$ is a real-analytic Riemannian manifold
then $\Phi$ will necessarily be real-analytic by \cite[Theorem 6.6.1]{Morrey},
so $\det\Phi$ (see \S \ref{sec:Period}) will be a complex-valued,
real-analytic function and if it vanishes on a non-empty open set, it will
vanish identically. However, we cannot constrain our parameters to be
real-analytic before applying the Sard-Smale theorem (as a space of
real-analytic parameters would not have a complete metric), so we make no
use of real-analyticity here.
\end{rmk}

\begin{lem}
\label{lem:vandQuadOrthogImages}
\cite{FL96}
If $v\in C^0(\La^+\otimes\su(E))$ and $\Phi\in C^0(W^+\otimes E)$
satisfy \eqref{eq:QuadPhiZero}, then $v,
\tau\rho^{-1}(\Phi\otimes\Phi^*)_{00}\in\Hom_\RR(\La^{+,*},\su(E))$ have
orthogonal images in $\su(E)$ at every point in $X$ and so $\Rank_\RR v(x) +
\Rank_\RR(\Phi(x)\otimes\Phi(x)^*)_{00} \leq 3$ at each $x\in X$.
\end{lem}

\begin{pf}
Since $\delta\tau\in C^r(\gl(\La^+))$ is arbitrary, we obtain
the pointwise identity
$$
\langle 
(\delta\tau_x)\tau\rho^{-1}(\Phi\otimes\Phi^*)_{00},v\rangle_x=0
\quad\text{for all $\delta\tau_x\in\gl(\Lambda^+|_x)$},
$$
and any $x\in X$.  The conclusion now follows from the second assertion in
Lemma \ref{lem:OrthogImage}.
\end{pf}

\begin{lem}
\label{lem:RankOfTau}
\cite[Lemma 2.21]{FL1}
If $\Phi \in C^0(W^+\otimes E)$ and $x\in X$, then the following hold: 
\begin{enumerate}
\item $\Rank_\RR(\Phi(x)\otimes\Phi^*(x))_{00} = 1$ 
if and only if $\Rank_\CC\Phi(x) = 1$,
\item $\Rank_\RR(\Phi(x)\otimes\Phi^*(x))_{00} = 3$ 
if and only if $\Rank_\CC\Phi(x) = 2$.
\end{enumerate}
\end{lem}

\begin{rmk}
Note that if $\Phi$ is rank two on $X$, then $(\Phi\otimes\Phi^*)_{00}$ is
rank three on $X$ and so \eqref{eq:PT} implies that $F_A^+$ is rank three
on $X$ when $(A,\Phi)$ is a $\PU(2)$ monopole: thus, $A$ cannot be a
reducible connection on $X$.
\end{rmk}

We can now complete the proof of Theorem \ref{thm:SmoothParamModuliSpace}
in either of two ways; the interest in the second method is explained in
the remark at the end of the proof.

\begin{proof}[First proof of Theorem \ref{thm:SmoothParamModuliSpace}]
Since $\Rank_\CC\Phi(x)=2$, for all $x\in U$, Lemma
\ref{lem:RankOfTau} implies that
$\Rank_\RR(\Phi(x)\otimes\Phi^*(x))_{00} = 3$, for all $x\in U$, and so
Lemma \ref{lem:vandQuadOrthogImages} implies that $\Rank_\RR v(x)= 0$, for
all $x\in U$, since $\su(E)$ is rank three and thus $v \equiv 0$ on $U$.
Hence, $(v,\psi)\equiv 0$ on $U$ and so $(v,\psi)\equiv 0$ on $X$ by unique
continuation for the Laplacian $D\ubarfS D\ubarfS^*$, as desired.
\end{proof}

Alternatively, we observe that it is enough to know that
$(\Phi\otimes\Phi^*)_{00}$ has at least rank two at some point of $X$:

\begin{proof}[Second proof of Theorem \ref{thm:SmoothParamModuliSpace}]
Since $\Rank_\CC\Phi(x)=2$, for all $x\in U$, Lemma \ref{lem:RankOfTau}
implies that $\Rank_\RR(\Phi(x)\otimes\Phi^*(x))_{00} \geq 2$, for all
$x\in U$, and so Lemma \ref{lem:vandQuadOrthogImages} implies that
$\Rank_\RR v(x) \leq 1$, for all $x\in U$, since $\su(E)$ is rank three.
On an open subset $U'\subset X$ where $v\neq 0$, we can write
$v=\varpi\otimes b$, where $\varpi\in C^r(U',\La^+)$ and $b\in
C^r(U',\su(E))$, with $|b|=1$ on $U'$. The argument of \cite[Lemma
4.3.25]{DK} and \cite[Proposition 3.4, p. 56]{FU} now implies that $d_Ab=0$ on
$U'$. The connection $A$ is thus reducible on the non-empty open subset
$U'\subset X$, which we may assume to be connected without loss of
generality. Theorem \ref{thm:LocalToGlobalReducible} implies that $A$, and
so the pair $(A,\Phi)$, is reducible on all of $X$ --- again contradicting
our assumption that $(A,\Phi)$ is irreducible. Hence, $v\equiv 0$ on $X$
and so $(v,\psi)\equiv 0$ on $X$.
\end{proof}

\begin{rmk}
In conjunction with the variations of the Dirac operator considered in \S
\ref{sec:Period}, a direct approach (avoiding the period map of \S
\ref{sec:Period}) 
to the proof of our main transversality theorem would work were it not for
the fact that it seems impossible (using these variations) to show that $v$
has at most rank one on some non-empty open subset of $X$, as used in the
second proof above.
\end{rmk}


\section{Variation of the Dirac operator}
\label{sec:VarDirac}
The second main step in the proof of Theorem \ref{thm:Transversality} is to
show, for generic parameters, that if $(A,\Phi)$ is a $\PU(2)$ monopole
then the element $\Phi$ of the kernel of the 
perturbed Dirac operator defined by $A$
cannot be rank one. In \S \ref{subsec:CliffordDiracVar} we introduce our
full parameter space of perturbations of the Dirac operator and compute the
differential of the resulting family, while in \S
\ref{subsec:ParamModuliDirac} we describe the corresponding perturbed
$\PU(2)$ monopole equations and the universal moduli space. 

\subsection{Clifford maps and Dirac-operator variations}
\label{subsec:CliffordDiracVar}  
In order to describe our variations, it is helpful to have a construction of
the Dirac operator at our disposal which is as simple as possible. The
minimal, axiomatic approach employed by 
Kronheimer-Mrowka \cite{KMContact} and Mrowka-Ozsv\'ath-Yu
\cite{MrowkaOzsvathYu} is extremely useful for this purpose, so this is the
approach we shall follow here. 

Recall that a real-linear map
$\rho_+:T^*X\to\Hom_\CC(W^+,W^-)$ defines a Clifford-algebra
representation $\rho:\Cl_\CC(T^*X)\to\End_\CC(W)$, with $W := W^+\oplus W^-$, 
if and only if
\cite{LM}, \cite{SalamonSWBook}
\begin{equation}
\rho_+(\alpha)^\dagger\rho_+(\alpha) = g(\alpha,\alpha)\id_{W^+},
\qquad \alpha\in C^\8(T^*X),
\label{eq:MinimalCliffordMap}
\end{equation}
where $g$ denotes the Riemannian metric on $T^*X$.
The real-linear map $\rho:T^*X\to\End_\CC(W^+\oplus W^-)$ is obtained by
defining a real-linear map 
$$
\rho_-:T^*X\to\Hom_\CC(W^-,W^+),
\qquad
\alpha \mapsto \rho_-(\alpha) := -\rho_+(\alpha)^\dagger,
$$
and setting
\begin{equation}
\rho(\alpha) 
:=
\begin{pmatrix}
0 & \rho_-(\alpha) \\ \rho_+(\alpha) & 0
\end{pmatrix}.
\label{eq:DefnFullCliffordMap}
\end{equation}
The Clifford algebra representation
$\rho:\Cl_\CC(T^*X)\to\End_\CC(W^+\oplus W^-)$ is 
uniquely determined by the map $\rho_+:T^*X\to\Hom_\CC(W^+,W^-)$ and
satisfies 
\begin{equation}
\rho(\alpha)^\dagger = - \rho(\alpha)
\quad\text{and}\quad
\rho(\alpha)^\dagger\rho(\alpha) = g(\alpha,\alpha)\id_W,
\qquad \alpha\in C^\8(T^*X).
\label{eq:CliffordMap}
\end{equation}
We may vary $\rho_+$ by automorphisms of $T^*X$.
If $f \in C^\8(\GL(T^*X))$, then
\begin{align*}
\rho_+(f(\alpha))^\dagger \rho_+(f(\alpha))
&= g(f(\alpha),f(\alpha))\id_{W^+}
\\
&= (f^*g)(\alpha,\alpha)\id_{W^+},
\end{align*}
and we obtain a Clifford map $\rho_{f,+} := \rho_+\circ f
:T^*X\to\Hom_\CC(W^+,W^-)$ which is compatible with the Riemannian metric
$f^*g$ on $T^*X$. Thus, for $f \in C^\8(\Or(T^*X))$, the Clifford map
$\rho_{f,+}$ is compatible with the given Riemannian metric $g$.

Recall from \cite[p. 89]{MorganSWNotes} (or Lemma
\ref{lem:MorganLemma}), that the complexification of $\rho_+$
yields an isomorphism
\begin{equation}
\rho_+:(T^*X)\otimes_\RR\CC \to \Hom_\CC(W^+,W^-).
\label{eq:RhoIsomorphism}
\end{equation}
The space $C^r(\GL(T^*X))$ is a
Banach Lie group with Lie algebra $C^r(\gl(T^*X))$
\cite{FU}, \cite{PalaisFoundationGlobal}. Then, for all $e\in
\Omega^1(X,\RR)$ and 
$\Phi \in \Omega^0(W^+)$, we have
\begin{equation}
((D\rho_{f,+})(\delta f))(e)\Phi = \rho_+((\delta f)e)\Phi,
\label{eq:VarRhoAutTX}
\end{equation}
where $\delta f \in C^r(\gl(T^*X))$.

For the Riemannian metric $g$ on $T^*X$, let $\cov^g$ be an $\SO(4)$
connection on $T^*X$. A unitary connection $\cov$ on $W$ is called
{\em spinorial with respect to $\cov^g$\/} if it induces the $\SO(4)$
connection $\cov^g$ on $T^*X$, or, equivalently, if the covariant
derivative $\cov$ on $C^\8(W)$ is a derivation with respect to
Clifford multiplication, $\rho:C^\8(\Lambda^\bullet(T^*X)\otimes W)\to
C^\8(W)$, so \cite{LM}, \cite{SalamonSWBook}
\begin{equation}
\cov_\eta(\rho(\beta)\Phi) 
= 
\rho(\cov_\eta^g\beta)\Phi + \rho(\beta)\cov_\eta\Phi,
\label{eq:RhoCompatibility}
\end{equation}
for all $\eta \in C^\8(TX)$, $\beta \in \Omega^\bullet(X,\RR)$, and $\Phi
\in C^\8(W)$.  (Our
convention differs from that of \cite[Definition 4.0.23]{MrowkaOzsvathYu}.)
The unitary connection $\cov$ on $W$ uniquely determines a unitary
connection on $\det W^+
\simeq \det W^-$ in the standard way \cite{Kobayashi}. Any two unitary
connections on $W$, which are both spinorial with respect to $\cov^g$,
differ by an element of $\Omega^1(X,i\RR)$. Conversely, a unitary
connection $\cov$ on a Hermitian two-plane bundle $W = W^+\oplus W^-$
over an oriented four-manifold $X$ is uniquely determined by
\begin{itemize}
\item A Clifford map $\rho_+:T^*X\to\Hom_\CC(W^+,W^-)$ satisfying
\eqref{eq:MinimalCliffordMap} for the Riemannian metric $g$ on $T^*X$,
\item An $\SO(4)$ connection $\cov^g$ on $T^*X$ for the metric $g$, which
need not be torsion free, and, 
\item A $\U(1)$ connection $B$ on $\det W^+$.
\end{itemize}
The resulting connection $\cov$ on $W$ is then spinorial with respect to
$\cov^g$. 

Digressing slightly, we recall that the local connection matrix one-form of
$\cov$ may be expressed in terms of those of the connections on $T^*X$ and
$\det W^+$.  To see this, let $\{e^i\}$ be an oriented, $g$-orthonormal
local frame for $T^*X$ with dual frame $\{e_i\}$ for $TX$ and let
${}^j\omega_{kl}$, $j=1,2$, be the corresponding $\so(4)$ connection
matrices for the $\SO(4)$ connections ${}^j\cov^g$ over an open subset
$U\subset X$. Let $B \in
\Omega^1(U,i\RR)$ also denote the local connection one-form for the $\U(1)$
connection on $\det W^+$, with respect to a trivialization for $\det W^+$
induced from that of $W^+$ over $U$. {}From \cite[p. 4]{Hitchin} (see also 
\cite[Theorem II.4.14]{LM}), the local
connection matrix one-forms for the \spinc connections ${}^j\cov$ on $W$
defined by $(\rho,{}^j\cov^g,B)$ are given by
\begin{align*}
{}^j\omega^{\Spinc(4)}
&=
{}^j\omega^{\Spin(4)} + \half B\,\id_{W^+}
\\
&=
\quarter\sum_{k,l}{}^j\omega_{kl}\otimes\rho(e^k\wedge e^l)
+ \half B\,\id_{W^+}.
\end{align*}
Although the formula of \cite{Hitchin} refers only to spin connections,
locally we may write $W|_U = S\otimes_\CC N$, where $S$ is a spin bundle
for the local spin structure and $N$ is a Hermitian line
bundle such that $N^{\otimes 2} = \det
W^+|_U$. We may then write the \spinc connections on $W|_U$ as tensor
products of spin connections on $S$ with connection matrix one-forms
${}^j\omega^{\Spin(4)}$ and a $\U(1)$ connection
on $N$ with connection one-form $\half B$, as above.

Given a unitary connection $A$ on an auxiliary Hermitian two-plane bundle
$E$, we let $\cov_A$ denote the induced unitary connection on $W\otimes E$.
The corresponding Dirac operator \cite[\S 3.4]{Gilkey},
\begin{equation}
\Dirac_A := \rho_+\circ\cov_A,
\label{eq:DiracDefn}
\end{equation}
is defined by the composition of the covariant derivative,
$\cov_A:C^\8(W^+\otimes E)\to C^\8(T^*X\otimes W^+\otimes E)$, and Clifford
multiplication, $\rho_+:C^\8(T^*X\otimes W^+\otimes E) \to C^\8(W^-\otimes
E)$. If $\{v^i\}$ is a (not necessarily $g$-orthonormal) local frame for
$T^*X$ with dual frame $\{v_i\}$ for $TX$, defined by
$v^i(v_j)=\delta_{ij}$, then $\cov_A = \sum_{i=1}^4v^i\otimes\cov_{A,v_i}$
and so $\Dirac_A$ has the familiar shape
$$
\Dirac_A = \sum_{i=1}^4\rho_+(v^i)\cov_{A,v_i}.
$$
It is very convenient to keep the unitary
connection $\cov$ on $W$ {\em fixed\/} throughout, so while $(\cov,\rho)$
necessarily induces an $\SO(4)$ connection $\cov^g$ on $T^*X$, we shall
{\em not\/} require that $\cov^g$ be torsion free, that is, we shall not
assume $\cov^g$ is the Levi-Civita connection ${}^{\LC}\cov^g$ for the
metric $g$. However, the relation between the Dirac operators defined by
$(\rho,B,A)$ and two different $\SO(4)$ connections for the metric $g$ on
$T^*X$ is easily determined:

\begin{lem}
\label{lem:DiffDiracSO(4)}
Let $X$ be an oriented four-manifold with Riemannian metric $g$ on $T^*X$,
compatible Clifford map $\rho_+:T^*X\to\Hom_\CC(W^+,W^-)$, unitary
connection $B$ on $\det W^+$, and unitary connection $A$ on a Hermitian
two-plane bundle $E$ over $X$. If ${}^j\cov^g$, $j=1,2$ are two $\SO(4)$
connections on $T^*X$ for the metric $g$, and ${}^j\cov$ are the unitary
connections on $W$ induced by $(\rho,{}^j\cov^g,B)$, and ${}^j\Dirac_A$ are
the Dirac operators defined by $(\rho,{}^j\cov,A)$, then
$$
{}^2\Dirac_A - {}^1\Dirac_A
= 
\rho_+(\sigma)\otimes\id_E,
$$
where 
$$
\sigma := {}^2\cov - {}^1\cov \in \Omega^1(X,\su(W^+))
$$
and $\rho_+(\sigma)\in
\Hom_\CC(W^+,W^-)$ is defined by the contraction 
$$
\rho_+:C^\8(T^*X\otimes
\End_\CC(W^+))\to C^\8(\Hom_\CC(W^+,W^-)).
$$ 
\end{lem}

\begin{proof}
Since ${}^2\cov$ and ${}^1\cov$ are unitary connections on $W^+$, then
$\sigma = {}^2\cov - {}^1\cov \in \Omega^1(X,\fu(W^+))$, and as they
both induce the connection $B$ on $\det W^+$, then $\tr\sigma =
0$, so $\sigma \in \Omega^1(X,\su(W^+))$. If $\{v^i\}$ is any local frame for
$T^*X$ with dual frame $\{v_i\}$ for $TX$, defined by
$v^i(v_j)=\delta_{ij}$, then the difference between the Dirac operators is
given by
\begin{align*}
{}^2\Dirac_A - {}^1\Dirac_A
&=
\sum_{i=1}^4\rho_+(v^i)({}^2\cov_{A,v_i} - {}^1\cov_{A,v_i})
\\
&=
\sum_{i=1}^4\rho_+(v^i)\sigma(v_i)\otimes\id_E
=
\rho_+(\sigma)\otimes\id_E,
\end{align*}
with $\rho_+(\sigma) \in \Hom_\CC(W^+W^-)$, as desired. Alternatively, the
conclusion can be inferred from the expression for the local connection
matrix one-forms for the connections ${}^j\cov$ given in the 
paragraphs preceding the statement of the lemma.
\end{proof}

Lemma \ref{lem:DiffDiracSO(4)} allows us to write the Dirac
operator in \eqref{eq:PT}, defined by a Riemannian but not necessarily
torsion-free connection $\cov^g$ as the Dirac operator for ${}^{\LC}\cov^g$
plus an element $\rho_+(\theta)$ of $\Hom_\CC(W^+,W^-) =
\rho_+(T^*X\otimes\CC)$.  Since we already include a perturbation term
$\rho_+(\vartheta)$ in our definition of the $\PU(2)$ monopole equations,
we may therefore assume at the conclusion of our transversality argument
that the Dirac operator corresponds to the Levi-Civita connection for the
Riemannian metric on $T^*X$ by absorbing the correction term
$\rho_+(\theta)$ into $\rho_+(\vartheta)$. {\em For the remainder of the
article, however, the unitary connection on $W$ will be held fixed and so it
induces a Riemannian but not always torsion-free connection on $T^*X$
via \eqref{eq:RhoCompatibility} as the metric on $T^*X$ is allowed to vary.\/} 

{}From \eqref{eq:DiracDefn} we obtain a family of Dirac operators
parametrized by $C^r(\GL(T^*X))$. Specifically,
let $\Dirac_{A,f}$ be the Dirac operator defined by the following data:
\begin{itemize}
\item Clifford map
$\rho_+\circ f$ compatible with the metric $f^*g$ on $T^*X$, where $\rho_+$ is
a Clifford map compatible with the metric $g$, and $f\in C^\8(\GL(T^*X))$,
\item Unitary connection $\cov$ on $W$, 
\item Unitary connection $A$ on a Hermitian two-plane bundle $E$. 
\end{itemize}
If $\{v^i\}$ is an oriented local frame for $T^*X$ with
dual frame $\{v_i\}$ for $TX$ defined by $v^j(v_i) = \delta_{ij}$, then
$$
\Dirac_{A,f} = \sum_{i=1}^4\rho_+(f(v^i))\cov_{A,v_i}
$$
is the corresponding family of Dirac operators.

The variation of the Dirac operator $\Dirac_{A,f}$
induced from \eqref{eq:VarRhoAutTX} is given by
\begin{equation}
(D\Dirac_A)_f(\delta f)\Phi 
= \sum_{i=1}^4\rho_+((\delta f)v^i)\cov_{A,v_i}\Phi,
\label{eq:VarDiracAutTX}
\end{equation}
where $\Phi \in \Omega^0(W^+\otimes E)$. For a given Hermitian metric and
unitary connection on $W^+\oplus W^-$, an orientation of $X$, and a
Riemannian metric $g$ on $T^*X$, we then take our family of perturbed Dirac
operators
$$
\Dirac_{A,f,\vartheta} : \Omega^0(W^+\otimes E) \to \Omega^0(W^-\otimes E)
$$
to be
\begin{equation}
\Dirac_{A,f,\vartheta} 
:= \sum_{i=1}^4\rho_+(f(v^i))\cov_{A,v_i} + \rho_+(f(\vartheta)),
\label{eq:PerturbedDirac}
\end{equation}
for any unitary connection $A$ on $E$. The unitary connection $\cov$ on
$W^+\oplus W^-$ and Clifford map $\rho_{f,+}:T^*X\to\Hom_\CC(W^+,W^-)$
induces, via \eqref{eq:RhoCompatibility}, an $\SO(4)$ (but necessarily
torsion-free) connection on $T^*X$ for the metric $f^*g$.

\begin{lem}
\label{lem:VarDirac}
Continue the notation of this section. Then 
$$
(D\Dirac)_{A,f,\vartheta}(\delta A,\delta f,\delta\vartheta)\Phi 
= \sum_{i=1}^4\rho_+((\delta f)v^i)\cov_{A,v_i}\Phi
+ \rho_+(f(\delta A))\Phi + \rho_+(f(\delta\vartheta))\Phi.
$$
\end{lem}

We note that a direct variation of the Dirac operator with respect to the
Riemannian metric $g$ has been computed by Bourguignon and Gauduchon
\cite{BG}, using a rather different and more elaborate approach than the one
considered here. Their computation, with some enhancements, was recently
used by Maier to prove certain generic metrics results for Dirac operators
on low-dimensional spin or \spinc manifolds \cite{Maier}. It is not known
at present whether there exists a purely `generic metrics' transversality
result for the Seiberg-Witten equations (analogous to the celebrated
generic metrics theorem of Freed and Uhlenbeck
\cite{DK}, \cite{FU} for the anti-self-dual equation), although there have been
some attempts in this direction \cite{BryantTaubes}.

Recall that the space $\Met(X)$ of all $C^\8$ Riemannian metrics on $X$ is
a contractible, open cone inside the space $S^2(T^*X)$ of symmetric,
rank-two, contravariant tensor fields on $X$ and so, at any metric $g\in
\Met(X)$, it has tangent space $T_g\Met(X) = S^2(T^*X)$. The technique
employed by Bourguignon and Gauduchon constructs spinor bundles which
depend implicitly on the Riemannian metric (as is standard
\cite{LM}). Thus, in order to compute the differential of their family of Dirac
operators parametrized by the space of metrics, the family must be pulled
back to an equivalent family of operators acting on a {\em fixed\/} Hilbert
space of sections of a {\em fixed\/} \spinc bundle. This is the technique
employed by Bourguignon-Gauduchon and Maier; a similar one was used by the
authors in our proof of Theorem 5.11 in \cite{FL1}.  It is important to
note that the resulting family of operators in \cite{BG}
is not necessarily a family of
Dirac operators. In any event, their key result, namely \cite[Theorem
21]{BG} (see also \cite[Proposition 2.4]{Maier}), could be used in place of
the simpler variational formulas
\eqref{eq:VarDiracAutTX} which we shall employ
in the sequel. The principal difference for our application would be that
variations of $g$ alone do not span $\gl(TX)$, so we would need to consider
auxiliary variations of the Clifford map $\rho_+$ by elements $f\in
C^\8(\SO_g(T^*X))$.

\subsection{The parametrized moduli space of $\PU(2)$ monopoles}
\label{subsec:ParamModuliDirac}
Although the curvature equation in \eqref{eq:PT} will only play an
auxiliary role in the proof of Theorem \ref{thm:NoRankOneIrreducibles}
(for example, through the local-to-global unique continuation result, Theorem
\ref{thm:LocalToGlobalReducible}, for reducible $\PU(2)$ monopoles), we
will need a parametrized moduli space as a domain of definition for our
infinite-dimensional period map. In this section we define $\PU(2)$
monopole equations containing the full set of perturbations
$(f,\tau,\vartheta)$ and define the corresponding parametrized moduli
space.

We fix a $C^\8$ Riemannian metric $g$ on $T^*X$ and set
\begin{equation}
\sP := C^r(\GL(T^*X))\oplus C^r(\Lambda^1\otimes\CC).
\label{eq:DiracParamSpace}
\end{equation}
As we remarked in \S \ref{subsec:TopStratumTransv}, we could of course
include the space $C^r(\GL(\Lambda^+))$ of parameters $\tau$ in our
definition of $\sP$ in \eqref{eq:DiracParamSpace}; {\em we simply omit them
for convenience here, as they play no role in \S \ref{sec:VarDirac} or \S
\ref{sec:Period}\/}, just as we omitted the parameters $f$ from the
definition of $\sP$ in \S \ref{sec:Ranktwo}.

For $f\in C^r(\GL(T^*X))$,
we let $\rho_f:T^*X\to\End_\CC(W^+\oplus W^-)$ be the Clifford map defined
by $\rho_{f,+}$ and \eqref{eq:DefnFullCliffordMap}. The map $\rho_f$ 
is compatible with the Riemannian metric $f^*g$ on $T^*X$ 
in the sense of \eqref{eq:CliffordMap} and thus
extends, in the usual way, to give an isomorphism of real three-plane
bundles,
\begin{equation}
\rho_f:\Lambda^{+,f^*g} \to \su(W^+).
\label{eq:ParamLambda+SuW+Isom}
\end{equation}
Let $\{e^i\}$ be an oriented, $g$-orthonormal local frame for
$T^*X$ and observe that $\{f(e^i)\}$ is an oriented, $f^*g$-orthonormal
local frame for $T^*X$, so we have induced isomorphisms
$$
f^{-1}:\Lambda^\pm\to\Lambda^{\pm,f^*g},
\qquad
\omega\mapsto f^{-1}(\omega),
$$
taking $e^1\wedge e^2+e^3\wedge e^4 \mapsto 
f(e^1)\wedge f(e^2)+f(e^3)\wedge f(e^4)$, and similarly for the remaining
two elements of the standard local frame for $\Lambda^+ := \Lambda^{+,g}$
and the three elements of the local frame for $\Lambda^- := \Lambda^{-,g}$.
Note that if $P_+(g) := \half(1 +*_g):\Lambda^2=\Lambda^+\oplus \Lambda^-\to
\Lambda^+$ is the projection onto $g$-self-dual two-forms, then
$$
P_+(f^*g) = f^{-1}\circ P_+(g)\circ f : \Lambda^2 \to \Lambda^{+,f^*g}
$$
is the projection onto $f^*g$-self-dual two-forms. (We find it
convenient to use automorphisms $f$ of $T^*X$ to vary the metric $g$ on
$T^*X$, while Freed-Uhlenbeck \cite[pp. 51--52]{FU} 
vary the metric on $TX$ by automorphisms of $TX$.)

Suppose $\omega := \sum_{i<j}\om_{ij}e^i\wedge e^j \in
\Omega^2(X,\RR)$. For any $\Phi\in C^\8(W^+)$, we see that
\begin{align*}
\rho_f(\omega)\Phi 
&=
\sum_{i<j}\om_{ij}\rho_f(e^i)\rho_f(e^j)\Phi
=
\sum_{i<j}\om_{ij}\rho(f(e^i))\rho(f(e^j))\Phi
\\
&=
\sum_{i<j}\om_{ij}\rho(f(e^i)\wedge f(e^j))\Phi
=
\rho\left(\sum_{i<j}\om_{ij}f(e^i\wedge e^j)\right)\Phi
\\
&=
\rho(f(\omega))\Phi = (\rho\circ f)(\omega)\Phi.
\end{align*}
Conversely, since $(\Phi\otimes\Phi^*)_{00} \in \Omega^0(\su(W^+))$,
$\tau\in \Omega^0(\GL(\Lambda^+))$,
and $\tau_f := f^{-1}\circ\tau\circ f \in \Omega^0(\GL(\Lambda^{+,f^*g}))$, we
have 
$$
\tau_f\rho_f^{-1}(\Phi\otimes\Phi^*)_{00}
=
f^{-1}\tau\rho^{-1}(\Phi\otimes\Phi^*)_{00}
\in \Omega^0(\Lambda^{+,f^*g}).
$$
For the metric $f^*g$ on $T^*X$,
compatible Clifford map $\rho_f:\Lambda^{+,f^*g}\to \su(W^+)$, 
and compatible automorphism $\tau_f$ of $\Lambda^{+,f^*g}$, the curvature
equation in \eqref{eq:PT} is given by
$$
P_+(f^*g)F_A -\tau_f\rho_f^{-1}(\Phi\otimes\Phi^*)_{00} = 0
\in \Omega^{+,f^*g}(\su(E)),
$$
that is,
$$
(f^{-1}\circ P_+(g)\circ f)(F_A) -
f^{-1}\rho^{-1}\tau(\Phi\otimes\Phi^*)_{00} 
= 0,
$$
or equivalently,
$$
P_+(g)f(F_A) - \rho^{-1}\tau(\Phi\otimes\Phi^*)_{00} 
= 0 \in \Omega^+(\su(E)).
$$
We now define our $\ssG_E$-equivariant map 
$$
\ubarfS:\sP\times \tsC_{W,E}
\to L^2_{k-1}(\La^+\otimes(\su(E))\oplus L^2_{k-1}(W^-\otimes E))
$$
by setting
\begin{equation}
\ubarfS(f,\vartheta,A,\Phi) 
:= 
\begin{pmatrix}
P_+(g)f(F_A) -\tau\rho^{-1}(\Phi\otimes\Phi^*)_{00} \\
\Dirac_{A,f}\Phi+\rho(f(\vartheta))\Phi
\end{pmatrix}.
\label{eq:DiracParamMonopole}
\end{equation}
The group $\ssG_E$ acts trivially on the space of perturbations $\sP$ and
so $\ubarfS^{-1}(0)/\ssG_E$ is a subset of $\sP\times\sC_{W,E}$.  We now
let $\fM_{W,E}$ denote the parametrized (or universal) moduli space
$\ubarfS^{-1}(0)/\ssG_E$ for the augmented space $\sP$ of perturbation
parameters and let $\fM_{W,E}^{*,0}:=\fM_{W,E}\cap (\sP\times
\sC_{W,E}^{*,0})$.


\section{Period maps for Dirac operators}
\label{sec:Period}
In \S \ref{sec:Ranktwo} we proved that the moduli space of $\PU(2)$
monopoles is cut out transversely away from the loci of pairs $(A,\Phi)$
with either $A$ reducible or $\Phi$ rank less than or equal to one.  In
this section we prove the following 
analogue of Proposition 4.3.14 and Corollary
4.3.15 in \cite{DK} for the locus of pairs $(A,\Phi)$ with $A$ irreducible
and $\Phi$ rank one but not identically zero: 

\begin{thm}
\label{thm:NoRankOneIrreducibles}
Let $X$ be a closed, oriented, smooth four-manifold with \spinc structure
$(\rho,W^+,W^-)$ and Hermitian two-plane bundle $E$. Then there is a
first-category subset $\sP^\8_{\fc} \subset \sP^\8$ such that for all
$(f,\tau,\vartheta)$ in $\sP^\8 - \sP^\8_{\fc}$ the following holds:
The moduli space $M_{W,E}^{*,0}(f,\tau,\vartheta)$ contains no $\PU(2)$
monopoles $(A,\Phi)$ with both $A$ irreducible and $\Phi$ rank one. 
\end{thm}

The argument we describe here only works for $\PU(2)$ monopoles $(A,\Phi)$
with $A$ irreducible. Thus, even for generic parameters
$(f,\tau,\vartheta)$, the moduli space $M_{W,E}(f,\tau,\vartheta)$ will in
general contain points $[A,\Phi]$ with $A$ reducible and $\Phi$ rank one
and these will not necessarily be smooth points of
$M_{W,E}(f,\tau,\vartheta)$. However, it is reassuring to note that the
loci of reducible $\PU(2)$ monopoles --- corresponding to moduli spaces of
Seiberg-Witten monopoles --- cannot be perturbed away by the argument we
present here: the fact that the $\PU(2)$ monopole connections $A$ are
irreducible is used in an essential way in the proof of the key Proposition
\ref{prop:CokerDifferentialPeriod}.

Given Theorem \ref{thm:NoRankOneIrreducibles}, we can quickly dispose of
the proof of Theorem \ref{thm:Transversality}.

\begin{proof}[Proof of Theorem \ref{thm:Transversality}, given Theorem
\ref{thm:NoRankOneIrreducibles}]
As we remarked at the beginning of \S \ref{subsec:TopStratumTransv} and \S
\ref{subsec:ParamModuliDirac}, we may assume without loss of generality
that the same Fr\'echet space of $C^\8$ perturbation parameters,
$$
\sP^\8
:=
\Omega^0(\GL(T^*X)) \times \Omega^0(\GL(\Lambda^+)) 
\times \Omega^0(\Lambda^1\otimes\CC),
$$
has been used for the proofs of both Theorem \ref{thm:GenericProjection}
and Theorem \ref{thm:NoRankOneIrreducibles}.
By Theorem \ref{thm:GenericProjection},
there is a first-category subset $\sP_{\fc}' \subset \sP$ such that for all
$p = (f,\tau,\vartheta)\in \sP-\sP_{\fc}'$, the moduli space
$M_{W,E}^{*,\natural}(p)$ of irreducible, rank-two $\PU(2)$ monopoles is a
smooth manifold of the expected dimension. On the other hand, by Theorem
\ref{thm:NoRankOneIrreducibles}, there is a
first-category subset $\sP_{\fc}'' \subset
\sP$ such that for all $p = (f,\tau,\vartheta)\in \sP-\sP_{\fc}''$, 
the moduli space $M_{W,E}^{*,0}(p)$ contains no 
irreducible, rank-one $\PU(2)$
monopoles, that is, $M_{W,E}^{*,\natural}(p) = M_{W,E}^{*,0}(p)$.
Hence, for all $p \in \sP-\sP_{\fc}$, where $\sP_{\fc} := \sP_{\fc}'\cup
\sP_{\fc}''$, the moduli space $M_{W,E}^{*,0}(p)$ is a smooth manifold of
the expected dimension.
\end{proof}

The remainder of \S 4 is taken up with the proof of Theorem
\ref{thm:NoRankOneIrreducibles}. 

\subsection{Infinite-dimensional Grassmann manifolds}
\label{subsec:Grassmann}
The period map described in \cite[\S 4.3.3 \& 4.3.5]{DK}
takes values in the Grassmann manifold
$\GG(H^2(X;\RR))$ of $b^-(X)$-dimensional subspaces of $H^2(X;\RR)$. We shall
need to consider an infinite-dimensional version of this
construction, so it is convenient at this point to recall some properties
of infinite-dimensional Grassmann manifolds \cite[\S 3.1.8]{AMR},
\cite[\S 1.1]{Koschorke}. 

If $\bE$ is a complex Banach space then the Grassmann manifold $\GG(\bE)$
is the set of splitting linear subspaces of $\bE$ and is a complex-analytic
manifold with tangent spaces
$$
T_K\GG(\bE) = \Hom_\CC(K, \bE/K).
$$
Recall that a closed, linear subspace $K\subset\bE$ {\em splits\/} if $K$
has a closed complement $K'\subset\bE$ such that $\bE\simeq K\times
K'$ and that any finite-dimensional subspace splits. For each integer
$\kappa\ge 1$, let
$$
\GG_\kappa(\bE) 
:= 
\{K\subset\bE: \text{$K$ is a $\kappa$-dimensional subspace of $\bE$}\},
$$
and set $\PP(\bE) = \GG_1(\bE)$. The spaces $\GG_\kappa(\bE)$ are connected
submanifolds of $\GG(\bE)$.

The Grassmann manifolds $\PP(\bE)$ and $\GG_\kappa(\bE)$ define a locus of
incident planes (or flag manifold) $\FF_\kappa(\bE)
\subset \PP(\bE)\times\GG_\kappa(\bE)$ by analogy with the usual
construction in the finite-dimensional case \cite{Harris}:
$$
\FF_\kappa(\bE) := \{(\ell,K) \in \PP(\bE)\times\GG_\kappa(\bE):
\ell \subset K\}.
$$
The locus $\FF_\kappa(\bE)$ is a smooth manifold with tangent spaces 
\begin{align}
\label{eq:TangentSpaceFlag}
&T_{(\ell,K)}\FF_\kappa(\bE)
\\
&\quad =
\{(\eta,\varphi) \in \Hom_\CC(\ell,\bE/\ell) \oplus \Hom_\CC(K,\bE/K):
\varphi|_\ell = \eta\pmod{K}\}.
\notag
\end{align}
Given a smooth submanifold $\fX \subset \PP(\bE)$, we define a locus
$\II_\kappa(\fX) \subset \GG_\kappa(\bE)$ of incident $\kappa$-planes by
setting 
$$
\II_\kappa(\fX) 
:= 
\pi_2(\pi_1^{-1}(\fX)) \subset \GG_\kappa(\bE),
$$
where $\pi_1$, $\pi_2$ are the projections from the flag manifold onto the
first and second factors:
$$
\begin{CD}
          @.                     @. \FF_\kappa(\bE)                     \\
          @.   {\pi_1}\swarrow   @.                @. \searrow{\pi_2}   \\
\PP(\bE)  @.                     @.                @.@. \GG_\kappa(\bE) \\
\end{CD}
$$
We then have the following straightforward observations:

\begin{claim}
The map 
\begin{equation}
\pi_1:\FF_\kappa(\bE) \to \PP(\bE),
\qquad
(\ell,K)\mapsto \ell,
\label{eq:ProjFlagToProjSpace}
\end{equation}
is a submersion. 
\end{claim}

\begin{proof}
{}From the description of $T_{(\ell,K)}\FF_\kappa(\bE)$ in
\eqref{eq:TangentSpaceFlag} and of $T_\ell\PP(\bE)$ as
$\Hom_\CC(\ell,\bE/\ell)$, we see that for any $(\ell,K)\in
\FF_\kappa(\bE)$, the differential (which is again projection
onto the first factor),
$$
(D\pi_1)_{(\ell,K)}:T_{(\ell,K)}\FF_\kappa(\bE) \to T_\ell\PP(\bE),
\qquad
(\eta,\varphi)\mapsto \eta,
$$
is surjective. Indeed, if $\eta\in \Hom_\CC(\ell,\bE/\ell)$,
then we can choose $\varphi\in \Hom_\CC(K,\bE/K)$ such that
$\varphi(\ell)=\eta$ and $\varphi$ is defined arbitrarily on $K/\ell$, so 
$(D\pi_1)_{(\ell,K)}(\eta,\varphi) = \eta$.
\end{proof}

Hence, the locus $\tII_\kappa(\fX) := \pi_1^{-1}(\fX) \subset
\FF_\kappa(\bE)$ is a smooth submanifold with tangent spaces
$$
T_{(\ell,K)}\tII_\kappa(\fX)
=
\{(\eta,\varphi) \in T_\ell\fX \oplus \Hom_\CC(K,\bE/K):
\varphi|_\ell = \eta\pmod{K}\},
$$
where we view $T_\ell\fX$ as a subspace of $\Hom_\CC(\ell,\bE/\ell)$, and
codimension
$$
\codim(\tII_\kappa(\fX);\FF_\kappa(\bE)) = \codim(\fX;\PP(\bE)).
$$ 
In the same vein, we have:

\begin{claim}
The map 
\begin{equation}
\pi_2:\FF_\kappa(\bE) \to \GG_\kappa(\bE),
\qquad
(\ell,K)\mapsto K,
\label{eq:ProjFlagToGrassmann}
\end{equation}
is a submersion.
\end{claim}

\begin{proof}
We see that the differential (which is again projection
onto the second factor)
$$
(D\pi_2)_{(\ell,K)}:T_{(\ell,K)}\FF_\kappa(\bE) \to T_K\GG_\kappa(\bE)
\qquad
(\eta,\varphi)\mapsto \varphi,
$$
is surjective from the description of $T_{(\ell,K)}\FF_\kappa(\bE)$ 
in \eqref{eq:TangentSpaceFlag} and of
$T_K\GG_\kappa(\bE)$ as $\Hom_\CC(K,\bE/K)$. 
Indeed, if $\varphi\in \Hom_\CC(K,\bE/K)$,
then we can choose $\eta\in \Hom_\CC(\ell,\bE/\ell)$ by setting
$\eta=\varphi|_\ell$, so $(D\pi_2)_{(\ell,K)}(\eta,\varphi) = \varphi$.
\end{proof}

Over each point $K\in
\GG_\kappa(\bE)$, the projection \eqref{eq:ProjFlagToGrassmann}
has fibers 
$$
\pi_2^{-1}(K) 
= 
\{\ell\in\PP(\bE):\ell\subset K\}
=
\PP(K),
$$
given by finite-dimensional, complex projective spaces $\PP(K) \simeq
\PP^{\kappa-1}$.

The image $\II_\kappa(\fX)\subset\GG_\kappa(\bE)$ of
$\tII_\kappa(\fX)\subset \FF_\kappa(\bE)$ under the projection map $\pi_2:
\PP(\bE)\times\GG_\kappa(\bE)\to\GG_\kappa(\bE)$ is not necessarily a
smooth submanifold. For this reason, we restrict our attention henceforth
to the smooth submanifold $\tII_\kappa(\fX)\subset \FF_\kappa(\bE)$.

\begin{rmk}
If the Banach space $\bE$ were finite-dimensional, so $\bE = \CC^{n+1}$,
$\PP(\bE) = \PP^n$, and $\GG_\kappa(\bE) = \GG(\kappa,n+1)$, and
$\fX\subset \PP^n$ were a complex projective variety, then the incidence
locus $\II_\kappa(\fX)\subset
\GG(\kappa,n+1)$ would be a complex projective subvariety with smooth locus
$\II_\kappa(\fX)_{\sm}$ of codimension \cite[Example 16.6]{Harris}
$$
\codim(\II_\kappa(\fX)_{\sm};\GG(\kappa,n+1)) 
= \codim(\fX_{\sm};\PP^n) - (\kappa-1),
$$
where $\II_\kappa(\fX)$ has a standard stratification by virtue of its
status as a complex projective variety. 
When $\bE$ is an infinite-dimensional Banach space, as in our
application, it would suffice for our purposes to show that the image
$\II_\kappa(\fX)$ of $\tII_\kappa(\fX) = \pi_1^{-1}(\fX)$ under the
projection to $\GG_\kappa(\bE)$, 
is contained in a countable union of smooth submanifolds of
$\GG_\kappa(\bE)$ with infinite codimension. 
\end{rmk}

\subsection{Spaces of Fredholm operators}
\label{subsec:Fredholm}
For a non-negative integer $n$, let $\Fred_n(\bE^+,\bE^-)$ denote the open
subset of the Banach space of bounded operators, $\Hom_\CC(\bE^+,\bE^-)$,
consisting of bounded, index-$n$ Fredholm operators, where $\bE^+ :=
L^2_k(W^+\otimes E)$ and $\bE^- := L^2_{k-1}(W^-\otimes E)$
\cite{Koschorke}.  While the direct analogue of the period map of
\cite{DonConn}, \cite{DK}, \cite{FU} 
would use the infinite-dimensional Grassmann manifold $\GG(\bE^+)$ as
a target space, this is somewhat less suitable for our purposes since
the dimension of $\Ker\Dirac_{A,f,\vartheta}$ might jump due to
spectral flow as the point $[A,\Phi,f,\vartheta]$ varies in
$\sC_{W,E}^{*,0}\times\sP$; consequently, the assignment
$(A,\Phi,f,\vartheta)\mapsto \Ker\Dirac_{A,f,\vartheta}$ might not
define a smooth map from $\tsC_{W,E}^{*,0}\times\sP$ to
$\GG(\bE^+)$. However, the assignment $(A,\Phi,f,\vartheta)\mapsto
\Ker\Dirac_{A,f,\vartheta}$ does give a smooth map from
$\tsC_{W,E}^{*,0}\times\sP$ to $\Fred_n(\bE^+,\bE^-)$. Indeed, a map
of this form is used by Maier
\cite{Maier} (with domain $\Met(X)$, the space of all $C^r$ Riemannian
metrics on $T^*X$).

Passing temporarily to a more general setting, let $\bE^+$, $\bE^-$ be
complex Hilbert spaces. The subsets
$$
\Fred_{\kappa,n}(\bE^+,\bE^-)
:=
\{B\in \Fred_n(\bE^+,\bE^-): \dim_\CC\Ker B = \kappa\}
$$
are locally-closed submanifolds of $\Hom_\CC(\bE^+,\bE^-)$. 
The normal bundle $N\Fred_{\kappa,n}(\bE^+,\bE^-)$
of the submanifold $\Fred_{\kappa,n}(\bE^+,\bE^-)$
relative to $\Fred_n(\bE^+,\bE^-)$, and so $\Hom_\CC(\bE^+,\bE^-)$, has fiber
$$
N_B\Fred_{\kappa,n}(\bE^+,\bE^-)
=
\Hom_\CC(\Ker B,\Coker B)
$$
over a point $B \in \Fred_{\kappa,n}(\bE^+,\bE^-)$.
The submanifold $\Fred_{\kappa,n}(\bE^+,\bE^-)$ thus has complex codimension
$\kappa(\kappa-n)$ in $\Fred_n(\bE^+,\bE^-)$ and we have a smooth
stratification 
$$
\Fred_n(\bE^+,\bE^-) = \bigcup_{\kappa\ge n}\Fred_{\kappa,n}(\bE^+,\bE^-),
$$
with top stratum $\Fred_{n,n}(\bE^+,\bE^-)$. 

Associated with each $\Fred_{\kappa,n}(\bE^+,\bE^-)$, we have a
`flag manifold'
\begin{align*}
\Flag_{\kappa,n}(\bE^+,\bE^-)
:&=
\{(\ell,B) \in \PP(\bE^+)\times\Fred_{\kappa,n}(\bE^+,\bE^-):
\ell \in \Ker B\}
\\
&\subset
\PP(\bE^+)\times\Fred_{\kappa,n}(\bE^+,\bE^-).
\end{align*}
The corresponding projection
\begin{equation}
\pi:\Flag_{\kappa,n}(\bE^+,\bE^-)
\to
\Fred_{\kappa,n}(\bE^+,\bE^-),
\qquad
(\ell,B)\mapsto B,
\label{eq:FlagToFred}
\end{equation}
has fiber over a point $B \in \Fred_{\kappa,n}(\bE^+,\bE^-)$,
$$
\pi^{-1}(B) 
=
\{\ell\in\PP(\bE^+):\ell\subset \Ker B\}
=
\PP(\Ker B),
$$
given by a finite-dimensional, complex projective space $\PP(\Ker B) \simeq
\PP^{\kappa-1}$.

Furthermore, we have smooth maps
\begin{equation}
\pi:\Fred_{\kappa,n}(\bE^+,\bE^-) \to \GG_\kappa(\bE^+),
\qquad
B \mapsto \Ker B,
\label{eq:FredToGrassmann}
\end{equation}
which are easily seen to be submersions:

\begin{lem}
\label{lem:FredToGrassmannIsSubmersion}
The canonical map $\pi:\Fred_{\kappa,n}(\bE^+,\bE^-) \to \GG_\kappa(\bE^+)$
is a submersion.
\end{lem}

\begin{proof}
We first compute the differential $(D\pi)_B$ 
of $\pi$ at a point $B\in
\Fred_{\kappa,n}(\bE^+,\bE^-)$. Let $B_t$ be a smooth path in
$\Fred_{\kappa,n}(\bE^+,\bE^-)$ through $B_0 := B$. We then have a smooth
path of operators $B_t^\dagger B_t \in \End_\CC(\bE^+)$ giving isomorphisms
$B_t^\dagger B_t: (\Ker B_t)^\perp \to (\Ker B_t)^\perp$, where $\Ran
B_t^\dagger = (\Ker B_t)^\perp$ since $B_t$ is Fredholm and thus has closed
range. Let
$$
G_t := (B_t^\dagger B_t)^{-1}B_t^\dagger \in \Hom_\CC(\bE^-,\bE^+)
$$
be the corresponding smooth path of Green's operators and let 
$$
\Pi_t := \id_{\bE^+} - G_t B_t \in \End_\CC(\bE^+)
$$
be the resulting smooth path of projections from $\bE^+$ onto $\Ker
B_t$. Differentiating this path yields
\begin{align*}
\frac{d\Pi_t}{dt}
&=
-\frac{dG_t}{dt} B_t - G_t\frac{dB_t}{dt}
\\
&=
G_t\frac{dB_t}{dt} G_tB_t - G_t\frac{dB_t}{dt}.
\end{align*}
So, if $\delta B := (d B_t/dt)|_{t=0}$ and $G := G_0$ and $\Pi := \Pi_0$,
we have 
\begin{align*}
(D\pi)_B
&=
\left.\frac{d\Pi_t}{dt}\right|_{t=0}
=
G(\delta B)(GB-\id_{\bE^+})
\\
&=
-G(\delta B)\Pi
\in 
\Hom_\CC(\Ker B,(\Ker B)^\perp),
\end{align*}
with $\delta B \in T_B\Fred_{\kappa,n}(\bE^+,\bE^-)$. Any
operator $A \in \Hom_\CC(\bE^+,\bE^-)$ may be decomposed as a block-matrix
$$
A = \begin{pmatrix} A_{11} & A_{12} \\ A_{21} & A_{22} \end{pmatrix}:
(\Ker B)^\perp\oplus \Ker B \to \Ran B\oplus (\Ran B)^\perp,
$$
with $A_{22} \in \Hom_\CC(\Ker B,(\Ran B)^\perp)$, the normal space to 
$T_B\Fred_{\kappa,n}(\bE^+,\bE^-)$ in $\Hom_\CC(\bE^+,\bE^-)$. Then,
$$
A_{12} \in \Hom_\CC(\Ker B,\Ran B) \subset
T_B\Fred_{\kappa,n}(\bE^+,\bE^-),
$$
where $A_{12}$ is regarded as an operator in $\Hom_\CC(\bE^+,\bE^-)$
by extending by zero, so $A_{12} := 0$ on $(\Ker B)^\perp$.
Indeed, $A_{12}$ is a compact operator (since $\dim \Ker B < \8$, so
$A_{12}$ is finite rank) and thus $\Ind(B+A_{12}) = \Ind B$. Moreover,
$\Ran(B + A_{12}) = \Ran B$, so $\Coker(B + A_{12}) = \Coker B$ and thus 
$\dim\Coker(B + A_{12}) = \dim\Coker B$ implies 
$\dim\Ker(B + A_{12}) = \dim\Ker B = \kappa$. Hence, $B_t := B + tA_{12} \in
\Fred_{\kappa,n}(\bE^+,\bE^-)$ for all $t\in\RR$ and, in particular,
$(dB_t/dt)|_{t=0} = A_{12} \in T_B\Fred_{\kappa,n}(\bE^+,\bE^-)$. Since
$G:\Ran B\to (\Ker B)^\perp$ is an isomorphism and we may choose any
$\delta B \in \Hom_\CC(\Ker B,\Ran B)$, we see that
\begin{align*}
\Ran(D\pi)_B
&=
\{G(\delta B)\Pi: \delta B \in T_B\Fred_{\kappa,n}(\bE^+,\bE^-)\}
\\
&=
\Hom_\CC(\Ker B,(\Ker B)^\perp)
=
T_{\Ker B}\GG(\bE^+),
\end{align*}
and so $(D\pi)_B$ is surjective, as desired.
\end{proof}

Consequently, we have:

\begin{lem}
\label{lem:FredFlagToGrassmannFlagSubmersion}
The space $\Flag_{\kappa,n}(\bE^+,\bE^-)$ is a smooth submanifold 
of $\PP(\bE^+)\times \Fred_{\kappa,n}(\bE^+,\bE^-)$ and the
canonical map $\pi:\Flag_{\kappa,n}(\bE^+,\bE^-) \to \FF_\kappa(\bE^+)$
is a submersion.
\end{lem}

\begin{proof}
Lemma \ref{lem:FredToGrassmannIsSubmersion} implies that the projection,
$$
\pi:\PP(\bE^+)\times\Fred_{\kappa,n}(\bE^+,\bE^-)
\to
\PP(\bE^+)\times\GG_\kappa(\bE^+),
\qquad
(\ell,B)\mapsto (\ell,\Ker B),
$$
is a submersion. So, as $\FF_\kappa(\bE^+)$ is a smooth submanifold of 
$\PP(\bE^+)\times\GG_\kappa(\bE^+)$, we see that the preimage
$$
\Flag_{\kappa,n}(\bE^+,\bE^-)
=
\pi^{-1}(\FF_\kappa(\bE^+))
$$
is a smooth submanifold of $\PP(\bE^+)\times
\Fred_{\kappa,n}(\bE^+,\bE^-)$ and that the restriction of the projection
map is a submersion.
\end{proof}

\subsection{The locus of rank-one sections}
\label{subsec:RankOneLoci}
With the general setting of \S \ref{subsec:Grassmann} and \S
\ref{subsec:Fredholm} at hand, we can now describe the locus of rank-one
sections in $L^2_{k-1}(W^+\otimes E)$ and the Dirac-operator period map
which is to be transverse, in a suitable sense, to this locus.

Assume $k\ge 4$, so that $\bE^+ = L^2_{k-1}(W^+\otimes E)$ is contained in 
$C^0(W^+\otimes E)$. Since $C^0(W^+\otimes E) = C^0(\Hom_\CC(W^{+,*},E))$,
we have a determinant map
\begin{equation}
\det: C^0(W^+\otimes E) \to C^0(\det E\otimes\det W^+),
\qquad
\Phi\mapsto \det\Phi,
\label{eq:DetMap}
\end{equation}
recalling that $\det W^+ = \Lambda^2(W^+)$ and $\det E =
\Lambda^2(E)$. If $\{\xi_1,\xi_2\}$ is a local unitary frame for $E$ and
$\{\phi_1,\phi_2\}$ is a local unitary frame for $W^+$, with dual coframe
$\{\phi_1^*,\phi_2^*\}$ for $W^{+,*}$ defined by the Hermitian metric
via $\phi_j^* := \langle\cdot,\phi_j\rangle$, then
the section $\det\Phi$ of $\det E\otimes\det W^+
\simeq \Hom_\CC((\det W^+)^*,\det E)$ is
defined in the usual way with respect to these frames by 
$$
\Phi(\phi_1^*)\wedge\Phi(\phi_2^*)
=
(\det\Phi)\phi_1^*\wedge\phi_2^* \in \CC\,\xi_1\wedge\xi_2.
$$
A section $\Phi \in C^0(W^+\otimes E)$ then has {\em rank one\/} if
$\det\Phi\equiv 0$ and $\Phi\not\equiv 0$ on $X$. Locally, any rank-one
section $\Phi \in C^0(W^+\otimes E)$ may be written as $\Phi =
\phi\otimes\xi$, for local sections $\phi\in C^0(W^+)$ and $\xi\in C^0(E)$,
though it is generally not possible to write a rank-one section globally in
this form. 

We can now conveniently define the `rank-one locus' $\fX\subset\PP(\bE^+)$ by
\begin{equation}
\fX
:=
\{[\Phi]\in\PP(\bE^+): \det\Phi = 0\},
\label{eq:RankOneLocusInProjSpace}
\end{equation}
where $[\Phi]$ denotes the line $\CC\cdot\Phi\subset\bE^+$. Unless we impose
further conditions on the Riemannian four-manifold $X$ and the bundles
$W^+$ and $E$ --- such as requiring them to be real analytic --- the locus
of rank-one sections will not be a smooth subvariety of
$\PP(\bE^+)$. However, as we shall see in the sequel, it suffices to work
with the locus $\fX_{\sm} \subset \fX$ of smooth points of $\fX$. The
following lemma provides a simple criterion for  point $[\Phi]\in\fX$ to be
smooth.

\begin{lem}
\label{lem:SmoothProjVariety}
Suppose $[\Phi]\in \PP(\bE^+)$ is point in the zero locus $\fX =
\det^{-1}(0)$ such that $\{\Phi\neq 0\}$ is a dense open subset
of $X$. Then the map $\det:C^0(E\otimes W^+)\to C^0(\det E\otimes\det
W^+)$ vanishes transversely at $\Phi$ and $[\Phi]$ is a smooth point of
$\fX$. The tangent space $T_{[\Phi]}\fX$ has both infinite dimension and
infinite codimension in $T_{[\Phi]}\PP(\bE^+)$, so
$$
\codim(\fX_{\sm};\PP(\bE^+)) = \8.
$$
\end{lem}

\begin{proof}
Choose a local unitary frame $\{\xi_1,\xi_2\}$ for $E$ and a local unitary
frame $\{\phi_1,\phi_2\}$ for $W^+$, with dual coframe
$\{\phi_1^*,\phi_2^*\}$ for $W^{+,*}$ defined by the Hermitian metric on
$W^+$ via $\phi_j^* := \langle\cdot,\phi_j\rangle$.
With respect to these local frames over an open subset $U\subset X$, the
section $\Phi$ and determinant map are represented by
\begin{align*}
&\det:C^\8(U,\gl(2,\CC)) \to C^\8(U,\CC),
\\
&\quad
\Phi 
= 
\begin{pmatrix}\varphi_{11} & \varphi_{12} \\ 
\varphi_{21} & \varphi_{22} \end{pmatrix}
\mapsto 
\det\Phi = \varphi_{11}\varphi_{22}-\varphi_{12}\varphi_{21},
\end{align*}
with differential
$$
(D\det)_{\Phi}(\delta\Phi)
=
(\delta\varphi_{11})\varphi_{22} 
+ \varphi_{11}(\delta \varphi_{22}) 
- (\delta \varphi_{12})\varphi_{21} 
- \varphi_{12}(\delta \varphi_{21}),
$$
where
$$
\delta\Phi 
:= 
\begin{pmatrix}\delta\varphi_{11} & \delta\varphi_{12} \\ 
\delta\varphi_{21} & \delta\varphi_{22} \end{pmatrix}
\in 
C^\8(U,\gl(2,\CC)).
$$
Since $\{\Phi\neq 0\}$ is a dense open subset of $X$,
the union of the complements of each of the zero-sets of the functions
$\varphi_{ij}$ is a dense open subset of $U$.  Now suppose that $\eta \in
C^0(\det E\otimes\det W^+)$ lies in $\Coker(D\det)_\Phi$, so
$((D\det)_{\Phi}(\delta\Phi),\eta)_{L^2(X)} = 0$ for all $\delta\Phi
\in C^0(W^+\otimes E)$. Since $\{\Phi\neq 0\}$ is a dense open subset of
$U$, we have $\eta\equiv 0$ on $U$ and, as $U$ was an arbitrary open subset
of $X$, we have $\eta \equiv 0$ on $X$.
\end{proof}

Aronszajn's theorem then yields the following useful consequence of the
preceding lemma.

\begin{cor}
\label{cor:SmoothProjVariety} 
If $[A,\Phi,f,\vartheta]$ is a point in
$\fM_{W,E}^{*,0}$ such that $\det\Phi = 0$, then $[\Phi]$ is a smooth point
of the zero locus $\fX = \det^{-1}(0) \subset \PP(\bE^+)$, that is,
$$
\pi(\tfM_{W,E}^{*,0})\subset\fX_{\sm},
$$
where $\pi:\tfM_{W,E}^{*,0}\to\PP(\bE^+)$ is the projection.
\end{cor}

\begin{proof}
The hypotheses imply that $\Phi \in \Ker \Dirac_{A,f,\vartheta}$. Thus,
Aronszajn's theorem \cite{Aronszajn}
implies that $\{\Phi\neq 0\}$ is a dense open subset of
$X$, so $[\Phi]$ is a smooth point of $\fX$ by Lemma
\ref{lem:SmoothProjVariety}.
\end{proof}

For each integer $\kappa \geq n$, 
the variety $\fX\subset\PP(\bE^+)$ defines `rank-one loci of
incident planes', 
$$
\tII_\kappa(\fX) := \pi_1^{-1}(\fX) \subset \FF_\kappa(\bE^+)
\quad\text{and}\quad
\II_\kappa(\fX) := \pi_2(\pi_1^{-1}(\fX)) \subset \GG_\kappa(\bE^+).
$$
Corollary \ref{cor:SmoothProjVariety} implies that we shall
only need to consider the incident loci $\tII_\kappa(\fX_{\sm})$ and
$\II_\kappa(\fX_{\sm})$ of the smooth
part of $\fX$. In this case, $\tII_\kappa(\fX_{\sm})$ is a smooth
submanifold of $\FF_\kappa(\bE^+)$ with codimension
\begin{equation}
\codim(\tII_\kappa(\fX_{\sm});\FF_\kappa(\bE^+))
=
\codim(\fX_{\sm};\PP(\bE^+))
=
\8,
\label{eq:CodimIncLocusInFlag}
\end{equation}
although $\II_\kappa(\fX_{\sm})$ is not necessarily a smooth submanifold of 
$\GG_\kappa(\bE^+)$.

As we saw in \S \ref{subsec:Fredholm}, the Grassmann manifolds
$\GG_\kappa(\bE^+)$ or $\GG(\bE^+)$ do not provide suitable target
manifolds for our Dirac-operator period map, so we instead consider the
preimages of the rank-one loci in $\Fred_{\kappa,n}(\bE^+,\bE^-)$ using
the projection \eqref{eq:FredToGrassmann} and then, in turn, in
$\Flag_{\kappa,n}(\bE^+,\bE^-)$. Define
\begin{equation}
\JJ_\kappa(\fX) 
:= 
\pi^{-1}(\II_\kappa(\fX)) 
\subset
\Fred_{\kappa,n}(\bE^+,\bE^-)
\subset
\Fred_n(\bE^+,\bE^-),
\end{equation}
where $\pi:\Fred_{\kappa,n}(\bE^+,\bE^-)\to \GG_\kappa(\bE^+)$ is the canonical
map. The space $\JJ_\kappa(\fX)$ need not be a smooth manifold, so we also
define
\begin{equation}
\tJJ_\kappa(\fX) 
:= 
\pi^{-1}(\tII_\kappa(\fX)) 
\subset
\Flag_{\kappa,n}(\bE^+,\bE^-),
\end{equation}
which is a smooth submanifold since the canonical map
$\pi:\Flag_{\kappa,n}(\bE^+,\bE^-)\to \FF_\kappa(\bE^+)$ is a submersion
according to Lemma \ref{lem:FredFlagToGrassmannFlagSubmersion}. This last
observation also implies
\begin{equation}
\codim(\tJJ_\kappa(\fX_{\sm});\Flag_{\kappa,n}(\bE^+,\bE^-))
=
\codim(\tII_\kappa(\fX_{\sm});\FF_\kappa(\bE^+))
=
\8,
\label{eq:CodimIncLocusInFlagOp}
\end{equation}
courtesy of \eqref{eq:CodimIncLocusInFlag}.

\subsection{The Dirac-operator period map and its differential}
\label{subsec:DiffPeriodMap}
Recall that the principal goal of \S \ref{sec:Period} is to show that for
generic parameters $(f,\vartheta)$ and any point $[A,\Phi]$ in
$M_{W,E}^{*,0}(f,\tau,\vartheta)$, the section $\Phi$ in $\Ker
D_{A,f,\vartheta}$ is not rank one. Hence, it is natural to define a
`period map', using the Dirac operator, to detect the rank-one loci that we
wish to avoid, namely $\JJ_\kappa(\fX)$ in $\Fred_{\kappa,n}(\bE^+,\bE^-)$
or $\tJJ_\kappa(\fX)$ in $\Flag_{\kappa,n}(\bE^+,\bE^-)$.

We therefore define a smooth, $\ssG_E$-equivariant map
\begin{equation}
\ubarP:\tsC_{W,E}\times\sP 
\to
\Fred_n(\bE^+,\bE^-),
\quad
(A,\Phi,f,\vartheta)
\mapsto
\Dirac_{A,f,\vartheta},
\label{eq:DefnPeriodMap}
\end{equation}
by analogy with \cite[\S VI]{DonConn}, \cite[Eq. (4.3.14)]{DK}.  The
map $\ubarP$ could obviously be defined simply on
$\sA_E\times\sP$, but we shall soon need to consider it as a map on the
universal moduli space $\tfM_{W,E}^{*,0}\subset\tsC_{W,E}^{*,0}\times\sP$.
The map \eqref{eq:DefnPeriodMap} has differential
$$
(D\ubarP)_{A,\Phi,f,\vartheta}: T_{A,\Phi}\tsC_{W,E}\oplus T_{f,\vartheta}\sP 
\to
\Hom_\CC(\bE^+,\bE^-),
$$
given by Lemma \ref{lem:VarDirac},
$$
(a,\phi,\delta f,\delta\vartheta)
\mapsto
(D\Dirac)_{A,f,\vartheta}(a,\delta f,\delta\vartheta),
$$
noting that $T_{\ubarP(A,\Phi,f,\vartheta)}\Fred_n(\bE^+,\bE^-)
= \Hom_\CC(\bE^+,\bE^-)$.
We then have the following analogue of Proposition 4.3.14 in \cite{DK}.

\begin{prop}
\label{prop:SurjectiveDiffPeriod}
Suppose $(A,\Phi,f,\vartheta)$ represents a point in 
the universal moduli space $\fM_{W,E}^{*,0} \subset
\sC_{W,E}^{*,0}\times\sP$. Then the following partial differential is
surjective: 
$$
(D\ubarP)_{A,\Phi,f,\vartheta}(0,\cdot)
: \{0\}\oplus T_{f,\vartheta}\sP 
\to 
T_{\ubarP(A,\Phi,f,\vartheta)}\Fred_n(\bE^+,\bE^-),
$$
where $(D\ubarP(A,\Phi,\cdot))_{f,\vartheta} = 
(D\ubarP)_{A,\Phi,f,\vartheta}(0,\cdot)$.
\end{prop}

Proposition \ref{prop:SurjectiveDiffPeriod} plays a crucial role in the
proof of Theorem \ref{thm:NoRankOneIrreducibles}. As we shall see below, it
follows easily from

\begin{prop}
\label{prop:CokerDifferentialPeriod}
Assume $(A,\Phi)$ is an irreducible, non-zero-section $\PU(2)$ monopole on $X$
for the perturbation parameters $f,\vartheta$ (and some $\tau$). If
$\phi\in\Om^0(W^+\otimes E)$ and $\psi \in \Om^0(W^-\otimes E)$ satisfy
\begin{equation}
\Real\left( (D\Dirac_A)_{f,\vartheta}(\delta f,\delta\vartheta),
\psi\otimes\phi^* \right)_{L^2(X)} =0,
\label{eq:RealDerivDiracOrthogImage}
\end{equation}
for all $(\delta f,\delta\vartheta)$,
then $\psi\otimes\phi^* \equiv 0$ on $X$.
\end{prop}

It is important to remember that 
$$
\Ran(D\ubarP)_{A,\Phi,f,\vartheta}
\subset
\Hom_\CC(\bE^+,\bE^-)
$$
is only a {\em real\/} subspace and so if there are elements of the
latter tangent space which do not lie in $\Ran(D\ubarP)_{A,\Phi,f,\vartheta}$,
then we can only assume there are tangent vectors obeying the real
$L^2$-orthogonality condition of Proposition
\ref{prop:CokerDifferentialPeriod}. Our proof of Proposition
\ref{prop:CokerDifferentialPeriod} ultimately relies on the
following crucial `unique continuation' result for reducible $\PU(2)$
monopoles:

\begin{thm} 
\label{thm:LocalToGlobalReducible}
\cite[Theorem 5.11]{FL1}
Suppose $(A,\Phi)$ is a solution to the perturbed $\PU(2)$ monopole equations
\eqref{eq:PT} over a connected, oriented, smooth
four-manifold $X$ with smooth Riemannian metric $g$ such that 
$(A,\Phi)$ is reducible
on a non-empty open subset $U\subset X$. Then $(A,\Phi)$ is reducible on
$X$ if
\begin{itemize}
\item $\Phi\not\equiv 0$ on $X$, or
\item $\Phi\equiv 0$, and $M_E^{\asd}(g)$ contains no {\em twisted
reducibles\/} or $U$ is {\em suitable\/}.
\end{itemize}  
\end{thm}

See \cite[p. 586]{KMStructure} for the definition of `twisted reducibles' and
see \cite[p. 589]{KMStructure} for the definition of a `suitable' open set.
A detailed proof of Theorem \ref{thm:LocalToGlobalReducible} is given in
\cite{FL1}. Our argument relies on the Agmon-Nirenberg unique continuation
theorem for an ordinary differential equation on a Hilbert space
\cite{Agmon}, \cite{AgmonNirenberg} and it generalizes the method used by
Donaldson and Kronheimer 
to prove the corresponding unique continuation result for reducible
anti-self-dual connections \cite[Lemma 4.3.21]{DK}. The proof of
transversality via holonomy perturbations given in \cite{FL1} relies on a
refined version of this unique continuation property.

\begin{proof}[Proof of Proposition \ref{prop:CokerDifferentialPeriod}]
Suppose $\psi\otimes\phi^* \not\equiv 0$ on $X$.  According to Lemma
\ref{lem:HermitianOrthog}, varying $\delta\vartheta$ alone implies that
both $\phi$ and $\psi$ have pointwise complex-orthogonal images in $E$ and
so they have at most complex rank one at each point of the subset
$U:=\{\phi\neq 0\}\cap\{\psi\neq 0\}\subset X$, which we assume is
non-empty. Over the open subset $U\subset X$, the section $\phi$ defines a
Hermitian line subbundle $L_1 := \CC\cdot\phi|_U \subset E|_U$.  Let $L_2
\subset E|_U$ be the Hermitian line subbundle given by $L_2 := L_1^\perp \simeq
(E|_U)/L_1$, so $L_2 \simeq (\det E|_U)\otimes L_1^*$ and we have a unitary
splitting of Hermitian vector bundles: $E|_U = L_1 \oplus L_2$.
Furthermore, noting that $\phi$, $\psi$ have pointwise complex-orthogonal
images in $E$, we have $\Ran\psi\subset L_2$.

With respect to this splitting, the unitary connection $A|_U$ on $E|_U$ may
be written as
$$
\cov_A 
=
\begin{pmatrix}
\cov_{A_1} & -b^*          \\
b       & \cov_{A_2} 
\end{pmatrix}
:
\Omega^0(U,L_1\oplus L_2) \to \Omega^1(U,L_1\oplus L_2),
$$
where $A_i$ is a unitary connection on $L_i$, $i=1,2$, and
$b\in\Omega^1(U,L_2\otimes L_1^*)$ is the second fundamental form of the
pair $(A_1,L_1)$ with respect to $(A,E)|_U$, and $b^*$ is the conjugate
transpose of $b$ \cite[\S I.6]{Kobayashi}. 
(In terms of holomorphic bundles over a complex surface
$X$, the second fundamental form $b$ may be identified with an element of
$\Ext^1(L_1,L_2)$.)  Clearly, $A|_U$ is a reducible connection with respect
to the splitting $E|_U=L_1\oplus L_2$ if and only if $b \equiv 0$. Since
$(A,\Phi)$ is irreducible and non-zero-section by hypothesis, Theorem
\ref{thm:LocalToGlobalReducible} implies that $(A,\Phi)|_U$ cannot be reducible
and so $b\not\equiv 0$ on the non-empty open subset $U\subset X$.

{}From Lemma \ref{lem:VarDirac} and the splitting $A|_U = A_1\oplus A_2$ on
$E|_U = L_1\oplus L_2$, we have (over $U\subset X$)
$$
\cov_A\phi 
=
\begin{pmatrix}
\cov_{A_1} & -b^*          \\
b       & \cov_{A_2} 
\end{pmatrix}
\begin{pmatrix}
\phi \\
0
\end{pmatrix}
=
\begin{pmatrix}
\cov_{A_1}\phi \\
b\phi
\end{pmatrix}
\in
\Omega^1(U,W^+\otimes L_1)\oplus \Omega^1(U,W^+\otimes L_2).
$$
Therefore,
\begin{align*}
(D\Dirac_A)_{f,\vartheta}(\delta f,\delta\vartheta)\phi
&=
\sum_{i=1}^4\rho_+((\delta f)e^i)\cov_{A,e_i}\phi
+ \rho_+(f(\delta\vartheta))\phi
\\
&=
\sum_{i=1}^4\rho_+((\delta f)e^i)\cov_{A_1,e_i}\phi
+ \sum_{i=1}^4\rho_+((\delta f)e^i)b(e_i)\phi
+ \rho_+(f(\delta\vartheta))\phi,
\end{align*}
where $\{e_i\}$ is a local oriented, orthonormal frame for $TX$, with dual
coframe $\{e^i\}$ for $T^*X$. Thus,
\begin{equation}
(D\Dirac_A)_{f,\vartheta}(\delta f,\delta\vartheta)\phi
=
\sum_{i=1}^4\rho_+((\delta f)e^i)\cov_{A_1,e_i}\phi
+ \rho_+(f(\delta\vartheta))\phi
+ \rho_+((\delta f)b)\phi. 
\label{eq:CovDerivSplit}
\end{equation}
The first two terms on the right-hand side of \eqref{eq:CovDerivSplit}
are in $\Om^0(U,W^-\otimes L_1)$,
while the last term is in $\Om^0(U,W^-\otimes L_2)$. {}From
\eqref{eq:CovDerivSplit}, the fact that the 
sections $\phi$, $\psi$ have complex-orthogonal images in $E$ at each point
of $X$, and the real $L^2$-orthogonality condition
\eqref{eq:RealDerivDiracOrthogImage}, we obtain
$$
\Real((D\Dirac_A)_{f,\vartheta}(\delta f,\delta\vartheta)\phi,\psi)_{L^2(X)} 
=
\Real(\rho_+((\delta f)b)\phi,\psi)_{L^2(X)} = 0, 
$$
for all $\delta f\in C^\8(\gl(T^*X))$. Since $\phi$,
$i\psi$ also have pointwise complex-orthogonal images in $E$, we may
replace $\psi$ by $i\psi$ above; combining the resulting two real
$L^2$-orthogonality equations yields the complex $L^2$-orthogonality
identity,
$$
(\rho_+((\delta f)b)\phi,\psi)_{L^2(X)} = 0, 
$$
for all $\delta f\in C^\8(\gl(T^*X))$,
which in turn, therefore, yields the pointwise identity:
\begin{equation}
\langle\rho_+((\delta f)_xb_x)\phi_x,\psi_x\rangle_x = 0, 
\qquad x\in U,
\label{eq:ComplexOrthogPointwisePhiPsi}
\end{equation}
for all $(\delta f)_x\in \gl(T^*X)_x$. Note that 
\begin{align*}
\rho_+((\delta f)b)
&\in \Hom_\CC(W^+\otimes_\CC L_1,W^-\otimes_\CC L_2) \\
&\simeq
\Hom_\CC(W^+,W^-)\otimes_\CC\Hom_\CC(L_1,L_2)
\end{align*}
and so
$$
\rho_+((\delta f)_xb_x)
\in 
\Hom_\CC(W^+,W^-)_x\otimes_\CC\Hom_\CC(L_1,L_2)_x
\simeq
\gl(2,\CC).
$$
Since $\rho_+((\delta f)_xb_x)\phi_x$ and $\psi_x$ are orthogonal with
respect to the full Hermitian inner product by
\eqref{eq:ComplexOrthogPointwisePhiPsi}, we may use the fact that
the complexification $\rho_+(T^*X)_x\otimes_\RR\CC $ is equal to
$\Hom_\CC(W^+,W^-)_x \simeq \gl(2,\CC)$ \cite[p. 89]{MorganSWNotes}. 
(See the proof in Lemma \ref{lem:MorganLemma}.)
Therefore, the identity \eqref{eq:ComplexOrthogPointwisePhiPsi} would yield
$$
\psi_x\otimes\phi_x^* = 0 \in \Hom_\CC(W^+,W^-)_x,
$$
for all points $x\in U$ at which $b_x\neq 0$.
Since $\psi_x\otimes\phi_x^* \ne 0$ for any point $x\in U$ by assumption,
this would imply that $b \equiv 0$ on $U$ and so $A|_U$ would be reducible,
a contradiction. Hence, the open set $U\subset X$ must be empty and thus
$\psi\otimes\phi^* \equiv 0$ on $X$, as desired.
\end{proof}

We can now conclude the proof of Proposition
\ref{prop:SurjectiveDiffPeriod}:

\begin{proof}[Proof of Proposition \ref{prop:SurjectiveDiffPeriod}]
Suppose $(A,\Phi,f,\vartheta)$ represents a point in the universal moduli
space $\fM_{W,E}^{*,0}$. For convenience, we denote $p := (f,\vartheta)$
and $\delta p := (\delta f,\delta\vartheta)$. If $\{\phi_a\}_{a\in\NN}$ is
an orthonormal basis for the Hilbert space $\bE^+$ and
$\{\varphi_b\}_{b\in\NN}$ is an orthonormal basis for the Hilbert space
$\bE^-$, then $\{\varphi_b\otimes\phi_a^*\}_{(a,b)\in\NN\times\NN}$ is an
orthonormal basis for the Hilbert space $\Hom_\CC(\bE^+,\bE^-)$, where
$\phi_a^* = \langle\cdot,\phi_a\rangle$.
If $\Ran (D\ubarP)_{A,\Phi,p}(0,\cdot) = \Ran (D\ubarP(A,\Phi,\cdot))_p
\subsetneqq T_{\ubarP(A,\Phi,p)}\Fred_n(\bE^+,\bE^-)$, then there are
sections $\phi \in \bE^+$ and $\varphi \in \bE^-$, so
$$
\varphi\otimes\phi^*
\in
\Hom_\CC(\bE^+,\bE^-)
= 
T_{\Dirac_{A,p}}\Fred_n(\bE^+,\bE^-),
$$
with $\varphi\otimes\phi^* \not\equiv 0$ on $X$ such that
$$
\Real((D\ubarP(A,\Phi,\cdot))_p(\delta p),
\varphi\otimes\phi^*)_{L^2(X)} = 0,
$$
for all $\delta p$. Hence, for all $\delta p$ we have
$$
\Real((D\Dirac_A)_p(\delta p)\phi,\varphi)_{L^2(X)} = 0
$$
and then Proposition \ref{prop:CokerDifferentialPeriod} implies that
$\varphi\otimes\phi^* = 0$, a contradiction. Hence, the differential
$(D\ubarP(A,\Phi,\cdot))_p$ is surjective, as claimed.
\end{proof}

\subsection{The Sard-Smale theorem and transversality}
\label{subsec:SardSmale}
We shall need a special form of the Sard-Smale theorem \cite{SmaleSard} for our
proof of Theorem \ref{thm:NoRankOneIrreducibles}. The result is well-known
and is essentially Proposition 4.3.10 in \cite{DK}, but we shall require a
more detailed statement than that given there, so we summarize the relevant
discussion from
\cite[\S 4.3.1 \& \S 4.3.2]{DK} and prove the precise version we need
here. 

Let $\sC$, $\sP$, and $\sF$ be $C^\8$ Banach manifolds.  Suppose $\sM
\subset\sC\times\sP$ is a $C^\8$ Banach submanifold and that the
restriction $\pi_{\sM;\sP}:\sM\to\sP$ of the projection
$\pi_\sP:\sC\times\sP\to\sP$ onto the second factor is Fredholm. Let
$$
\ubarP:\sM \subset\sC\times\sP \to \sF
$$ 
be a $C^\8$ map which is transverse to a $C^\8$ Banach submanifold
$\sJ\subset \sF$.

\begin{prop}
\label{prop:SmalePerturbMap}
Continue the notation of the preceding paragraph.  Then there is a
first-category subset $\sP_{\fc}\subset \sP$ such that the following holds. 
For all $p$ in $\sP - \sP_{\fc}$, we have
\begin{itemize}
\item $M := \pi_{\sM;\sP}^{-1}(p)$ is a $C^\8$ manifold of dimension
$\Ind(\pi_{\sM;\sP})_p < \8$, 
\item $P := \ubarP(\cdot,p):M\to\sF$ is transverse
to the submanifold $\sJ\subset\sF$, 
\item $Z := P^{-1}(\sJ) \subset M$ is a $C^\8$
submanifold of codimension \hfill\break$\codim(Z;M) = \codim(\sJ;\sF)$.
\end{itemize}
\end{prop}

\begin{rmk}
When $\sC$ is finite dimensional and $\sJ$ has finite codimension in $\sF$,
the preceding result is proved in \cite[\S 3.6B]{AMR}.
\end{rmk}

\begin{proof}
The proof is similar to those of Propositions 4.3.10 and 4.3.11 in
\cite{DK} (see \cite[p. 143]{DK}). The first item follows immediately from
the Sard-Smale theorem \cite{SmaleSard}, for some first-category subset
$\sP_{\fc} \subset \sP$. The hypotheses imply that the preimage
$$
\sZ: = \ubarP^{-1}(\sJ) \subset \sM \subset \sC\times\sP
$$
is a $C^\8$ Banach submanifold of $\sM$. Let $p$ be a regular value of
the projections
\begin{align*}
&\pi_{\sZ;\sP}:\sZ \subset \sC\times\sP \to \sP,
\\
&\pi_{\sM;\sP}:\sM \subset \sC\times\sP \to \sP,
\end{align*}
so $p \in \sP-\sP_{\fc}$, for some possibly larger first-category subset
$\sP_{\fc} \subset \sP$, and let $(c,p)$ be any point in the fiber 
$$
Z := \pi_{\sZ;\sP}^{-1}(p) = \ubarP(\cdot,p)^{-1}(\sJ) = P^{-1}(\sJ),
$$
and let $h = \ubarP(c,p) \in \sF$.
Similarly, let $M := \pi_{\sM;\sP}^{-1}(p)$ and note that $Z$ and $M$ are
$C^\8$ manifolds. 

Note that $T_{(c,p)}\sZ$ and $T_{(c,p)}\sM$ are subspaces of $T_c\sC\oplus
T_p\sP$. By choice of $p\in \sP-\sP_{\fc}$ and by hypothesis, respectively,
the maps  
\begin{align}
(D\pi_{\sZ;\sP})_{(c,p)} &: T_{(c,p)}\sZ \to T_p\sP,
\label{eq:ProjSurj}
\\
(D\ubarP)_{(c,p)} &: T_{(c,p)}\sM \to T_h\sF \to T_h\sF/T_h\sJ,
\label{eq:PeriodSurj}
\end{align}
are surjective and
\begin{align*}
T_c M\oplus\{0\} 
&= (D\pi_{\sM;\sP})_{(c,p)}^{-1}(0) \subset T_{(c,p)}\sM, 
\\
T_{(c,p)}\sZ &= (D\ubarP)_{(c,p)}^{-1}(T_h\sJ) \subset T_{(c,p)}\sM,
\end{align*}
where $T_c M \simeq T_c M \oplus\{0\} \subset T_c\sC\oplus\{0\}$. 
We claim that the map $P:M\to\sF$ is transverse to $\sJ\subset\sF$, so the
preimage $Z = P^{-1}(\sJ)$ is a $C^\8$ manifold of codimension
$\codim(Z;M) = \codim(\sJ;\sF)$.

Suppose $\delta h \in T_h\sF$.  By \eqref{eq:PeriodSurj} there is a tangent
vector $(\delta c,\delta p) \in T_{(c,p)}\sM$ such that
$(D\ubarP)_{(c,p)}(\delta c,\delta p) = \delta h
\pmod{T_h\sJ}$. According to \eqref{eq:ProjSurj}, given $\delta p$,  there
is a tangent vector of the form $(\delta c',\delta p) \in T_{(c,p)}\sZ$ and
so 
$$
(\delta c - \delta c', 0)
=
(\delta c,\delta p) - (\delta c',\delta p)
\in
T_c M\oplus\{0\}
\subset
T_{(c,p)}\sM.
$$
Since $(D\ubarP)_{(c,p)}(\delta c',\delta p) \in T_h\sJ$ (because
$\ubarP^{-1}(\sJ) = \sZ$) and $P = \ubarP(\cdot,p)$, this gives
\begin{align*}
(DP)_c(\delta c - \delta c')
&=
(D\ubarP)_{(c,p)}(\delta c - \delta c', 0) \\
&=
(D\ubarP)_{(c,p)}(\delta c,\delta p) - (D\ubarP)_{(c,p)}(\delta c',\delta p) \\
&=
\delta h \pmod{T_h\sJ}.
\end{align*}
Hence, the composition of the differential and quotient map,
$$
(DP)_p: T_c M \to T_h\sF \to T_h\sF/T_h\sJ,
$$
is surjective, as desired.
\end{proof}

\subsection{The Sard-Smale theorem and semi-Kuranishi models}
\label{subsec:Kuranishi}
We shall use Proposition \ref{prop:SmalePerturbMap} to prove Theorem
\ref{thm:NoRankOneIrreducibles}. As the detailed argument is 
technical and lengthy, we shall first outline the basic idea 
(under some simplifying assumptions) and establish a few
notational conventions that will be useful in this section. 

The space $\tsC_{W,E}^{*,0}\times\sP$ is a principal
$\ssG_E$-bundle over the base manifold $\sC_{W,E}^{*,0}\times\sP$. The
parametrized period map $\ubarP$ may then be viewed as a section of the
associated Banach vector bundle
$$
\begin{CD}
(\tsC_{W,E}^{*,0}\times\sP\times\Fred_n(\bE^+,\bE^-))/\ssG_E \\
@VVV \\
\sC_{W,E}^{*,0}\times\sP
\end{CD}
$$
In order to obtain a Fredholm projection map onto the parameter space
$\sP$, we need to restrict the base of the  
preceding bundle to the universal moduli
space $\fM_{W,E}^{*,0}\subset \sC_{W,E}^{*,0}\times\sP$. Of course, the
space $\fM_{W,E}^{*,0}$ is not necessarily a smooth manifold, so in
the argument below we replace small open neighborhoods in $\fM_{W,E}^{*,0}$
by `thickened universal moduli spaces' containing the smooth locus of
$\fM_{W,E}^{*,0}$ as finite-dimensional submanifolds. Temporarily assuming
$\fM_{W,E}^{*,0}$ is smooth, for the purposes of this outline,
we thus consider the Banach vector bundle
together with a Fredholm projection map:
$$
\begin{CD}
(\tfM_{W,E}^{*,0}\times\Fred_n(\bE^+,\bE^-))/\ssG_E \\
@VVV \\
\fM_{W,E}^{*,0}
\end{CD}
\qquad\text{and}\qquad
\begin{CD}
\fM_{W,E}^{*,0} \\
@VV{\pi_{\fM;\sP}}V \\
\sP
\end{CD}
$$
It is rather cumbersome to phrase the transversality argument in terms of
sections of vector bundles rather than maps of Banach manifolds, although
this can be done. So, before proceeding further, we restrict to slice
neighborhoods in $\sC_{W,E}^{*,0}$ and $\fM_{W,E}^{*,0}$, which for
simplicity we continue to write as $\sC_{W,E}^{*,0}$ and $\fM_{W,E}^{*,0}$
here, abusing notation slightly, rather than introduce new symbols. With
respect to such a trivialization, we then have a parametrized period
map and a Fredholm projection:
\begin{align*}
\ubarP: &\fM_{W,E}^{*,0}\subset \sC_{W,E}^{*,0}\times\sP 
\to \Fred_n(\bE^+,\bE^-),
\\
\pi_{\fM;\sP}: &\fM_{W,E}^{*,0}\subset \sC_{W,E}^{*,0}\times\sP \to \sP.
\end{align*}
Fix an integer $\kappa \ge n:= \Ind\sD_{A,f,\vartheta}$ and temporarily
assume, again for the purposes of this outline, that the rank-one locus
$\JJ_\kappa(\fX_{\sm})\subset \Fred_n(\bE^+,\bE^-)$ is a smooth submanifold. 
If the parametrized period map $\ubarP:\fM_{W,E}^{*,0}\to
\Fred_n(\bE^+,\bE^-)$ is transverse to the locus
$\JJ_\kappa(\fX_{\sm})\subset 
\Fred_n(\bE^+,\bE^-)$, then Proposition \ref{prop:SmalePerturbMap} implies
that the 
period map $P := \ubarP(\cdot, p)$ from the fiber $\pi_{\fM;\sP}^{-1}(p) =
M_{W,E}^{*,0}(p)$,
$$
P:M_{W,E}^{*,0}(p)\to\Fred_n(\bE^+,\bE^-),
$$
will be transverse to $\JJ_\kappa(\fX_{\sm})\subset
\Fred_n(\bE^+,\bE^-)$ for all $p\in\sP-\sP_{\fc}$, for some first-category
subset $\sP_{\fc}$, and so
$$
\codim(P^{-1}(\JJ_\kappa(\fX_{\sm});M_{W,E}^{*,0}(p)) 
= 
\codim(\JJ_\kappa(\fX_{\sm});\Fred_n(\bE^+,\bE^-)) = \8.
$$
Since $M_{W,E}^{*,0}(p)$ is a finite-dimensional smooth manifold and
$P^{-1}(\JJ_\kappa(\fX_{\sm}))$ is a smooth submanifold with infinite
codimension, the subset $P^{-1}(\JJ_\kappa(\fX_{\sm}))$ is necessarily empty.
Hence, after repeating this argument for each integer $\kappa\ge n$
and observing that the subset
$$
\bigcup_{\kappa\ge n}P^{-1}(\JJ_\kappa(\fX_{\sm}))
=
\{[A,\Phi] \in M_{W,E}^{*,0}(p):\text{$\Ker D_{A,p}$ contains some 
rank-one $\Psi$}\}
$$
is empty for generic $p\in\sP$,
we will have shown that for generic $p \in \sP$, the moduli space
$M_{W,E}^{*,0}(p)$ contains no rank-one pairs.

Given this outline, we now proceed to the detailed proof of Theorem
\ref{thm:NoRankOneIrreducibles}. 

\begin{proof}[Proof of Theorem \ref{thm:NoRankOneIrreducibles}] 
Recall that in order to obtain a Fredholm
projection map --- thus permitting an application of the Sard-Smale theorem
\cite[\S 4.3.1 \& \S 4.3.2]{DK}, \cite{SmaleSard} 
in the form of Proposition \ref{prop:SmalePerturbMap} --- we need to
restrict our attention to the universal moduli space of irreducible,
non-zero-section $\PU(2)$ monopoles:
$$
\fM_{W,E}^{*,0}
:= 
\{(A,\Phi,p) \in \tsC_{W,E}^{*,0}\times\sP: \fS(A,\Phi,p) = 0\}/\ssG_E.
$$
Let $(A_0,\Phi_0,p_0)$ be any point in $\fM_{W,E}^{*,0}$.
The space $\fM_{W,E}^{*,0}$ is not necessarily smooth so we first construct a
{\em semi-Kuranishi\/} model for an open neighborhood
$(A_0,\Phi_0,p_0)$ in $\fM_{W,E}^{*,0}$. Let
$$
\bV := L^2_{k-1}(\Lambda^+\otimes\su(E))\oplus L^2_{k-1}(W^+\otimes E)
$$
and observe that the map
$$
\ubarfS:\sC_{W,E}^{*,0}\times\sP \to \bV
$$
is {\em right semi-Fredholm\/}, in the sense that its differentials are
right semi-Fredholm operators, having closed range and finite-dimensional
cokernel (though not finite-dimensional kernels) \cite[Definition XI.2.3 \&
Proposition XI.2.10]{Conway}: 

\begin{claim}
The map $\ubarfS:\sC_{W,E}^{*,0}\times\sP \to \bV$ is right semi-Fredholm.
\end{claim}

\begin{proof}
One can see this by observing that the differentials
$$
(D\ubarfS)_{(A,\Phi,p)}:
T_{(A,\Phi)}\sC_{W,E}^{*,0} \oplus T_p\sP \to \bV
$$
have right Green's operators defined by
$$
G_{(A,\Phi,p)}
=
(D\ubarfS)_{(A,\Phi,p)}^\dagger
\left((D\ubarfS)_{(A,\Phi,p)}(D\ubarfS)_{(A,\Phi,p)}^\dagger\right)^{-1}:
\bV\to T_{(A,\Phi)}\sC_{W,E}^{*,0} \oplus T_p\sP.
$$
Let $\Pi_{(A,\Phi,p)}$ denote the $L^2$-orthogonal projection from
$\bV$ onto the finite-dimensional subspace
$$
\bH_{(A,\Phi,p)}^2 
:= 
\left(\Ran (D\ubarfS(\cdot,p))_{(A,\Phi)}\right)^\perp \subset \bV,
$$
and let $\Pi_{(A,\Phi,p)}^\perp := \id - \Pi_{(A,\Phi,p)}$ be the
$L^2$-orthogonal projection from $\bV$ onto the subspace $\Ran
(D\ubarfS(\cdot,p))_{(A,\Phi)}$.  Then,
$$
(D\ubarfS)_{(A,\Phi,p)}\circ G_{(A,\Phi,p)}
= \id - \Pi_{(A,\Phi,p)}:\bV\to\bV
$$
and so is right semi-Fredholm by \cite[Definition XI.2.3]{Conway}, as
$\Pi_{(A,\Phi,p)}$ is a finite-rank operator and thus compact.
\end{proof}

We now consider the pair of maps
\begin{align}
\Pi_{(A_0,\Phi_0,p_0)}^\perp\circ\ubarfS &:\sC_{W,E}^{*,0}\times\sP 
\to (\bH_{(A_0,\Phi_0,p_0)}^2)^\perp\cap\bV,
\label{eq:StabilizedCanonicalMap} 
\\
\Pi_{(A,_0\Phi_0,p_0)}\circ\ubarfS       &:\sC_{W,E}^{*,0}\times\sP 
\to \bH_{(A_0,\Phi_0,p_0)}^2.
\label{eq:ObstructMap}
\end{align}
By construction, the differential of the stabilized canonical map
\eqref{eq:StabilizedCanonicalMap} is surjective at the point
$(A_0,\Phi_0,p_0)$ and thus also on some open neighborhood
$\sU_{(A_0,\Phi_0,p_0)}$ of the point $(A_0,\Phi_0,p_0)$ in
$\sC_{W,E}^{*,0}\times\sP$. An open neighborhood of $(A_0,\Phi_0,p_0)$ in
$\fM_{W,E}^{*,0}$ is then given by the zero locus of the finite-rank
obstruction map \eqref{eq:ObstructMap} in the $C^\8$ Banach manifold
\begin{equation}
\sM_{(A_0,\Phi_0,p_0)}
:=
\sU_{(A_0,\Phi_0,p_0)}
\cap
(\Pi_{(A_0,\Phi_0,p_0)}^\perp\circ\ubarfS)^{-1}(0)
\subset
\sC_{W,E}^{*,0}\times\sP,
\label{eq:ThickenendParamMod}
\end{equation}
comprising a {\em thickened, parametrized moduli space\/}.  Let
$$
\pi_{\sM;\sP}:\sM_{(A_0,\Phi_0,p_0)} \to \sP
$$ 
denote the restriction to
$\sM_{(A_0,\Phi_0,p_0)}$ of the projection from $\sC_{W,E}^{*,0}\times\sP$
onto the second factor, and let
\begin{equation}
\label{eq:DefnH1}
\begin{aligned}
\bH_{(A_0,\Phi_0,p_0)}^1
:&=
\Ker(D\ubarfS(\cdot,p_0))_{(A_0,\Phi_0)}
\\
&= 
\Ker(\Pi_{(A_0,\Phi_0,p_0)}^\perp\circ D\ubarfS(\cdot,p_0))_{(A_0,\Phi_0)}
\\
&=
\Ker(D\pi_{\sM;\sP})_{(A_0,\Phi_0,p_0)}
\subset T_{(A_0,\Phi_0)}\sC_{W,E}^{*,0}
\end{aligned}
\end{equation}
be the Zariski tangent space to the fiber
$M_{W,E}^{*,0}(p_0)$ at the point $(A_0,\Phi_0)$. Note that the differentials
\begin{equation}
(D\pi_{\sM;\sP})_{(A,\Phi,p)}:
T_{(A,\Phi,p)}\sM_{(A_0,\Phi_0,p_0)} \to T_p\sP
\label{eq:DiffProjThickenedModtoParam}
\end{equation}
are Fredholm and that $\bH_{(A_0,\Phi_0,p_0)}^1$ is finite-dimensional.

We now consider the period map:
\begin{equation}
\ubarP:\sM_{(A_0,\Phi_0,p_0)} \to \Fred_n(\bE^+,\bE^-).
\label{eq:ThickenedPeriodMap}
\end{equation}
The map $\ubarP(A_0,\Phi_0,\cdot):\sP\to \Fred_n(\bE^+,\bE^-)$ is a submersion 
by Proposition \ref{prop:SurjectiveDiffPeriod} and in
particular a submersion at the point $p_0$, but it does not follow that
the same is true for the map \eqref{eq:ThickenedPeriodMap} at the point
$(A_0,\Phi_0,p_0)$ since the tangent space
$T_{(A_0,\Phi_0,p_0)}\sM_{(A_0,\Phi_0,p_0)}$ does not necessarily contain
the subspace 
$$
\{0\}\oplus T_{p_0}\sP 
\subset 
T_{(A_0,\Phi_0)}\sC_{W,E}^{*,0} \oplus T_{p_0}\sP,
$$
as required by the hypotheses of Proposition
\ref{prop:SurjectiveDiffPeriod}.  To circumvent this problem and permit an
analysis of the preimage of $\JJ_\kappa(\fX_{\sm})$ in
$\sM_{(A_0,\Phi_0,p_0)}$, we use a second stabilization technique
similar to that employed by Donaldson for the anti-self-dual equation
\cite{DonConn}, \cite[\S 7.2.2]{DK}. 

\begin{claim}
\label{claim:EqualVectorSpaceAndFiniteCodim}
There is an equality of vector spaces,
\begin{equation}
\label{eq:EqualVectorSpace}
T_{(A_0,\Phi_0,p_0)}\sM_{(A_0,\Phi_0,p_0)} + (\{0\}\oplus T_{p_0}\sP)
=
\bH_{(A_0,\Phi_0,p_0)}^1\oplus T_{p_0}\sP,
\end{equation}
and an isomorphism of vector spaces,
\begin{equation}
\label{eq:VectorSpaceIsomorphism}
\left(\bH_{(A_0,\Phi_0,p_0)}^1\oplus T_{p_0}\sP\right)
/T_{(A_0,\Phi_0,p_0)}\sM_{(A_0,\Phi_0,p_0)}
\simeq
\Coker (D\pi_{\sM;\sP})_{(A_0,\Phi_0,p_0)},
\end{equation}
and so $T_{(A_0,\Phi_0,p_0)}\sM_{(A_0,\Phi_0,p_0)}$ has finite
codimension in $\bH_{(A_0,\Phi_0,p_0)}^1\oplus T_{p_0}\sP$.
\end{claim}

\begin{proof}
For every tangent vector 
$$
\delta p \in \Ran (D\pi_{\sM;\sP})_{(A_0,\Phi_0,p_0)} \subset T_{p_0}\sP,
$$
there is a corresponding tangent vector $(a,\phi,\delta p)$ in 
$T_{(A_0,\Phi_0,p_0)}\sM_{(A_0,\Phi_0,p_0)}$, so the equality
\eqref{eq:EqualVectorSpace} follows. 
Since the differential \eqref{eq:DiffProjThickenedModtoParam} is Fredholm,
the cokernel
$$
\left(\Ran (D\pi_{\sM;\sP})_{(A_0,\Phi_0,p_0)}\right)^\perp
\simeq \Coker (D\pi_{\sM;\sP})_{(A_0,\Phi_0,p_0)}
$$
is finite dimensional. Indeed, using the isomorphisms
\begin{align*}
T_{p_0}\sP
&\simeq
\Ran (D\pi_{\sM;\sP})_{(A_0,\Phi_0,p_0)} 
\oplus \Coker (D\pi_{\sM;\sP})_{(A_0,\Phi_0,p_0)},
\\
T_{(A_0,\Phi_0,p_0)}\sM_{(A_0,\Phi_0,p_0)}
&\simeq
(\Ker (D\pi_{\sM;\sP})_{(A_0,\Phi_0,p_0)})^\perp 
\oplus \bH_{(A_0,\Phi_0,p_0)}^1,
\\
&\simeq
\Ran (D\pi_{\sM;\sP})_{(A_0,\Phi_0,p_0)}\oplus \bH_{(A_0,\Phi_0,p_0)}^1,
\end{align*}
(the second isomorphism above being due to
\eqref{eq:DiffProjThickenedModtoParam}), we obtain the isomorphism
\eqref{eq:VectorSpaceIsomorphism}. The assertion on codimension follows
from the fact that $\dim\Coker (D\pi_{\sM;\sP})_{(A_0,\Phi_0,p_0)}<\8$,
since the projection $\pi_{\sM;\sP}$ is Fredholm.
\end{proof}

Proposition \ref{prop:SurjectiveDiffPeriod} asserts that the differential
$$
(D\ubarP)_{(A_0,\Phi_0,p_0)}:\{0\}\oplus T_{p_0}\sP 
\to
T_{\ubarP(A_0,\Phi_0,p_0)}\Fred_n(\bE^+,\bE^-)
$$
is surjective and so the same is true for
$$
(D\ubarP)_{(A_0,\Phi_0,p_0)}:\bH_{(A_0,\Phi_0,p_0)}^1\oplus T_{p_0}\sP
\to
T_{\ubarP(A_0,\Phi_0,p_0)}\Fred_n(\bE^+,\bE^-).
$$
Claim \ref{claim:EqualVectorSpaceAndFiniteCodim} and its proof imply that 
\begin{align*}
&\bH_{(A_0,\Phi_0,p_0)}^1\oplus T_{p_0}\sP
\\
&\quad\simeq
\bH_{(A_0,\Phi_0,p_0)}^1
\oplus (\Ker (D\pi_{\sM;\sP})_{(A_0,\Phi_0,p_0)})^\perp 
\oplus \Coker (D\pi_{\sM;\sP})_{(A_0,\Phi_0,p_0)}
\\
&\quad\simeq
T_{(A_0,\Phi_0,p_0)}\sM_{(A_0,\Phi_0,p_0)}
\oplus \Coker (D\pi_{\sM;\sP})_{(A_0,\Phi_0,p_0)}.
\end{align*}
Define a finite-dimensional subspace
\begin{align*}
F_{(A_0,\Phi_0,p_0)}
:&=
(D\ubarP)_{(A_0,\Phi_0,p_0)}
\left(\Coker (D\pi_{\sM;\sP})_{(A_0,\Phi_0,p_0)}\right)
\\
&\subset
T_{\ubarP(A_0,\Phi_0,p_0)}\Fred_n(\bE^+,\bE^-)
\end{align*}
and let
$$
\iota:F_{(A_0,\Phi_0,p_0)}\to 
T_{P(A_0,\Phi_0,p_0)}\Fred_n(\bE^+,\bE^-)
= \Hom_\CC(\bE^+,\bE^-)
$$
denote the inclusion. Although the differential
$$
(D\ubarP)_{(A_0,\Phi_0,p_0)}:
T_{(A_0,\Phi_0,p_0)}\sM_{(A_0,\Phi_0,p_0)} 
\to 
T_{\ubarP(A_0,\Phi_0,p_0)}\Fred_n(\bE^+,\bE^-)
$$
need not necessarily be surjective, the following linear map is surjective
by construction:
\begin{multline}
\label{eq:StabDifferentialThickenedPeriodMap}
\iota+(D\ubarP)_{(A_0,\Phi_0,p_0)}:
F_{(A_0,\Phi_0,p_0)}\oplus T_{(A_0,\Phi_0,p_0)}\sM_{(A_0,\Phi_0,p_0)} 
\\
\to 
T_{\ubarP(A_0,\Phi_0,p_0)}\Fred_n(\bE^+,\bE^-),
\end{multline}
where, for all $(\ff,\delta\fm) \in 
F_{(A_0,\Phi_0,p_0)}\oplus T_{(A_0,\Phi_0,p_0)}\sM_{(A_0,\Phi_0,p_0)}$,
we set
$$
\left(\iota+(D\ubarP)_{(A_0,\Phi_0,p_0)}\right)(\ff,\delta\fm)
:=
\iota(\ff)+(D\ubarP)_{(A_0,\Phi_0,p_0)}(\delta\fm).
$$
Define a stabilized period map 
\begin{equation}
\iota+\ubarP:F_{(A_0,\Phi_0,p_0)} \times \sM_{(A_0,\Phi_0,p_0)}
\to
\Fred_n(\bE^+,\bE^-)
\label{eq:StabThickenedPeriodMap}
\end{equation}
by setting $(\iota+\ubarP)(\ff,\fm) := \iota(\ff)+\ubarP(\fm)$, 
for all $(\ff,\fm)$ in 
$F_{(A_0,\Phi_0,p_0)} \times \sM_{(A_0,\Phi_0,p_0)}$. (We should write
$\Hom_\CC(\bE^+,\bE^-)$ for the image space in
\eqref{eq:StabThickenedPeriodMap}, but when $\ff$ is close enough to
the origin in $F_{(A_0,\Phi_0,p_0)}$, the point $\iota(\ff)+\ubarP(\fm)$
lies in the open subset $\Fred_n(\bE^+,\bE^-)$ since $\ubarP(\fm)$ lies in
$\Fred_n(\bE^+,\bE^-)$.) By construction, the differential 
\eqref{eq:StabDifferentialThickenedPeriodMap} of the
stabilized period map \eqref{eq:StabThickenedPeriodMap} is surjective at
the point $(0,A_0,\Phi_0,p_0)$ and so there are an open neighborhood of the
origin,
$$
O_{(A_0,\Phi_0,p_0)} \subset F_{(A_0,\Phi_0,p_0)},
$$
and, in the definition \eqref{eq:ThickenendParamMod} of
$\sM_{(A_0,\Phi_0,p_0)}$, a possibly smaller open neighborhood
$\sU_{(A_0,\Phi_0,p_0)}$ of the point $(A_0,\Phi_0,p_0)$ in
$\sC_{W,E}^{*,0}\times\sP$, such that the restriction
\begin{equation}
\iota+\ubarP:O_{(A_0,\Phi_0,p_0)} \times \sM_{(A_0,\Phi_0,p_0)}
\to
\Fred_n(\bE^+,\bE^-)
\label{eq:StabThickenedPeriodSubmersion}
\end{equation}
of the map \eqref{eq:StabThickenedPeriodMap} is a submersion.

Suppose $\kappa\ge n$ is an integer. Because the smooth map 
\eqref{eq:StabThickenedPeriodSubmersion} is a submersion, the preimage
\begin{equation}
\label{eq:DefnStratumOMSpace}
\sN_{(A_0,\Phi_0,p_0)}^\kappa
:=
(O_{(A_0,\Phi_0,p_0)} \times \sM_{(A_0,\Phi_0,p_0)})
\cap
(\iota+\ubarP)^{-1}(\Fred_{\kappa,n}(\bE^+,\bE^-))
\end{equation}
is a smooth submanifold of $O_{(A_0,\Phi_0,p_0)} \times
\sM_{(A_0,\Phi_0,p_0)}$, with (finite) real codimension
\begin{equation}
\label{eq:StratumofOMSpaceCodim}
\codim_\RR\left(\sN_{(A_0,\Phi_0,p_0)}^\kappa;
O_{(A_0,\Phi_0,p_0)} \times \sM_{(A_0,\Phi_0,p_0)}\right)
=
2\kappa(\kappa-n),
\end{equation}
and so
$$
O_{(A_0,\Phi_0,p_0)} \times \sM_{(A_0,\Phi_0,p_0)}
=
\bigcup_{\kappa\ge n} \sN_{(A_0,\Phi_0,p_0)}^\kappa
$$
is a countable, disjoint union of smooth submanifolds.

\begin{claim}
\label{claim:StabThickenedPeriodSubmersionStratum}
The following smooth map is a submersion:
\begin{equation}
\iota+\ubarP:\sN_{(A_0,\Phi_0,p_0)}^\kappa
\to
\Fred_{\kappa,n}(\bE^+,\bE^-).
\label{eq:StabThickenedPeriodSubmersionStratum}
\end{equation}
\end{claim}

\begin{proof}
Let $(A,\Phi,p)$ be a point in 
$\sN_{(A_0,\Phi_0,p_0)}^\kappa$. {}From the definition
\eqref{eq:DefnStratumOMSpace} of the submanifold
$\sN_{(A_0,\Phi_0,p_0)}^\kappa$ and the fact that the map
\eqref{eq:StabThickenedPeriodSubmersion} is a submersion, we see that
$$
T_{(A,\Phi,p)}\sN_{(A_0,\Phi_0,p_0)}^\kappa
=
(D(\iota+\ubarP))_{(A,\Phi,p)}^{-1}
\left(T_{(\iota+\ubarP)(A,\Phi,p)}\Fred_{\kappa,n}(\bE^+,\bE^-)\right),
$$
as the tangent space to the preimage is the preimage of the tangent space.
Hence, it follows trivially that
$$
(D(\iota+\ubarP))_{(A,\Phi,p)}
\left(T_{(A,\Phi,p)}\sN_{(A_0,\Phi_0,p_0)}^\kappa\right)
=
T_{(\iota+\ubarP)(A,\Phi,p)}\Fred_{\kappa,n}(\bE^+,\bE^-),
$$
which yields the claim.
\end{proof}

While the manifold $\Fred_{\kappa,n}(\bE^+,\bE^-)$ contains the rank-one
locus $\JJ_\kappa(\fX_{\sm})$, the latter space is not necessarily smooth
and so it is convenient at this point to shift our attention instead to the
smooth rank-one locus $\tJJ_\kappa(\fX_{\sm})$ in the flag manifold 
$\Flag_{\kappa,n}(\bE^+,\bE^-)$ covering $\Fred_{\kappa,n}(\bE^+,\bE^-)$.
Recall that the canonical projection
$$
\Flag_{\kappa,n}(\bE^+,\bE^-)
\to
\Fred_{\kappa,n}(\bE^+,\bE^-)
$$
has finite-dimensional fibers $\PP(\Ker B)\simeq \PP^{\kappa-1}$ over points
$$
B\in \Fred_{\kappa,n}(\bE^+,\bE^-). 
$$
Then the submersion
\eqref{eq:StabThickenedPeriodSubmersionStratum} then lifts to a smooth map
\begin{equation}
\iota+\tubarP:\sN_{(A_0,\Phi_0,p_0)}^\kappa
\to
\Flag_{\kappa,n}(\bE^+,\bE^-),
\label{eq:StabThickenedPeriodFlagStratum}
\end{equation}
though not necessarily to a submersion, given by 
$$
(\ff,A,\Phi,p) \mapsto (\iota+\tubarP)(\ff,A,\Phi,p), 
$$
where
$$
(\iota+\tubarP)(\ff,A,\Phi,p)
:=
([\Phi],\iota(\ff)+\ubarP(A,\Phi,p)) 
= 
([\Phi],\iota(\ff)+\Dirac_{A,p}).
$$
The proof of the next lemma is the key application of Proposition
\ref{prop:SmalePerturbMap}. 

\begin{lem}
\label{lem:ThickenedModNbhdNoRankOneIrred}
Continue the above notation. Then there is a first-category subset
$\sP_{\fc}\subset\sP$, depending on $(A_0,\Phi_0,p_0)$, such
that for all $p\in \sP-\sP_{\fc}$, the thickened moduli space
$\sM_{(A_0,\Phi_0,p_0)}|_p := \pi_{\sM;\sP}^{-1}(p)\cap\sM_{(A_0,\Phi_0,p_0)}$
contains no points $(A,\Phi,p)$ with $\Phi$ rank one.
\end{lem}

\begin{proof}
{}From its definition, a countable union of first-category subsets of $\sP$
is again a first-category subset, so it will suffice to consider a single
open neighborhood of a point $(\ff_1,A_1,\Phi_1,p_1)$ in
$\sN_{(A_0,\Phi_0,p_0)}^\kappa$, as the space
$\sN_{(A_0,\Phi_0,p_0)}^\kappa$ is paracompact and so we may repeat the
argument below for each element of a countable open cover.

We introduce a third stabilization, this time for the map
\eqref{eq:StabThickenedPeriodFlagStratum}, yielding a submersion onto
$\Flag_{\kappa,n}(\bE^+,\bE^-)$.  Let $\{\Phi_{1,\alpha}\}_{\alpha=1}^\kappa$
be an orthonormal basis for the kernel of $(\iota+\ubarP)(\ff_1,A_1,\Phi_1,p_1)
= \iota(\ff_1)+\Dirac_{A_1,p_1}$ and let
\begin{equation}
\label{eq:FamilyProjOntoStabDiracKernel}
\pi_{(\ff,A,\Phi,p)} :\bE^+ \to \Ker(\iota(\ff)+\Dirac_{A,p})
\end{equation}
be the smooth family of $L^2$-orthogonal projections from $\bE^+$ onto
$\Ker(\iota(\ff)+\Dirac_{A,p})$, parametrized by
the points $(\ff,A,\Phi,p)$ in  
$\sN_{(A_0,\Phi_0,p_0)}^\kappa$.
Let $\bz := (z_1,\dots,z_\kappa) \in \CC^\kappa$ and define
a smooth map
\begin{equation}
\iota+\hatubarP:
\CC^\kappa\times\sN_{(A_0,\Phi_0,p_0)}^\kappa
\to
\Flag_{\kappa,n}(\bE^+,\bE^-)
\label{eq:TwiceStabThickenedPeriodStratum}
\end{equation}
by the assignment
$(\bz,\ff,A,\Phi,p) \mapsto (\iota + \hatubarP)(\bz,\ff,A,\Phi,p)$,
where
\begin{align*}
(\iota + \hatubarP)(\bz,\ff,A,\Phi,p)
&:=
\iota(\ff) + \hatubarP(\bz,A,\Phi,p)
\\
&:=
([\Phi+\pi_{(\ff,A,\Phi,p)}
(\textstyle{\sum_{\alpha=1}^\kappa} z_\alpha\Phi_{1,\alpha})],
\iota(\ff) + \Dirac_{A,p}).
\end{align*}
(Again, we should really use an open neighborhood of the origin in
$\CC^\kappa$ rather than all of $\CC^\kappa$ in the definition of the map
\eqref{eq:TwiceStabThickenedPeriodStratum}, so the vector in $\bE^+$
is non-zero and defines an element of $\PP(\bE^+)$; 
however, the map is clearly well-defined for $\bz$ near zero
and that is all we shall need below.)
Note that when $(\bz,\ff) = (0,0)$, the term $(\iota +
\hatubarP)(\bz,\ff,A,\Phi,p)$ simplifies to
\begin{align}
(\iota + \hatubarP)(0,0,A,\Phi,p)
&=
([\Phi],\Dirac_{A,p})
\label{eq:SimplifyExtToOrigPeriodMap}
\\
&= \tubarP(A,\Phi,p)
\in
\PP(\Ker\Dirac_{A,p})\times \Fred_{\kappa,n}(\bE^+,\bE^-),
\notag
\end{align}
with $\PP(\Ker\Dirac_{A,p}) \subset \PP(\bE^+)$. 

\begin{claim}
\label{claim:TwiceStabThickenedPeriodStratumSubmersion}
The map $(\iota + \hatubarP)$ of
\eqref{eq:TwiceStabThickenedPeriodStratum} is a submersion at the point
$(0,\ff_1,A_1,\Phi_1,p_1)$.
\end{claim}

\begin{proof}
{}From \eqref{eq:TwiceStabThickenedPeriodStratum} we have
$(\iota + \hatubarP)(0,\ff_1,A_1,\Phi_1,p_1) = ([\Phi_1],\iota(\ff_1) + 
\Dirac_{A_1,p_1})$. The map 
\eqref{eq:StabThickenedPeriodSubmersionStratum}, given by
$(\ff,A,\Phi,p)\mapsto \iota(\ff) + \Dirac_{A,p}$, is a submersion at
$(\ff_1,A_1,\Phi_1,p_1)$. Hence, if $B_t$ is a smooth path in
$\Fred_{\kappa,n}(\bE^+,\bE^-)$ through $B_1 := \iota(\ff_1) + 
\Dirac_{A_1,p_1}$, we may choose a smooth path 
$(\ff_t,A_t,\Phi_t,p_t)$ in $\sN_{(A_0,\Phi_0,p_0)}^\kappa$, passing
through $(\ff_1,A_1,\Phi_1,p_1)$, which maps to the smooth path $B_t$ near
$t=1$ via \eqref{eq:StabThickenedPeriodSubmersionStratum}, so
$$
B_t = \iota(\ff_t) + \Dirac_{A_t,p_t}. 
$$
Using \eqref{eq:FamilyProjOntoStabDiracKernel}, define a smooth path of
projections onto $\Ker B_t\subset\bE^+$ by
$$
\pi_t := \pi_{(\ff_t,A_t,\Phi_t,p_t)}, \quad\text{$t$ near $1$.}
$$
Then, $\pi_1 = \id$ on $\Ker B_1$ and
$\{\pi_t(\Phi_{1,\alpha})\}_{\alpha=1}^\kappa$ is a basis for $\Ker B_t =
\Ker(\iota(\ff_t) + \Dirac_{A_t,p_t})$ near $t=1$, since
$\{\Phi_{1,\alpha})\}_{\alpha=1}^\kappa$ was chosen to be a basis for $\Ker
B_1$. Now suppose $\Psi_t$ is a smooth path in $\bE^+$ through $\Phi_1$
such that $B_t\Psi_t = 0$, so $([\Psi_t],B_t)$ is a smooth path in
$\Flag_{\kappa,n}(\bE^+,\bE^-)$.  Because $\Psi_t, \Phi_t \in \Ker B_t$,
near $t=1$ we can write
$$
\Psi_t - \Phi_t
= \sum_{\alpha=1}^\kappa(\Psi_t - \Phi_t,\pi_t(\Phi_{1,\alpha}))_{L^2}
\cdot\pi_t(\Phi_{1,\alpha})
=:
\sum_{\alpha=1}^\kappa z_\alpha(t)\pi_t(\Phi_{1,\alpha}),
$$
with $z_\alpha(t) = 0$ when $t=1$. Hence, near $t=1$, the path
$([\Psi_t],B_t)$ through $(\iota + \hatubarP)(0,\ff_1,A_1,\Phi_1,p_1)$ in
the image, $\Flag_{\kappa,n}(\bE^+,\bE^-)$, lifts to a smooth path
$(\bz_t,\ff_t,A_t,\Phi_t,p_t)$ in the domain,
$\CC^\kappa\times\sN_{(A_0,\Phi_0,p_0)}^\kappa$, through
$(0,\ff_1,A_1,\Phi_1,p_1)$ with
$$
(\iota + \hatubarP)(\bz_t,\ff_t,A_t,\Phi_t,p_t) 
= 
([\Psi_t],B_t).
$$
In particular, we have
$$
\frac{d}{dt}(\iota + \hatubarP)(\bz_t,\ff_t,A_t,\Phi_t,p_t)\Bigr|_{t=1}
=
\frac{d}{dt}([\Psi_t],B_t)\Bigr|_{t=1}.
$$
Since the smooth path $([\Psi_t],B_t)$ through $(\iota +
\hatubarP)(0,\ff_1,A_1,\Phi_1,p_1)$ was arbitrary, we see that 
$\iota + \hatubarP$ is a submersion at $(0,\ff_1,A_1,\Phi_1,p_1)$, as
desired. 
\end{proof}

By Claim \ref{claim:TwiceStabThickenedPeriodStratumSubmersion}, the map
\eqref{eq:TwiceStabThickenedPeriodStratum} is a submersion at the point
$(0,\ff_1,A_1,\Phi_1,p_1)$ and so is a submersion on an open 
neighborhood of this point in $\CC^\kappa\times\sN_{(A_0,\Phi_0,p_0)}^\kappa$. 
Rather than introduce further notation, we may
suppose without loss of generality 
(see the remark at the beginning of the proof of the lemma
concerning first-category sets) that the map
\eqref{eq:TwiceStabThickenedPeriodStratum} is a submersion on its entire domain
$\CC^\kappa\times\sN_{(A_0,\Phi_0,p_0)}^\kappa$.

Consequently, the locus 
$$
\sZ_{(A_0,\Phi_0,p_0)}^\kappa
:=
(\CC^\kappa\times\sN_{(A_0,\Phi_0,p_0)}^\kappa)
\cap(\iota+\hatubarP)^{-1}(\tJJ_\kappa(\fX_{\sm}))
$$
is a $C^\8$ Banach submanifold of
$\CC^\kappa\times \sN_{(A_0,\Phi_0,p_0)}^\kappa$. Let
\begin{align*}
\sM_{(A_0,\Phi_0,p_0)}^\kappa
:&=
\sN_{(A_0,\Phi_0,p_0)}^\kappa
\cap
(\{0\}\times\sM_{(A_0,\Phi_0,p_0)})
\\
&\subset
O_{(A_0,\Phi_0,p_0)}\times\sM_{(A_0,\Phi_0,p_0)},
\end{align*}
and observe that we have a countable disjoint union of subsets:
$$
\sM_{(A_0,\Phi_0,p_0)}
=
\bigcup_{\kappa\ge n} \sM_{(A_0,\Phi_0,p_0)}^\kappa.
$$
Of course, the unstabilized period map \eqref{eq:ThickenedPeriodMap} also
lifts to a map,
\begin{equation}
\tubarP:\sM_{(A_0,\Phi_0,p_0)}^\kappa
\to
\Flag_{\kappa,n}(\bE^+,\bE^-),
\label{eq:ThickenedPeriodFlagStratum}
\end{equation}
defined, as in \eqref{eq:StabThickenedPeriodFlagStratum}, by setting
$$
\tubarP(A,\Phi,p) := (\iota+\tubarP)(0,A,\Phi,p) = ([\Phi],\Dirac_{A,p})
= ([\Phi],\ubarP(A,\Phi,p)).
$$
(Note that the preimage $\ubarP^{-1}(\JJ_\kappa(\fX_{\sm}))$ in 
$\sM_{(A_0,\Phi_0,p_0)}^\kappa$ consists of points $(A,\Phi,p)$ such that 
$\Ker\Dirac_{A,p}$ contains {\em some\/} rank-one section, whereas 
the preimage $\tubarP^{-1}(\tJJ_\kappa(\fX_{\sm}))$ in 
$\sM_{(A_0,\Phi_0,p_0)}^\kappa$ consists of points $(A,\Phi,p)$ with $\Phi$
a rank-one element of $\Ker\Dirac_{A,p}$.)

Note that by
\eqref{eq:SimplifyExtToOrigPeriodMap} we have
\begin{equation}
\label{eq:ParamSetEquality}
\begin{aligned}
{}&\sM_{(A_0,\Phi_0,p_0)}^\kappa\cap\tubarP^{-1}(\tJJ_\kappa(\fX_{\sm}))
\\
&\quad =
(\{0\}\times\{0\}\times\sM_{(A_0,\Phi_0,p_0)}^\kappa)
\cap(\iota+\hatubarP)^{-1}(\tJJ_\kappa(\fX_{\sm}))
\\
&\quad =
(\{0\}\times\{0\}\times\sM_{(A_0,\Phi_0,p_0)}^\kappa)
\cap\sZ_{(A_0,\Phi_0,p_0)}^\kappa.
\end{aligned}
\end{equation}
Let $\pi_{C\times\sN;\sP}$ be the projection to the factor $\sP$ from the
$C^\8$ Banach submanifold
$$
\CC^\kappa\times \sN_{(A_0,\Phi_0,p_0)}^\kappa
\subset
\CC^\kappa \times F_{(A_0,\Phi_0,p_0)} \times \sC_{W,E}^{*,0}\times\sP,
$$
and similarly define $\pi_{\sN;\sP}:\sN_{(A_0,\Phi_0,p_0)}^\kappa\to \sP$.
Proposition \ref{prop:SmalePerturbMap} now implies that there is a
first-category subset $\sP_{\fc}\subset\sP$ such that for each
$p\in\sP-\sP_{\fc}$, the preimages
\begin{equation}
\label{eq:DefnRestriction1}
\CC^\kappa\times \sN_{(A_0,\Phi_0,p_0)}^\kappa|_p
:= 
\pi_{C\times\sN;\sP}^{-1}(p)
\cap(\CC^\kappa\times \sN_{(A_0,\Phi_0,p_0)}^\kappa),
\end{equation}
where $\sN_{(A_0,\Phi_0,p_0)}^\kappa|_p :=
\pi_{\sN;\sP}^{-1}(p)\cap\sN_{(A_0,\Phi_0,p_0)}^\kappa$, and 
\begin{equation}
\label{eq:DefnRestriction2}
\sZ_{(A_0,\Phi_0,p_0)}^\kappa|_p
:=
(\iota+\hatubarP(\cdot,p))^{-1}(\tJJ_\kappa(\fX_{\sm}))
\cap 
(\CC^\kappa\times \sN_{(A_0,\Phi_0,p_0)}^\kappa|_p)
\end{equation}
are smooth submanifolds; the parameters $p\in\sP-\sP_{\fc}$ are regular values
of the projections defining the above two preimages, which are $C^\8$ and
Fredholm.

\begin{claim}
\label{claim:EmptyZ}
The manifolds $\CC^\kappa\times \sN_{(A_0,\Phi_0,p_0)}^\kappa|_p$ and 
$\sZ_{(A_0,\Phi_0,p_0)}^\kappa|_p$ have the following dimension and
codimension, respectively:
\begin{align}
\begin{split}
&\dim_\RR\left(\CC^\kappa\times \sN_{(A_0,\Phi_0,p_0)}^\kappa\right)|_p
\\
&\quad =
\dim_\RR F_{(A_0,\Phi_0,p_0)} + 2\kappa + \dim
\bH_{(A_0,\Phi_0,p_0)}^1 - 2\kappa(\kappa-n) < \8,
\end{split}
\label{eq:FiniteFDimension}
\\
\begin{split}
&\codim_\RR\left(\sZ_{(A_0,\Phi_0,p_0)}^\kappa|_p;
\CC^\kappa\times\sN_{(A_0,\Phi_0,p_0)}^\kappa|_p\right)
\\
&\quad = \codim(\tJJ_\kappa(\fX_{\sm});\Flag_{\kappa,n}(\bE^+,\bE^-)) = \8.
\end{split}
\label{eq:InfiniteZCodimension}
\end{align}
In particular, $\sZ_{(A_0,\Phi_0,p_0)}^\kappa|_p$ is empty.
\end{claim}

\begin{proof}
According to \eqref{eq:DefnStratumOMSpace} and \eqref{eq:DefnRestriction1}
we have 
\begin{align*}
&\CC^\kappa\times\sN_{(A_0,\Phi_0,p_0)}^\kappa|_p
\\
&\quad =
(O_{(A_0,\Phi_0,p_0)} \times \sM_{(A_0,\Phi_0,p_0)}|_p)
\cap
(\iota+\ubarP(\cdot,p))^{-1}(\Fred_{\kappa,n}(\bE^+,\bE^-)).
\end{align*}
Hence, the expression for the dimension of $\CC^\kappa\times
\sN_{(A_0,\Phi_0,p_0)}^\kappa|_p$ follows from the facts that
$$
\codim_\RR\left(\Fred_{\kappa,n}(\bE^+,\bE^-);\Fred_n(\bE^+,\bE^-)\right)
=
2\kappa(\kappa-n),
$$
(just as in the derivation of \eqref{eq:StratumofOMSpaceCodim})
in conjunction with Proposition \ref{prop:SmalePerturbMap} 
applied to the submersion \eqref{eq:StabThickenedPeriodSubmersion} given by
$\iota+\ubarP$ on $O_{(A_0,\Phi_0,p_0)} \times
\sM_{(A_0,\Phi_0,p_0)}$ and the
genericity of $p$, the fact that $\dim_\RR
F_{(A_0,\Phi_0,p_0)} = \dim_\RR O_{(A_0,\Phi_0,p_0)}$, and the fiber
dimension formula implied by \eqref{eq:DefnH1}, namely
$$
\dim\sM_{(A_0,\Phi_0,p_0)}|_p
=
\dim\bH_{(A_0,\Phi_0,p_0)}^1.
$$  
The equality in \eqref{eq:InfiniteZCodimension} follows from the fact that
the map $\iota+\hatubarP(\cdot,p)$ is a submersion on
$\CC^\kappa\times\sN_{(A_0,\Phi_0,p_0)}^\kappa|_p$ according to Claim
\ref {claim:TwiceStabThickenedPeriodStratumSubmersion} and
Proposition \ref{prop:SmalePerturbMap}, for generic $p$.  
The infinite-codimension assertion
for $\sZ_{(A_0,\Phi_0,p_0)}^\kappa|_p$ follows from
\eqref{eq:CodimIncLocusInFlagOp}.  Thus, $\sZ_{(A_0,\Phi_0,p_0)}^\kappa|_p$
is a smooth manifold of negative dimension, by \eqref{eq:FiniteFDimension}
and \eqref{eq:InfiniteZCodimension}, and so is empty.
\end{proof}

{}From \eqref{eq:ParamSetEquality}, we see that (for any $p\in\sP$) we have
$$
\sM_{(A_0,\Phi_0,p_0)}^\kappa|_p
\cap\tubarP(\cdot,p)^{-1}(\tJJ_\kappa(\fX_{\sm}))
\subset
\sZ_{(A_0,\Phi_0,p_0)}^\kappa|_p,
$$
and so, as $\sZ_{(A_0,\Phi_0,p_0)}^\kappa|_p$ is empty
according to Claim \ref{claim:EmptyZ}, the intersection 
$\sM_{(A_0,\Phi_0,p_0)}^\kappa|_p\cap
\tubarP(\cdot,p)^{-1}(\tJJ_\kappa(\fX_{\sm}))$
must be empty as well. Therefore, 
the thickened moduli space $\sM_{(A_0,\Phi_0,p_0)}|_p$ contains
no points $(A,\Phi,p)$ with $\Phi$ rank one, when
$\dim_\CC\Ker\sD_{A,p}=\kappa$ and $p\in \sP_\fc$. 
Repeating this argument for all $\kappa\ge
n$ yields the desired conclusion,
recalling that a countable union of first-category subsets of $\sP$ is
again a first-category subset. This completes the proof of Lemma
\ref{lem:ThickenedModNbhdNoRankOneIrred}. 
\end{proof}

Since we have an inclusion
$$
\fM_{W,E}^{*,0}\cap \sM_{(A_0,\Phi_0,p_0)}
\subset \sM_{(A_0,\Phi_0,p_0)}
$$
of an open neighborhood of the point $(A_0,\Phi_0,p_0)$ in the parametrized
moduli space $\fM_{W,E}^{*,0}$ into the thickened, parametrized moduli
space $\sM_{(A_0,\Phi_0,p_0)}$, Lemma
\ref{lem:ThickenedModNbhdNoRankOneIrred} implies that 
$M_{W,E}^{*,0}(p)\cap \sM_{(A_0,\Phi_0,p_0)}|_p$
contains no points $(A,\Phi,p)$ with $\Phi$ rank one, when $p\in
\sP-\sP_{\fc}$. We now repeat this process for every point
$(A_1,\Phi_1,p_1)$ in $\fM_{W,E}^{*,0}$, using appropriate open neighborhoods
$\sU_{(A_1,\Phi_1,p_1)}$ of $(A_1,\Phi_1,p_1)$ in
$\sC_{W,E}^{*,0}\times\sP$, and obtaining a first-category subset for each
such neighborhood.  
Since $\sC_{W,E}^{*,0}\times\sP$ is a paracompact Banach
manifold, we may pass to a collection of open neighborhoods in
$\sC_{W,E}^{*,0}\times\sP$ which gives a countable open covering of
$\fM_{W,E}^{*,0}$. A countable union of first-category subsets of $\sP$ is
again a first-category subset, so this process yields the desired set
$\sP_{\fc}\subset \sP$ in the case of $C^r$ parameters.  That is, for each
$C^r$ parameter $p\in\sP-\sP_{\fc}$, the moduli space $M_{W,E}^{*,0}(p)$
contains no irreducible, rank-one $\PU(2)$ monopoles. The argument of
\cite[\S 5.1.2]{FL1} then allows us to reduce to the case of $C^\8$
parameters. This completes the proof of Theorem
\ref{thm:NoRankOneIrreducibles}. 
\end{proof}

\begin{rmk}
It is well-known that a Dirac operator $\Dirac_A:\Om^0(W^+\otimes F)\to
\Om^0(W^-\otimes F)$ defined by fixed \spinc connection on $W$ and unitary
connection $A$ on a Hermitian vector bundle $F$ over a closed four-manifold
has $\Coker \Dirac_A = 0$ when $\Ind \Dirac_A \geq 0$ and the connection $A$ is
generic; see, for example, \cite[Lemma 6.9.3]{MorganSWNotes} for a proof
using the Sard-Smale theorem when $F$ is a line bundle. This and related
vanishing results for harmonic spinors are proved directly in
\cite{Anghel}, without appealing to the Sard-Smale theorem.
\end{rmk}

\appendix
\section{PU(2) monopoles on K\"ahler surfaces and Spin${}^c$ polynomial
invariants} 
\subsection{Rank-one, irreducible PU(2) monopoles on K\"ahler surfaces}
\label{app:Teleman}
The observation that the statement and proof of Proposition I.3.5 in
\cite{PTDirac} is incorrect is an important one 
due to A. Teleman \cite{OTSurvey}, \cite{TelemanMonopole} for the theory of
\spinc 
polynomial invariants and for the developing theory of $\PU(2)$ monopoles.
The transversality results given in \cite[\S I]{PTDirac}, which underly the
Donaldson-Pidstrigach-Tyurin
theory of \spinc polynomial invariants, rely on \cite[Proposition
I.3.5]{PTDirac}. For completeness, we review Teleman's counterexample here. 

We first recall the form of the $\PU(2)$ monopole equations on a
K\"ahler manifold (see \cite{BradlowGP}, \cite{OTVortex},
\cite{OTQuaternion}, \cite{OTSurvey}, \cite{Witten}).  Let $(X,J,g)$
be a four-manifold with almost-complex structure $J$, Hermitian metric
$g$, and corresponding K\"ahler form $\omega$. The canonical \spinc
structure is defined by
$$
W^+_{\can} := \Lambda^{0,0}\oplus \Lambda^{0,2},
\qquad
W^-_{\can} := \Lambda^{0,1},
$$
where $\Lambda^{p,q} = \Lambda^{p,q}(T^*X)$, and $\rho:T^*X\to\End_\CC(W)$
is the standard Clifford multiplication of
\cite{Gilkey}, \cite{MorganSWNotes}, \cite{SalamonSWBook}.  
A straightforward modification of Witten's description of the solutions to
the $\U(1)$ monopole equations on a K\"ahler surface 
\cite{MorganSWNotes}, \cite{Witten} then
yields the following form of the (unperturbed) $\PU(2)$ monopole equations
\cite{OTVortex}, \cite{OTQuaternion} (see also \cite{BradlowGP}),
\begin{align}
F_A^{2,0} &= -\half(\alpha\otimes\barbeta)_0, \notag
\\
F_A^{0,2} &= \half(\beta\otimes\baralpha)_0, \label{eq:KaehlerMonopole}
\\
i\Lambda_g F_A 
&= -\half((\alpha\otimes\baralpha)_0 - *_g(\beta\otimes\barbeta)_0), \notag
\\
\barrd_A\alpha + \barrd_A^*\beta &= 0, \notag
\end{align}
for a pair $(A,\Phi)$, where $A$ is an $\SO(3)$ connection on $\su(E)$,
inducing a unitary connection on $E$ via the fixed unitary connection $A_e$
on $\det E$, and
$$
\Phi := (\alpha,\beta)
\in \Omega^0(W_{\can}^+\otimes E)
=
\Omega^{0,0}(E)\oplus \Omega^{0,2}(E).
$$
Now suppose $E$ is a Hermitian two-plane bundle with holomorphic structure
$\barrd_A$ over a K\"ahler surface $X$ with $H^0(E) \neq 0$ and that $\alpha$
is a holomorphic section of $E$. If the
holomorphic structure $\barrd_A$ on $E$ is indecomposable, then one finds
from \eqref{eq:KaehlerMonopole} that the triple $(A,\alpha,0)$
defines an irreducible, rank-one $\PU(2)$ monopole on $X$.

\subsection{The definition of spin${}^c$ polynomial invariants}
\label{app:SpincPolyInvariants}
The basic idea underlying the definition of \spinc polynomial
invariants is due to Donaldson \cite{DonArbeitstagung}. The
construction of \spinc polynomial invariants and their development and
application to smooth four-manifold topology has been carried out by
Pidstrigach and Tyurin in their series of articles \cite{PidVanDeVen},
\cite{PTDirac}, \cite{Tyurin}.  To the best of our knowledge, all
results concerning smooth four-manifold topology which have been
proved using \spinc polynomial invariants have been proved
independently using either Donaldson invariants or, more recently,
Seiberg-Witten invariants. Nonetheless, as the \spinc-ASD equations
can be viewed as precursor to the $\PU(2)$ monopole equations, the
\spinc polynomial invariants at least have some historical
interest. The
\spinc-ASD equations can be viewed (the original definition of
\cite{PTDirac} uses a slightly different trace condition for the unitary
connection on $\det E$) as the following variant of the $\PU(2)$ monopole
equations \eqref{eq:PT}:
\begin{align}
P_+(g)f(F_A) &= 0, \label{eq:SpincASD}
\\
(\Dirac_{A,f} + \rho(f(\vartheta)))\Phi &= 0. 
\notag
\end{align}
The moduli space 
$$
N_{W,E}(f,g,\vartheta) 
:= \{(A,\Phi): \textrm{Eq. \eqref{eq:SpincASD} holds}\}/\ssG_E
$$ 
of \spinc-ASD pairs is defined in the same way as the moduli space
$M_{W,E}(f,g,\tau,\vartheta)$ of $\PU(2)$ monopoles, with
$N_{W,E}^*(f,g,\vartheta)$ denoting the space of \spinc-ASD pairs
$(A,\Phi)$ where is $A$ irreducible. In \cite[\S I]{PTDirac} the moduli
space $N_{W,E}^*(f,g,\vartheta)$ is used to define invariants of smooth
four-manifolds.  As we remarked in \S \ref{app:Teleman}, Teleman has
pointed out that the proof given by Pidstrigach-Tyurin in \cite{PTDirac}
that $N_{W,E}^*(f,g,\vartheta)$ is a smooth manifold of the expected
dimension (see Proposition I.3.5 and Corollaries I.3.6. I.3.7, and I.3.8 in
\cite{PTDirac}), for generic $(g,\vartheta)$, is incorrect. We note
that variations of the Dirac operator $\Dirac_{A,f} + \rho(f(\vartheta))$
with respect to $f\in\Omega^0(\GL(T^*X))$ are {\em not\/} used in
\cite{PTDirac}, while the one-form $\vartheta$ is assumed to be purely
imaginary (a unitary perturbation of the $\U(1)$ connection on $\det W^+$),
rather than complex as we suppose here. Just as in the proof of our main
transversality result, Theorem \ref{thm:Transversality}, for the moduli
space of $\PU(2)$ monopoles, the principal difficulty one needs to address
is the possible presence in $N_{W,E}^*(f,g,\vartheta)$ of irreducible,
rank-one pairs. The proof of Theorem \ref{thm:NoRankOneIrreducibles}
carries over, with one slight change, to give: 

\begin{thm}
\label{thm:NoRankOneIrreducibleSpincASDPairs}
Let $X$ be a closed, oriented, simply-connected,smooth four-manifold with
$C^\8$ Riemannian metric $g$,
\spinc structure $(\rho,W^+,W^-)$, and Hermitian two-plane bundle $E$. Then
there is a first-category subset $\sP^\8_{\fc} \subset \sP^\8$ such that
for all $(f,\vartheta)$ in $\sP^\8 - \sP^\8_{\fc}$ the following holds:
The moduli space $N_{W,E}^{*,0}(f,g,\vartheta)$ contains no \spinc-ASD
pairs $(A,\Phi)$ with both $A$ irreducible and $\Phi$ rank one. 
\end{thm}

\begin{proof}
The only difference --- and this is the reason for the additional
constraint on $\pi_1(X)$ --- is that our unique
continuation result for reducible $\PU(2)$ monopoles, Theorem
\ref{thm:LocalToGlobalReducible}, must be replaced by the corresponding
unique continuation result for reducible 
anti-self-dual $\SO(3)$ connections \cite[Lemma 4.3.21]{DK}. The remainder
of the proof is otherwise identical to that of Theorem
\ref{thm:NoRankOneIrreducibles}. 
\end{proof}

Given Theorem \ref{thm:NoRankOneIrreducibleSpincASDPairs}, the remainder of
the argument of \cite[Proposition I.3.5]{PTDirac} yields:

\begin{thm}
\label{thm:TransversalitySpincASD}
Let $X$ be a closed, oriented, simply-connected,
smooth four-manifold with $C^\8$ Riemannian
metric $g$, \spinc structure $(\rho,W^+,W^-)$ with \spinc connection, and a
Hermitian line bundle $\det E$ with unitary connection. Then there is a
first-category subset $\sP^\8_{\fc}$ of the Fr\'echet space $\sP^\8$ of 
$C^\8$ perturbation parameters $(f,\vartheta)$ such that for all
$(f,\vartheta)$ in $\sP^\8 - \sP^\8_{\fc}$ the following holds: 
For each parameter $(f,\vartheta)$ in $\sP^\8 - \sP^\8_{\fc}$ 
and Hermitian two-plane bundle $E$ over $X$, the moduli space
$N^*_{W,E}(f,g,\vartheta)$ of \spinc-ASD pairs is a smooth
manifold of the expected dimension,
\begin{align*} 
\dim N^*_{W,E}(f,g,\vartheta)
&= 
-2p_1(\su(E))-\threehalf(e(X)+\sigma(X)) 
\\
&\quad + \half p_1(\su(E))+\half((c_1(W^+)+c_1(E))^2-\sigma(X))-1.
\end{align*}
\end{thm}

\begin{rmk}
\begin{enumerate}
\item The holonomy perturbations of \cite{FL1} can be used to give an
alternative, direct proof of Theorem \ref{thm:TransversalitySpincASD}, without
addressing the possible presence or absence of irreducible, rank-one
\spinc-ASD pairs in the moduli space $N^*_{W,E}(f,g,\vartheta)$.
\item The constraint on the topology of $X$ in the hypotheses of Theorem
\ref{thm:TransversalitySpincASD} can be relaxed by taking
account of the possible presence of the `twisted reducible'
anti-self-dual $\SO(3)$ connections of Kronheimer-Mrowka
\cite[\S 2]{KMStructure} with a little more care. 
\end{enumerate}
\end{rmk}

\subsection{Some linear algebra} 
\label{subsec:LinearAlgebra} 
When proving Theorem \ref{thm:Transversality}, our
transversality result for the $\PU(2)$ monopole equations
\eqref{eq:PT}, we used the following well-known, elementary fact:

\begin{lem}
\label{lem:MorganLemma}
\cite[p. 89]{MorganSWNotes}
Let $X$ be an oriented, Riemannian, smooth four-manifold with \spinc
structure $(\rho,W^+,W^-)$. Then the Clifford map $\rho:T^*X\to
\Hom_\CC(W^+,W^-)$ extends to an isomorphism of complex vector
bundles,
$$
\rho:T^*X\otimes_\RR\CC\to\Hom_\CC(W^+,W^-).
$$
\end{lem}

The applications arose, specifically, in the proofs of Lemma
\ref{lem:HermitianOrthog} and Proposition
\ref{prop:CokerDifferentialPeriod}. 
The lemma is also an essential ingredient in the standard proofs of
transversality for the perturbed Seiberg-Witten equations (see, for
example, \cite{KMThom}, \cite{MorganSWNotes}, \cite{Witten}). Hence,
though elementary, we include a proof here as we are not aware of a
reference.

\begin{proof}[Proof of Lemma \ref{lem:MorganLemma}]
Suppose $x\in X$. We need to show that the complex-linear map 
$$
\rho:(T^*X)_x\otimes_\RR\CC\to\Hom_\CC(W^+,W^-)_x
$$
is an isomorphism and, for this purpose, it is very convenient to use
the quaternion model for (complex) Clifford multiplication \cite{LM},
\cite{SalamonSWBook}. Thus, we employ the identifications of complex vector
spaces, 
$W^+|_x = \HH\oplus 0$ and $W^-|_x = 0\oplus\HH$. Moreover, with this
identification, $(T^*X)_x = \RR^4 = \HH$ acts on $W|_x = \HH\oplus\HH$ by
$$
\rho(v)(\phi^+,\phi^-) = (v\phi^- , -\bar{v}\phi^+),
\quad\text{for all }
v\in \HH \text{ and } (\phi^+,\phi^-) \in \HH\oplus\HH,
$$
while $i \in \CC$ acts by $-I \in \HH$ on the right, commuting with
Clifford multiplication on the left; we give $\HH$ its standard
basis $\{1,I,J,K\}$ over $\RR$.

Since $(T^*X)_x\otimes_\RR\CC$ and $\Hom_\CC(W^+,W^-)_x$ both have
complex dimension four, it suffices to show that $\rho$ is
surjective. Thus, by complex linearity,
it is enough to show that for a given unit-norm $\Phi
\in W^+|_x$ and $\Psi\in W^-|_x$, there is a $v \in
(T^*X)_x\otimes_\RR\CC$ such that $\rho(v) = \Psi\otimes
\langle\cdot,\Phi\rangle$ or, equivalently, that
\begin{equation}
\label{eq:QuaternionSurjectivity}
\rho(v)\Phi = \Psi.
\end{equation}
Suppose we are given $\Phi = (\phi,0) \in W^+|_x$ and $\Psi =
(0,\psi) \in W^-|_x$. Then,
$$
\rho(v)\Phi = v(\phi,0) = (0,-\bar{v}\phi),
\quad\text{for all }
v\in\HH.
$$
We want to show that $\rho(v)\Phi = \Psi$ for some $v\in \HH$, that is, 
\begin{equation}
\label{eq:QuaternionMultiplication}
(0,-\bar{v}\phi) = (0,\psi).
\end{equation}
But since $\HH$ acts transitively on itself by quaternionic
multiplication, we can find $v \in \HH$ so that
\eqref{eq:QuaternionMultiplication} --- and hence 
\eqref{eq:QuaternionSurjectivity} --- holds. This completes the proof.
\end{proof} 


\nocite{FrM}
\nocite{KazdanTopic}
\nocite{KadisonRingrose1}
\nocite{KotschickSW}
\nocite{MarcolliBook}
\nocite{Kuranishi}
\nocite{TauSymp}
\nocite{TelemanNonabelian}

\providecommand{\bysame}{\leavevmode\hbox to3em{\hrulefill}\thinspace}

\end{document}